\documentclass[a4paper,twoside,10pt]{article}
\usepackage[T1]{fontenc}
\usepackage{multirow,multicol}
\usepackage{stmaryrd}
\usepackage{enumerate}
\usepackage{graphicx,graphics,epstopdf}
\usepackage{amsmath,amssymb,subfigure}
\usepackage{amsthm}
\usepackage{color}
\usepackage {geometry}
\usepackage{ulem}
\usepackage{hyperref}
\hypersetup{
    colorlinks,%
    citecolor=blue,%
    filecolor=black,%
    linkcolor=red,%
    urlcolor=blue
}
\usepackage{natbib}

\newtheorem{theorem}{Theorem}[section]
\theoremstyle{definition}%remark
\newtheorem{remark}[theorem]{Remark}
\newtheorem{example}{Example}[section]

\numberwithin{equation}{section}

\usepackage{bbm}
\usepackage{dsfont}

\newcommand{\field}[1]{\mathds{#1}}
\def\etal{\mbox{et al.}}
\def\toprule{\\[-5.5pt]\hline\\[-3.5pt]}
\def\botrule{\\[-7.7pt]\hline}
\usepackage{xspace}
\newcommand{\pdf}{pdf}
\newcommand{\iid}{{\itshape i.i.d.}\xspace}
\renewcommand{\hat}{\widehat}

\DeclareMathOperator*{\supp}{supp}
\DeclareMathOperator*{\card}{Card}

\title{\textbf{\scshape On adaptive wavelet estimation \\ of a class of weighted densities}}
\author{{\scshape Fabien Navarro$^{1,2}$, Christophe Chesneau$^{1}$ and Jalal Fadili$^{2}$}\\\\
{\itshape $^{1,2}$LMNO CNRS-Universit\'e de Caen, 14032 Caen Cedex, France}\\
{\itshape $^{2}$GREYC CNRS-ENSICAEN-Universit\'e de Caen, 14050 Caen Cedex, France}\\}
\date{}

\usepackage{fancyhdr}
\begin{document}

\maketitle

\begin{abstract}
We investigate the estimation of a weighted density taking the form $g=w(F)f$, where $f$ denotes an unknown density, $F$  the associated distribution function and $w$ is a known (non-negative) weight. Such a class encompasses many examples, including those arising in order statistics or when $g$ is related to the maximum or the minimum of $N$ (random or fixed) independent and identically distributed (\iid) random variables. We here construct a new adaptive non-parametric estimator for $g$ based on a plug-in approach and the wavelets methodology. For a wide class of models, we show that it attains fast rates of convergence under the $\field{L}_p$ risk with $p\ge 1$ (not only for $p = 2$ corresponding to the mean integrated squared error) over Besov balls. Our estimator is also simple to implement and fast. We also report an extensive simulation study to support our findings.\vspace{0.5cm}

\noindent{\bf Key words and phrases:} Reliability, weighted density, density estimation, plug-in approach, wavelets, block thresholding, series system, parallel system.
\end{abstract}

%non-parametric functional estimation, order statistics, Lifetime distributions, Extreme value
%%%%%%%%%%%%%%%%%%%%%%%%%%%%%%%%%%%%%%%%%%%%%%%%%%%%%%%%%%%%%%%%%%%%%%%%%%%%%%%%%%%%%%%%%%%%%%%%%%%%%%%%%%%%%%%%%%%%%%%%%%%%
%%%%%%%%%%%%%%%%%%%%%%%%%%%%%%%%%%%%%%%%%%%%%%%%%%%%%%%%%%%%%%%%%%%%%%%%%%%%%%%%%%%%%%%%%%%%%%%%%%%%%%%%%%%%%%%%%%%%%%%%%%%%
\section{Introduction}\label{sec:intro}
%%%%%%%%%%%%%%%%%%%%%%%%%%%%%%%%%%%%%%%%%%%%%%%%%%%%%%%%%%%%%%%%%%%%%%%%%%%%%%%%%%%%%%%%%%%%%%%%%%%%%%%%%%%%%%%%%%%%%%%%%%%%
\subsection{Problem statement}
 Let $(\Omega, \mathcal{A},\field{P})$ be a probability space, $X$ be a real random variable with unknown density $f$ and $Y$ be a random variable having the unknown weighted density 
\begin{equation}\label{dens1}
g(x)=w(F(x))f(x), \qquad x\in \field{R},
\end{equation}
where $w$ denotes a known (of course non-negative) weight and $F$ denotes the distribution function of $f$. The goal we pursue here is to estimate $g$ from a random number $N$ of \iid sample $X_1,\ldots,X_N$ of $X$.

Such an estimation problem arises in many situations, typically when $g$ is related to the maximum\footnote{Since $\min (X_1,\dots,X_N)=-\max (-X_1,\dots,-X_N)$ the results can be easily reformulated for the sample minimum.} of $N$ \iid random variables, where $N$ is a discrete random number in $\field{N}^*$ which is independent of the $X_i$'s. Application fields cover hydrology, meteorology, reliability, investment, management science, insurance business, etc.. For example, when the $X_i$ are non-negative, the random variable $Y=\max(X_1,\dots,X_N)$ (or $Y=\min(X_1,\dots,X_N)$) arises naturally in reliability theory as the lifetime of a parallel (series) system with a random number $N$ of identical components with lifetimes $X_1,\dots,X_N$.

To make things clearer to the reader, we next give some illustrative examples.

\subsection{Motivating examples}
\begin{example}[Order statistics]\label{ex1} 
Let $X_1,\ldots,X_m$ be \textit{i.i.d.} random variables with absolutely continuous distribution function $F$ and probability density function (\pdf) $f$. Let $X_{(1)}\leq\ldots\leq X_{(m)}$ denote the corresponding order statistics.  Then, the \pdf~ $g_{X(j)}$ of the $j$-th order statistic is 
\begin{equation*}
g_{X(j)}(x)=\frac{m!}{(j-1)!(m-j)!}(F(x))^{j-1}(1-F(x))^{m-j}f(x), \qquad x\in \field{R}.
\end{equation*}
Thus, $X_{(m)}$, for example, is the random variable representing the largest observation of a sample of $n$ and corresponds to the sample maximum and the density $g_{X(m)}$ of $X_{(m)}=\max(X_1,\ldots, X_m)$ is given by
\begin{equation*}
g_{X(m)}(x)=m\left(F(x)\right)^{m-1}f(x), \qquad x\in \field{R}.
\end{equation*}
The aim is to estimate $g_{X(j)}$ from a $m$ \iid sample $X_1,\ldots,X_m$ of $X$. 
\end{example}

\begin{example}[Maximum of a random number $N$ of  \iid random variables]\label{ex2}
Let $X$ be a random variable with density $f$,  $\lbrace X_n\rbrace_{n\in\field{N}^*}$ be a sequence of \iid random variables with density $f$ and $N$ be a discrete random variable taking values in $\field{N}^*$ with a known probability mass function. Then the density of $Y=\max(X_1,\ldots, X_N)$  is 
\begin{equation}\label{densityy}
g(x)=w(F(x))f(x), \qquad x\in \field{R},
\end{equation}
where 
\[
w(u)=\sum_{k=1}^{\infty}ku^{k-1}\field{P}(N=k), \qquad u\in [0,1].
\]
The goal is again to estimate $g$ from an $n$ \iid sample $X_1,\ldots,X_n$ of $X$. 
\end{example}

\begin{example}[Pile-up model]\label{ex4}
 Let us now present the ``pile-up model''.  Let $\lbrace Y_n\rbrace_{n\in\field{N}^*}$ be a sequence of \iid random variables with density $g$, $N$ be a discrete random variable in $\field{N}^*$ as in the previous example, and let $X=\min(Y_1,\ldots, Y_N)$ with density $f$. Then the density of $Y_1$ is 
 \begin{equation*}
 g(x)=w(F(x))f(x), \qquad x \in \field{R},
 \end{equation*}
 where 
 \[
 w(u)=\frac{1}{M'(M^{-1}(1-u))}, \qquad u\in [0,1],
 \]
 $M(u)=\field{E}(u^N)$ and $M'(u)=\field{E}(Nu^{N-1})$. We are seeking an estimate of $g$ from a $n$ \iid sample $X_1,\ldots,X_n$ of $X$. 
\end{example}

\subsection{Previous work} 
Some distributional properties of the maximum and minimum of random variables have been extensively studied in the literature (see, e.g.,~\cite{ragpat},~\cite{shak} and~\cite{shawong}). In addition, the literature on order statistics contains a huge work about the maximum. In the context of extreme value theory, various statistical properties and (real data) applications can be found in~\cite{ada},~\cite{louz} and the references therein. 

The estimation of the density function of the maximum of two independent random variables has been considered by~\cite{chen} via kernel methods. The Pile-up model has been considered by~\cite{comtereb} via model selection methods. 

\subsection{Contributions and relation to prior work}
In this paper, we develop a new \textit{non-linear} adaptive estimator for $g$ in model \eqref{dens1} based on a plug-in method, wavelets and the block thresholding rule introduced by~\cite{cai}. Wavelet-based thresholding estimators are attractive for non-parametric function estimation because of their virtues from the viewpoints of spatial adaptivity, computational efficiency and asymptotic optimality properties.  In the case of simple density estimation, wavelet thresholding is probably one of the most attractive nonlinear methods. We refer to e.g.,~\cite{anto2},~\cite{hardle} and~\cite{vid} for a detailed discussion of the performances of wavelet estimators and some of their advantages over traditional methods such as kernel-based or projection estimators. 

We here explore the theoretical performance of our estimator under the $\mathbb{L}_p$ risk with $p\ge 1$ over a very rich class of function spaces, namely Besov spaces. Sharp rates of convergence are obtained. Application of our estimator to Example~\ref{ex2} above is described in detail. Finally, extensive simulation experiments are carried out to illustrate the practical performance of our estimator. In particular, the numerical tests indicate that our block thresholding estimator, which is simple to implement and fast, compares very favorably to standard kernel-based methods.  

\subsection{Paper organization}
The paper is structured as follows. Our wavelet estimator is described in Section~\ref{wav}. Section~\ref{results}  presents our estimator convergence rates. Simulations are detailed in Section~\ref{sec:simu}. %The proofs are postponed to Section~\ref{proofs}. 

%%%%%%%%%%%%%%%%%%%%%%%%%%%%%%%%%%%%%%%%%%%%%%%%%%%%%%%%%%%%%%%%%%%%%%%%%%%%%%%%%%%%%%%%%%%%%%%%%%%%%%%%%%%%%%%%%%%%%%%%%%%%
%%%%%%%%%%%%%%%%%%%%%%%%%%%%%%%%%%%%%%%%%%%%%%%%%%%%%%%%%%%%%%%%%%%%%%%%%%%%%%%%%%%%%%%%%%%%%%%%%%%%%%%%%%%%%%%%%%%%%%%%%%%%
\section{Wavelet estimators}\label{wav}
First of all, we briefly recall some key facts on wavelets and Besov spaces that will be essential to us in the sequel. Then we develop our nonlinear adaptive wavelet block thresholding estimator.

%%%%%%%%%%%%%%%%%%%%%%%%%%%%%%%%%%%%%%%%%%%%%%%%%%%%%%%%%%%%%%%%%%%%%%%%%%%%%%%%%%%%%%%%%%%%%%%%%%%%%%%%%%%%%%%%%%%%%%%%%%%%
\subsection{Wavelets and Besov balls}
Let $b>0$, $p>0$ and $\field{L}_p([-b,b])=\left\lbrace h : [-b,b]\rightarrow \field{R}; \ \|h\|_p^p=\int_{-b}^{b}|h(x)|^pdx < +\infty\right\rbrace$.

For the purposes of this paper, we use compactly supported wavelet bases on $[-b,b]$. More precisely, we consider the Daubechies family $\mathrm{db}_{2R}$ with the scaling and wavelet functions $\phi$ and $\psi$, where $R \ge 2$ is a fixed integer, see e.g.,~\cite{mallat}. Define the scaled and translated version of the $\phi$ and $\psi$
\[
\phi_{j,k}(x)=2^{j/2}\phi(2^jx-k),  \qquad \psi_{j,k}(x)=2^{j/2}\psi(2^jx-k).
\]
Then there exists an integer $\tau$ and a set of consecutive integers $\Lambda_j$ such that $\card(\Lambda_j)=C2^{j}$ for a $C>0$ and, for any integer $\ell\ge \tau$, the collection 
\[
\mathcal{B}=\{ \phi_{\ell,k}, \ k\in \Lambda_{\ell}; \ \psi_{j,k}; \ j \in \field{N}-\{0,\ldots, \ell-1\} ,\ k\in \Lambda_j\},
\]
is an orthonormal basis of $\field{L}_2([-b,b])$. 

Consequently, for any integer $ \ell \ge \tau$, any $h\in \field{L}_2([-b,b])$ can be expanded on $\mathcal{B}$ as
\begin{equation*}
h(x)= \sum_{k\in \Lambda_{\ell}}\alpha_{\ell,k}\phi_{\ell,k}(x)  +\sum_{j= \ell}^{\infty}  \sum_{k\in \Lambda_j}\beta_{j,k}\psi_{j,k}(x), 
\end{equation*}
where 
\begin{equation}\label{avion}
\alpha_{\ell,k}=\int_{-b}^{b} h(x)\phi_{\ell,k}(x)dx, \qquad  \beta_{j,k}=\int_{-b}^{b} h(x)\psi_{j,k}(x)dx.
\end{equation}

As is traditional in the wavelet estimation literature, we will investigate the performance of our estimator by assuming that the unknown density $f$ belongs to a Besov ball. The Besov norm for a function can be related to a sequence space norm on its wavelet coefficients. More precisely, let $M>0$, $s\in (0,R)$, $q \ge 1$ and $r \ge 1$.  A function $h \in \field{L}_p([-b,b])$ belongs to $B^s_{q,r}(M)$ if and only if there exists a constant $M^*>0$ (depending on $M$) such that  the associated wavelet coefficients \eqref{avion} satisfy
\begin{equation*}
\left(\sum_{k\in \Lambda_{\tau}}|\alpha_{\tau,k}|^{q}\right)^{1/q}+ \left(\sum_{j=\tau}^{\infty} \left(2^{j(s+1/2-1/q)}\left(\sum_{k\in \Lambda_j}|\beta_{j,k}|^{q}\right)^{1/q}\right)^{r}\right)^{1/r}\le  M^*,
\end{equation*}
with the usual modifications if $q=\infty$ or $r=\infty$.

In this expression, $s$ is a smoothness parameter and $q$ and $r$ are norm parameters.  They include many traditional smoothness spaces such as H\"older and Sobolev spaces. A comprehensive account on Besov spaces can be found in e.g.,~\cite{devore,meyer,hardle}.

%%%%%%%%%%%%%%%%%%%%%%%%%%%%%%%%%%%%%%%%%%%%%%%%%%%%%%%%%%%%%%%%%%%%%%%%%%%%%%%%%%%%%%%%%%%%%%%%%%%%%%%%%%%%%%%%%%%%%%%%%%%%
\subsection{Plug-in block wavelet estimator}\label{bart}
Let us consider the general statistical framework described in Section \ref{sec:intro} with \iid sample $X_1,\cdots,X_n$. First of all, we investigate the estimation of $f$ via the so-called wavelet block hard thresholding estimator. We suppose that $\supp(f)\subseteq [-b,b]$ with $b>0$.% and $h\in \field{L}_2([-b,b])$: since $f$ is a pdf, this assumption is obvious fulfilled.  

Let $p\ge 1$, and $j_1$ and $j_2$ be the integers corresponding to the finest and coarsest scales defined as 
\[
j_{1}=  \lfloor (\max(p,2)/2) \log_2 (\log n) \rfloor, \qquad j_2=\lfloor \log_2 (n/\log n) \rfloor,
\]
where $\lfloor a \rfloor$ denotes the whole number part of $a \in \field{R}^+$. For any $j\in\{j_1,\ldots,j_2\}$, let ${A}_j$ and $U_{j,K}$ be given such that $\cup_{K\in {A}_j}U _{j,K} =\Lambda_j$, $|U_{j,K}|=L=\lfloor (\log n)^{\max(p,2)/2} \rfloor$ and $U_{j,K}\cap U_{j,K'}=\emptyset$ for any $K \not = K'$ with $K$, $K' \in {A}_j$. In a nutshell, at each scale $j$, each $U_{j,K}$ is the set containing position indices of the wavelet coefficients inside block $K \in A_j$.\\

We define the wavelet block hard thresholding estimator of $f$ by
\begin{equation}\label{hinz2}
\hat f(x)=\sum_{k\in \Lambda_{j_1}}\hat \alpha_{j_1,k}\phi_{j_1,k}(x)+
\sum_{j=j_1}^{j_2} \sum_{K\in A_j }\sum_{k\in
U_{j,K}} \hat \beta_{j,k}{\mathds{1}}_{\left \lbrace  \left(\sum_{k\in
U_{j,K}}|\hat\beta_{j,k}|^p/L\right)^{1/p}\ge
{\kappa} n^{-1/2} \right\rbrace}  \psi_{j,k}(x), \qquad x\in \lbrack -b,b \rbrack ,
\end{equation}
where $\mathds{1}_{\lbrace\cdot\rbrace}$ denotes the indicator function, and
\begin{equation*}
\hat \alpha_{j_1,k}= \frac{1}{n}\sum_{i=1}^n \phi_{j_1,k}(X_i), \qquad \hat \beta_{j,k}= \frac{1}{n}\sum_{i=1}^n \psi_{j,k}(X_i)
\end{equation*}
%and $\kappa$ denotes a large enough constant. 
and $\kappa>0$ is a threshold parameter to be discussed later.

The estimator $\hat f$ was initially developed by~\cite{cai} for the regression model under the $\mathbb{L}_2$-risk  with equispaced deterministic samples.  The $\mathbb{L}_p$ risk version was studied in~\cite{picard} for the standard density estimation problem and by~\cite{chesneau2} for the biased density estimation problem. 

The idea underlying $\hat f$ \eqref{hinz2} is to operate a group/block selection: it keeps intact the large groups of unknown wavelet coefficients of $f$ \eqref{avion} and removes the others. Wavelet block thresholding is one of the most attractive non-linear thresholding methods, since it is both numerically straightforward to implement and asymptotically optimal for a large variety of Sobolev or Besov classes. Detailed references on the subject for various models include, but are not limited to,~\cite{cai,cai2},~\cite{li,li2},~\cite{picard} and~\cite{chesneau,chesneau2}.

The performance of Block thresholding estimators depends on the threshold level $\kappa$. In the non-parametric regression setting, in order to choose this key parameter,~\cite{CaiZhou} proposed an adaptive James-Stein block thresholding estimator whose parameters (including the threshold) minimize the Stein's unbiased risk estimate (SURE) and established its minimax rates of convergence under the mean squared error over Besov balls. Other selection strategies have been developed in the literature (see e.g.~\cite{Nason} which considered wavelet estimators based on cross-validation to choose the thresholding parameter in practice). In this work, we focus on the universal threshold proposed by~\cite{donoho}. The reason for this choice is twofold. First, it is the one consistent with the theoretical convergence rates established in Section~\ref{results}. Secondly, it allows to remain fair when comparing to the other methods of the literature tested in Section~\ref{sec:simu}.
%Note that these improvements can be easily implemented.

Finally, plugging \eqref{hinz2} into \eqref{dens1} leads to the following estimator of $g$:
\begin{equation}\label{plug2}
\hat g(x)=w(\hat F(x))\hat f (x), \qquad x\in [-b,b],
\end{equation} 
where 
\begin{equation}\label{F}
\hat F(x)=\frac{1}{n}\sum_{i=1}^n\mathds{1}_{\lbrace X_i\le x\rbrace}.
\end{equation}

The rest of the paper explores  the theoretical and practical performances of $\hat g$. 

%%%%%%%%%%%%%%%%%%%%%%%%%%%%%%%%%%%%%%%%%%%%%%%%%%%%%%%%%%%%%%%%%%%%%%%%%%%%%%%%%%%%%%%%%%%%%%%%%%%%%%%%%%%%%%%%%%%%%%%%%%%%
%%%%%%%%%%%%%%%%%%%%%%%%%%%%%%%%%%%%%%%%%%%%%%%%%%%%%%%%%%%%%%%%%%%%%%%%%%%%%%%%%%%%%%%%%%%%%%%%%%%%%%%%%%%%%%%%%%%%%%%%%%%%
\section{Estimator convergence rates}\label{results}
In this section, we discuss the asymptotic properties of the proposed estimator. Rates of $\mathbb{L}_p$ convergence are investigated under the following assumptions:
\begin{enumerate}[(\bf{A.}1)]
\item Compact support: $\supp(f)\subseteq [-b,b]$ with $b>0$.
\label{hyp:supp}
\item Uniform boundedness: there exists a constant $C_1>0$ such that
\label{hyp:linfty}
\begin{equation}\label{boundff}
\sup_{x\in [-b,b]}f(x)\le C_1 ~.
\end{equation}
\item Uniform Lipschitz continuity of $w$ (with Lipschitz constant $C_2$): %there exists a constant $C_2>0$ such that
\label{hyp:lip}
\begin{equation}\label{boundw}
|w(u)-w(v)|\le C_2 |u-v|, \qquad \mbox{for all}~ (u,v)\in [0,1]^2,
\end{equation}
\end{enumerate}

Assumption A.\ref{hyp:supp} is a usual one in the wavelet density estimation framework (see e.g.~\cite{donoho}).  Extension to non-compactly supported densities might be possible and ideas from~\cite{Juditsky} might be inspiring, although these authors considered a model different from ours, and their results were valid only for the case where $f$ is in the H\"older class. Such an extension is however beyond the scope of this paper and we leave it for a future work.

\medskip
%\cite{Juditsky} proposed an adaptive wavelet procedure for density estimation on $\field{R}$ when $f$ belongs to a H\"{o}lder spaces. They propose a data-driven wavelet thresholding estimator and established its minimax rates of convergence (up to a logarithmic factor) under the $\field{L}_p$-risk. A possible future extension of our estimator to density function on the whole real line need further investigations that we leave for a future work.

%%%%%%%%%%%%%%%%%%%%%%%%%%%%%%%%%%%%%%%%%%%%%%%%%%%%%%%%%%%%%%%%%%%%%%%%%%%%%%%%%%%%%%%%%%%%%%%%%%%%%%%%%%%%%%%%%%%%%%%%%%%%
%\subsection{Estimator convergence rates}

Theorem \ref{danzig2} studies the $\mathbb{L}_p$ risk of $\hat g$ \eqref{plug2} over Besov balls and Assumptions~A.\ref{hyp:supp}-A.\ref{hyp:lip} on $f$ and $w$. 

\begin{theorem}\label{danzig2}
Consider the general statistical framework described in Section \ref{sec:intro} (the estimation of $g$ \eqref{dens1} is of interest). Suppose that Assumptions~A.\ref{hyp:supp}-A.\ref{hyp:lip} hold. Let $p\ge 1$ and $\hat g$ be given by \eqref{plug2}. Then, for any $f\in B^{s}_{q,r}(M)$, $q\ge 1$, $r\ge 1$ and $s\in (1/q,R)$, there exists a constant $C>0$ such that 
\begin{equation*}
\field{E}\left(\Vert \hat g-g\Vert_{p}^{p}\right)\le C\varphi_{n,p},
\end{equation*}
where 
\begin{equation}\label{rate}
\varphi_{n,p} =
\begin{cases}  
n^{-sp/(2s+1)}, & {\text{if}} \ \    q\ge p,\\
\left( \displaystyle \frac{\log n}{n}\right)^{sp/(2s+1)}, & {\text{if}} \ \    \{p>q, \ q s>(p-q)/2\},\\
\left( \displaystyle \frac{\log n}{n}\right)^{(s-1/q+1/p)p/(2(s-1/q)+1)}, & {\text{if}} \ \   q s<(p-q)/2 \ \   {\text{or}} \  \ \{q s=(p-q)/2, \  p\le q/r\}, \\ 
\left( \displaystyle \frac{\log n}{n}\right)^{(s-1/q+1/p)p/(2(s-1/q)+1)}(\log n)^{p-q/r}, & {\text{if}} \  \  \{q s=(p-q)/2, \  p> q/r\}.
\end{cases}
\end{equation}
\end{theorem}
Note that the rate $\varphi_{n,p}$ in \eqref{rate} is the near optimal (or optimal in some regimes) one in the minimax sense for $f$. See, e.g.,~\cite{donoho} and~\cite{hardle}. In plain words, the near optimality of the estimator of $f$ is transferred to that of $g$ through the plug-in principle.\\ 

The proof of Theorem \ref{danzig2} uses a suitable decomposition of the $\mathbb{L}_p$ risk and capitalizes on  results on the performances of $\hat f$ \eqref{hinz2} and $\hat F$ \eqref{F} established in~\cite{chesneau2}. %See Section~\ref{proofs} for details.

\begin{proof}[\textbf{Proof of Theorem \ref{danzig2}}]\label{proof1}
%\noindent {\bf Proof of Theorem \ref{danzig2}.} 
Observe that
\begin{equation*}
\hat g(x)-g(x) = w(\hat F(x))(\hat f(x)- f(x))+ f (x)(w(\hat F(x))- w(F(x))).
\end{equation*}		
By Assumption~A.\ref{hyp:lip} implying $\sup_{x\in  [0,1]}w(x)\le C$, together with Assumptions~\ref{hyp:linfty}, we have
\begin{align*}
|\hat g(x)-g(x) | & \le  C (  |\hat f(x)- f(x)|+|w(\hat F(x))- w(F(x))|)\\
& \le  C ( |\hat f(x)- f(x)|+|\hat F(x)- F(x)|).
\end{align*}	

By the Jensen inequality, we have
\begin{equation*}
\field{E}(\Vert \hat g-g\Vert_p^p)\le C(\field{E}(\Vert \hat f-f\Vert_p^p)+\field{E}(\Vert \hat F-F\Vert_p^p)).
\end{equation*}

It follows from \cite[Theorem 4.1 with $w(x)=1=\mu$]{chesneau2} that
\begin{equation*}
\field{E}(\Vert \hat f-f\Vert_p^p)\le C \varphi_{n,p},
\end{equation*}
where $\varphi_{n,p}$ is given by \eqref{rate}.

Now note that
\[
\hat F(x)-F(x)=\frac{1}{n}\sum_{i=1}^n U_i(x),
\]
with $U_i(x)=\mathds{1}_{\lbrace X_i\le x\rbrace}-F(x)$. Since $U_1(x),\ldots,U_n(x)$ are \iid with $\mathbb{E}(U_1(x))=0$, $|U_1(x)|\le 2$ and $\mathbb{E}( (U_1(x))^2)\le 4$, the Rosenthal inequality (see~\cite{rosen}) yields
\begin{equation*}
\field{E}(\Vert \hat F-F\Vert_p^p) \le  C\sup_{x\in \mathbb{R}} \mathbb{E}\left( (\hat F(x)-F(x))^{p} \right)\le  C \frac{1}{n^{p/2}}\le C \varphi_{n,p}.
\end{equation*}
Combining the inequalities above, we obtain the desired result i.e., 
\begin{equation*}
\field{E}(\Vert \hat g-g\Vert_p^p)\le C\varphi_{n,p}.
\end{equation*}
\end{proof}

%%%%%%%%%%%%%%%%%%%%%%%%%%%%%%%%%%%%%%%%%%%%%%%%%%%%%%%%%%%%%%%%%%%%%%%%%%%%%%%%%%%%%%%%%%%%%%%%%%%%%%%%%%%%%%%%%%%%%%%%%%%%
\subsection{An illustrative application}\label{applications1}
Let's recall Example~\ref{ex2}, where $\lbrace X_n\rbrace_{n\in\field{N}^*}$ is a sequence of \iid random variables with \pdf~ $f$ and $N$ be a discrete random variable of values in $\field{N}^*$ independent of this sequence. The density  of $Y=\max(X_1,\ldots, X_N)$ is given by \eqref{densityy}.

Suppose that Assumptions~A.\ref{hyp:supp}-A.\ref{hyp:linfty} hold.  Thus, several examples for the distribution of $N$ can be considered.% , which is entirely determined by the distribution of $N$.
\begin{enumerate}[(a)]
\item Degenerate distribution. $\field{P}(N=m)=1$. Then % Constant random variable
\begin{equation}\label{fixed}
w(u)=mu^{m-1}, \qquad u\in [0,1],
\end{equation}
\item Geometric distribution. $N\sim G(\eta)$ ($\field{P}(N=k)=\eta(1-\eta)^{k-1}$, $k\in\field{N}^*$). Then
\begin{equation}\label{geo}
w(u)=\frac{\eta}{(1-u(1-\eta))^2}, \qquad u\in [0,1],
\end{equation}
\item Poisson plus $1$ distribution. $N=P+1$ with $P\sim \mathcal{P}(\lambda)$  ($\field{P}(N=k)=e^{-\lambda} \frac{\lambda^{k-1}}{(k-1)!}$, $k\in\field{N}^*$).
Then 
\[
w(u)=e^{-\lambda}e^{\lambda u}(1+\lambda u), \qquad u\in [0,1].
\]
\end{enumerate}
\begin{remark}
In examples (a)--(c) above it is clear that Assumption~\ref{hyp:lip} is satisfied; more precisely, in example (a), we have $C_2=m(m-1)$; in example (b), we have $C_2=2(1-\eta)/\eta^2$; in example (c), we have $C_2=\lambda(2+\lambda)$.
\end{remark}

In this context, Theorem \ref{danzig2} can be applied. Let $p\ge 1$ and $\hat g$ be the estimator given in \eqref{plug2}. Then, for any $f\in B^{s}_{q,r}(M)$, $q\ge 1$, $r\ge 1$ and $s\in (1/q,R)$ there exists a constant $C>0$ such that 
\begin{equation*}
\field{E}\left(\Vert \hat g-g\Vert_{p}^{p}\right)\le C\varphi_{n,p},
\end{equation*}
where $\varphi_{n,p}$ is given by \eqref{rate}.

\begin{remark}
Taking $m=2$, the obtained rate is similar to the one attained by the kernel estimators developed by~\cite{chen}; the only difference is the extra-logarithmic term $(\log n)^{2s/(2s +1)}$. However, unlike kernel estimators~\cite{chen}, our procedure $\hat g$ is adaptive and our rate of convergence holds for a wider class of functions $f$ including H\"older class, Sobolev class, etc..
\end{remark}

%%%%%%%%%%%%%%%%%%%%%%%%%%%%%%%%%%%%%%%%%%%%%%%%%%%%%%%%%%%%%%%%%%%%%%%%%%%%%%%%%%%%%%%%%%%%%%%%%%%%%%%%%%%%%%%%%%%%%%%%%%%%
\section{Simulation results}\label{sec:simu}

We now illustrate these theoretical results by a simulation study within the context described in Section \ref{applications1}. That is, we consider the problem of estimating the density $g$ of the maximum of a random number $N$ of \iid random variables. From a reliability study standpoint, this problem corresponds to a parallel system with $N$ identical components. Thereby, we have considered two numerical examples. They complement the asymptotic results of \ref{danzig2}.

\begin{figure}[t]
\centering
\includegraphics[width=0.245\textwidth]{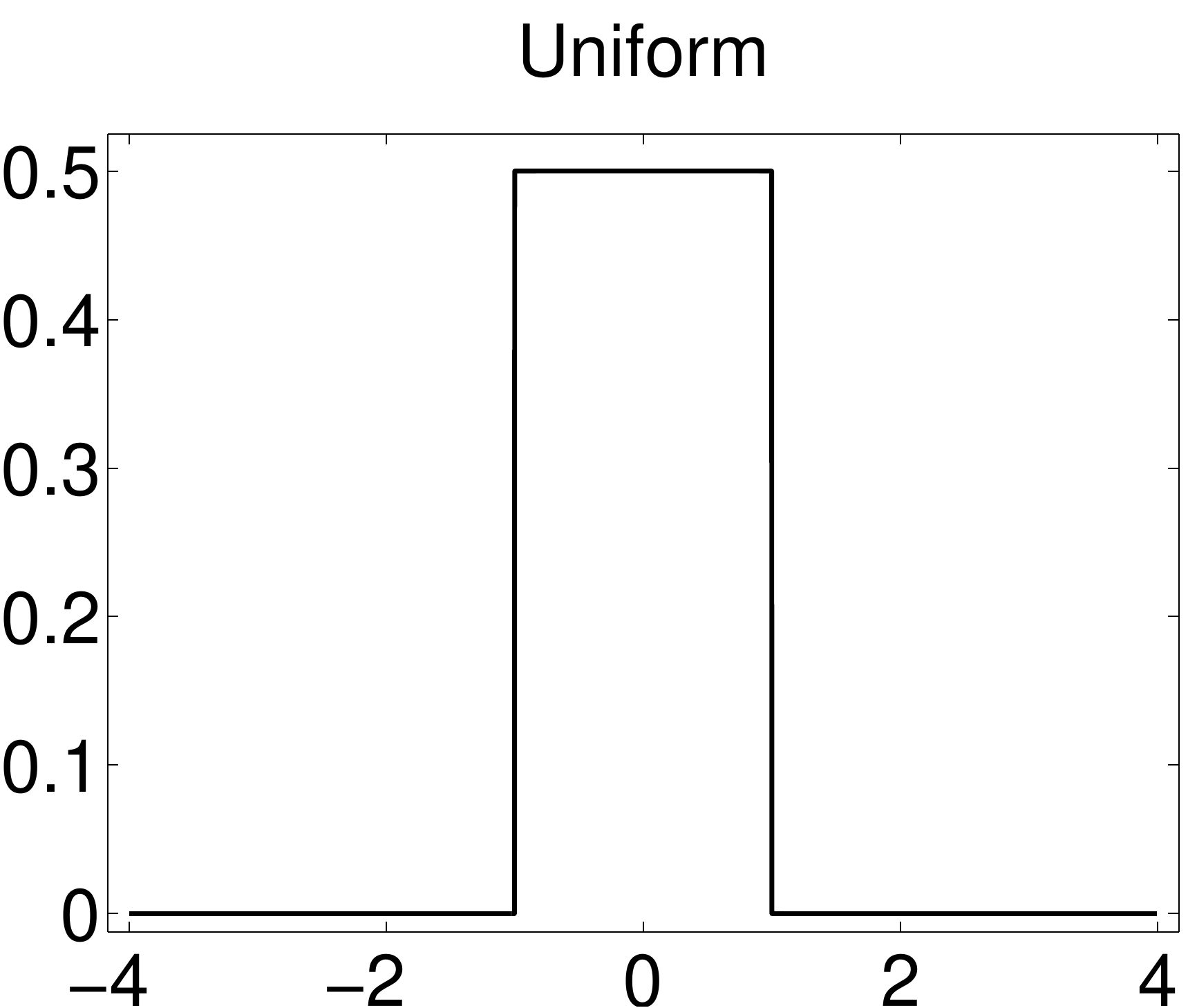}
\includegraphics[width=0.245\textwidth]{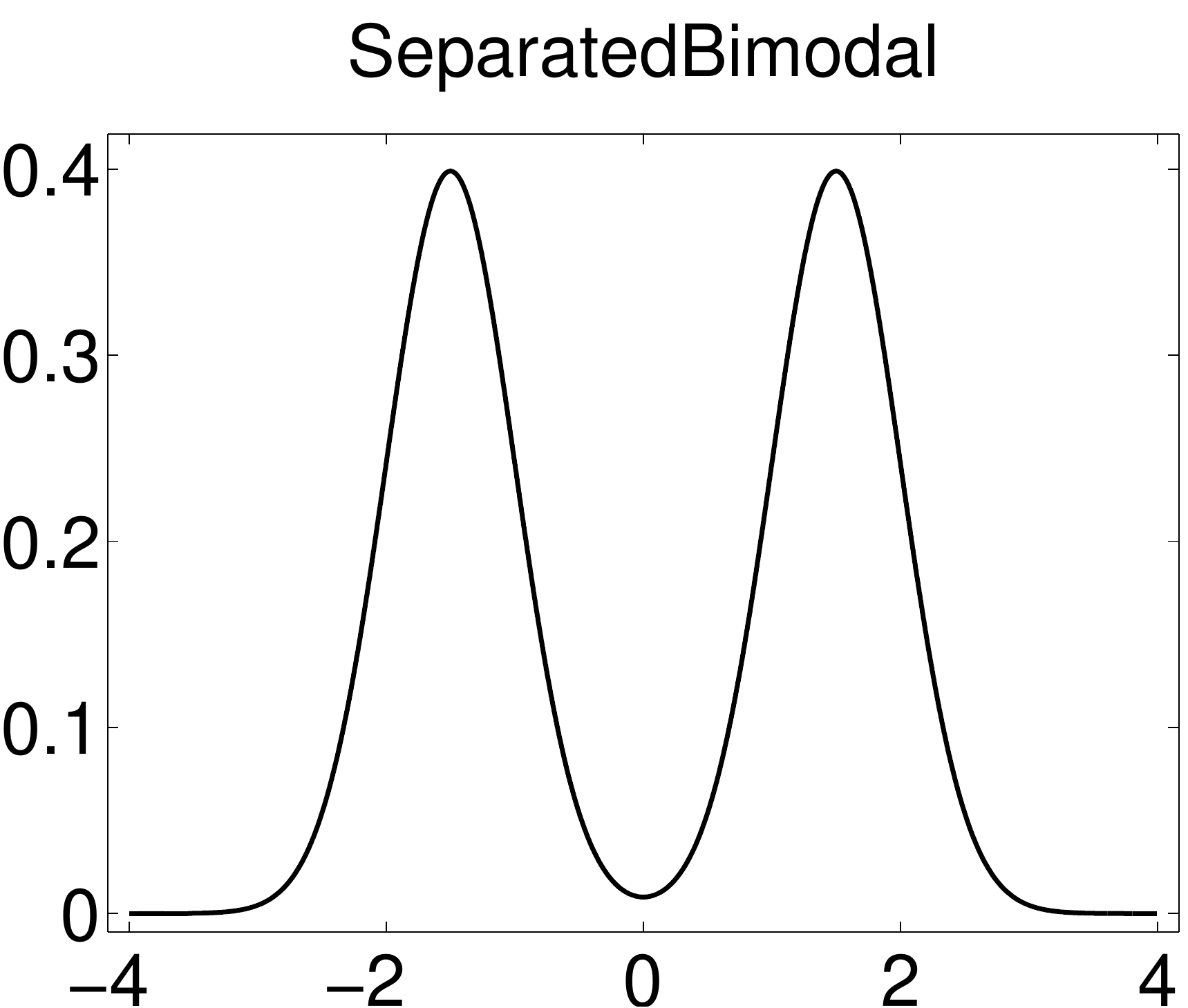}
\includegraphics[width=0.245\textwidth]{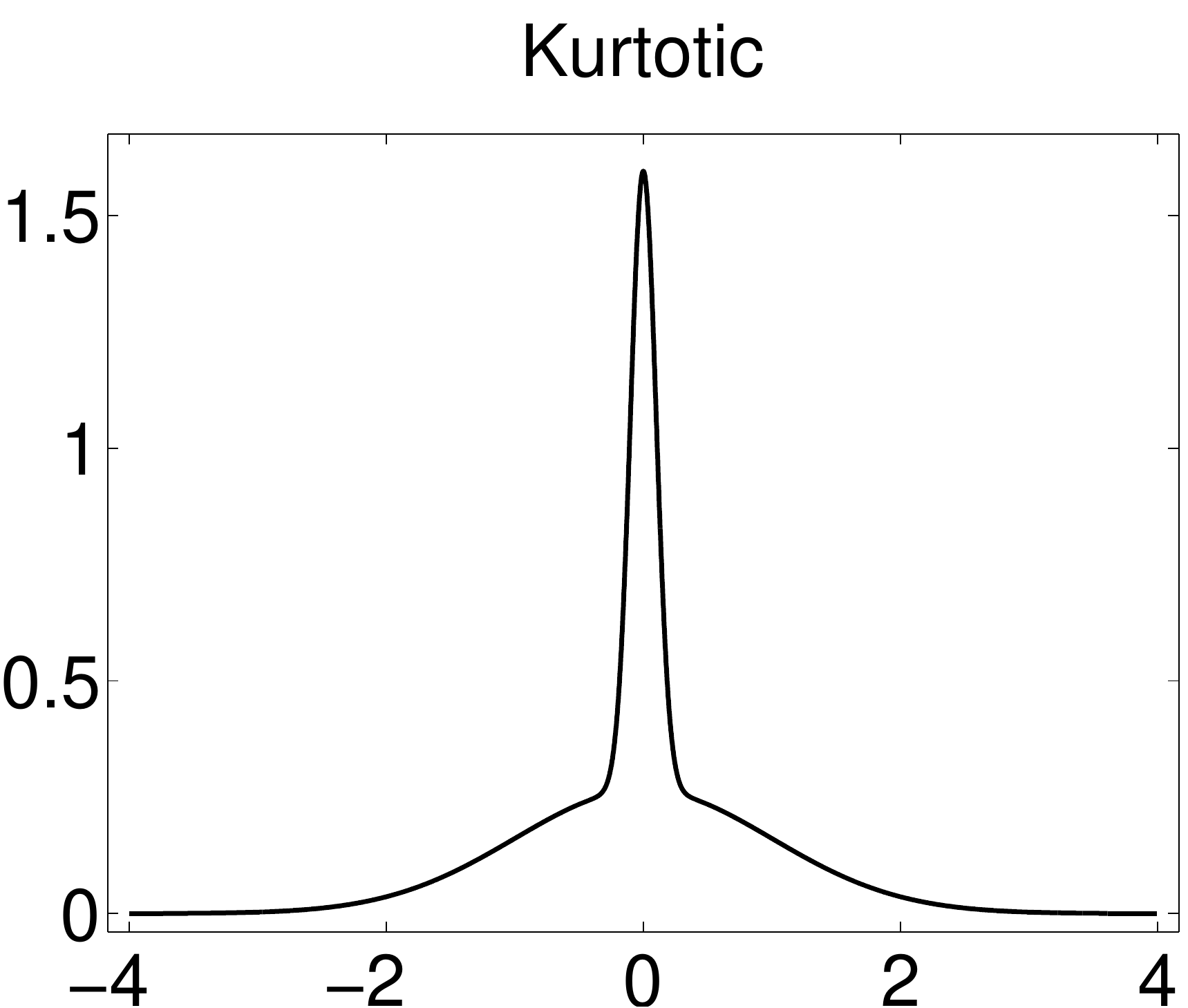}
\includegraphics[width=0.245\textwidth]{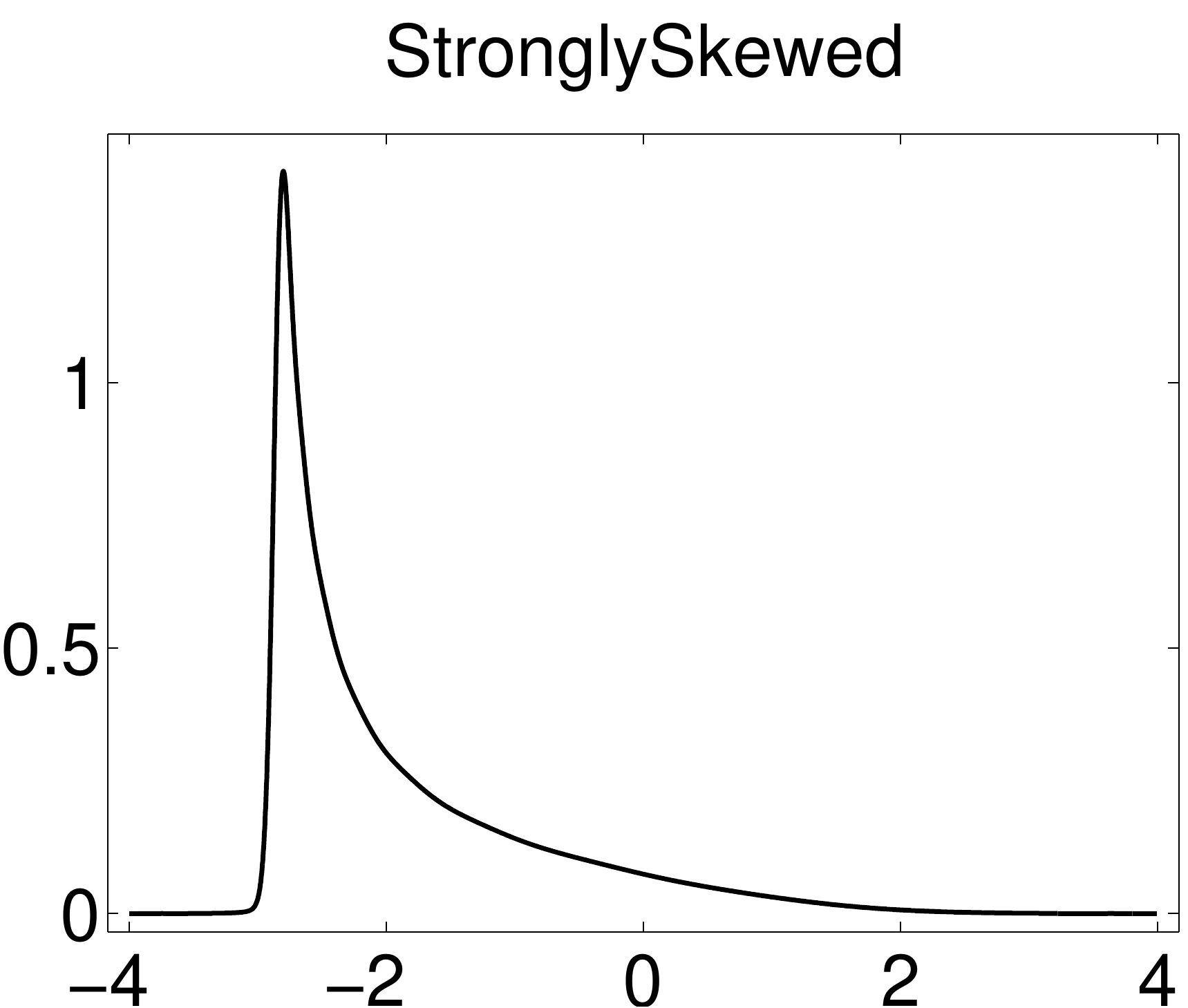}
\caption{Test densities.}
\label{fig:orig}
\end{figure}

\paragraph{Computational aspects.}
In the sequel, we will refer to our adaptive wavelet estimator \eqref{plug2} simply as Block. We have compared its performance to alternatives from the literature on several densities. We have considered the uniform distribution, as well as a family of normal mixture densities (``SeparatedBimodal'',  ``Kurtotic'' and ``StronglySkewed'', initially introduced in~\cite{MarronWand}) representing different degrees of smoothness (see Figure~\ref{fig:orig}). We assumed that the density function $f$ of the $X_i$'s has a compact support included in $[-b, b]$. We have used the formulae given by~\cite{MarronWand} to simulate such densities so that
\begin{equation*}
 \min_{l}(\mu_l - 3\sigma_l)=-3,\quad \max_{l}(\mu_l+3\sigma_l)=3,
\end{equation*}
where $l=1,\dots,d$ with $d$ the number of densities in the mixture (see, \cite[Section 4, Table 1]{MarronWand}, for the values of the parameters). Thereby, it is very unlikely to have values outside the interval $[-4,4]$ and we loose little by assuming compact support . The same kind of assumption was made in the context of wavelet density estimation by~\cite{vanvid}. In order to simplify the presentation of the results, one can simply rescale the data such that they fall into $[-b, b]$ (which covered the full range of all observed data). Thus, the density was evaluated at $T=2^J$ equispaced points $t_i=2ib/T$, $i=-T/2,\ldots,T/2-1$ between $-b$ and $b$ , where $J$ is the index of the highest resolution level and $T$ is the number of discretization points. The primary level $j_1=3$, $T=512$ and the Symmlet wavelet with $6$ vanishing moments were used throughout all experiments. All simulations have been implemented under Matlab. 

\begin{figure}
\centering
\subfigure[]{\includegraphics[width=0.245\textwidth]{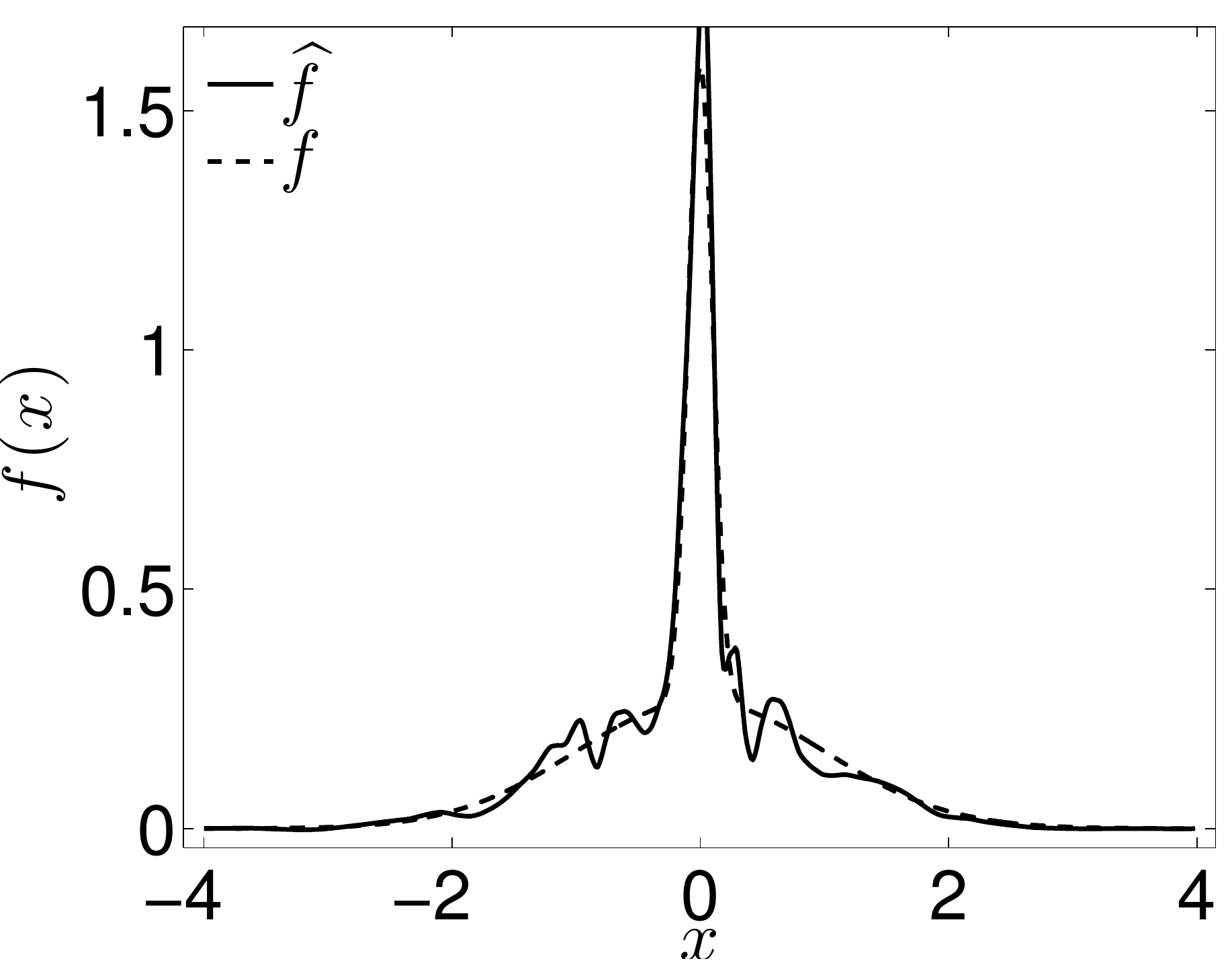}}
\subfigure[]{\includegraphics[width=0.245\textwidth]{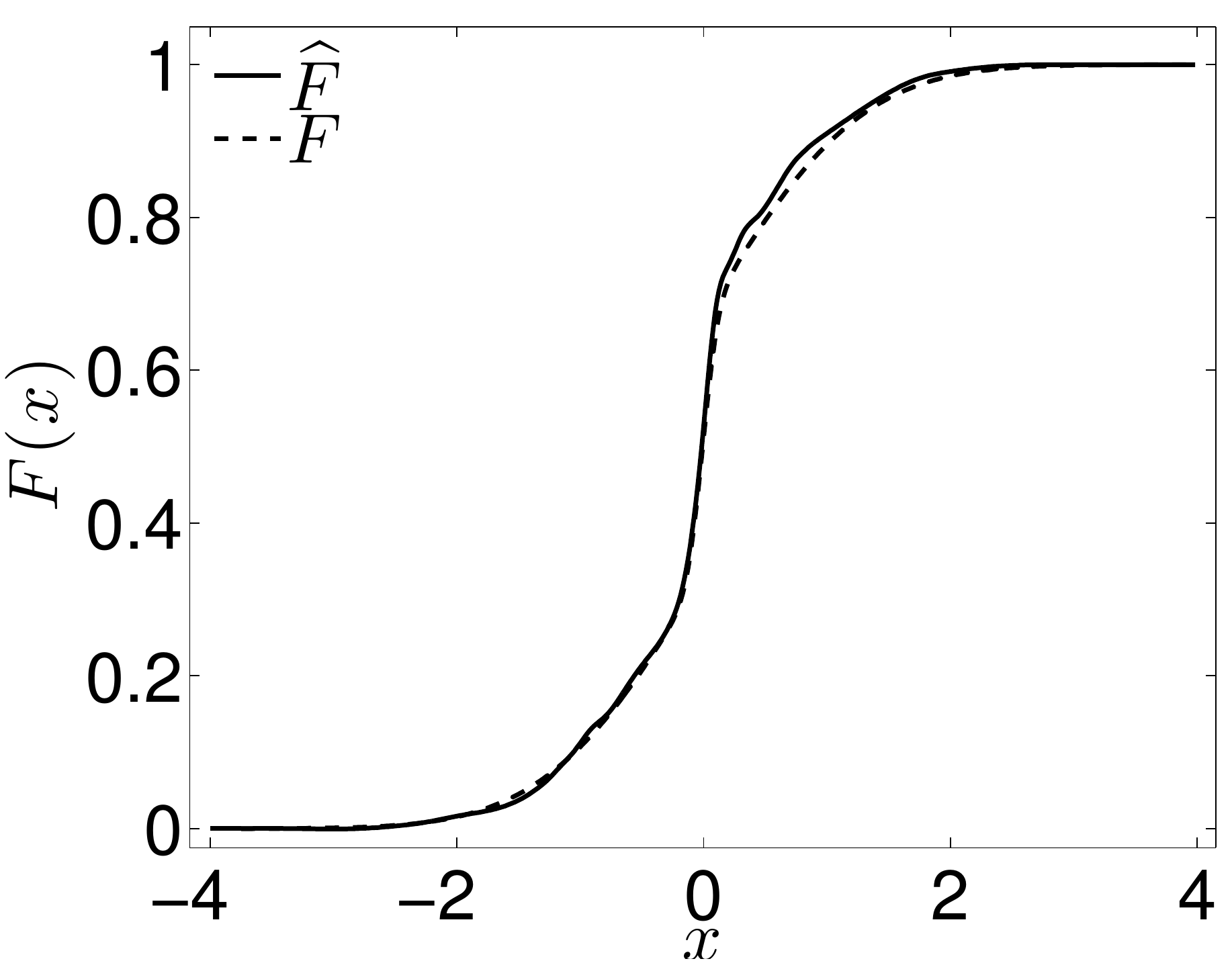}}
\subfigure[]{\includegraphics[width=0.235\textwidth]{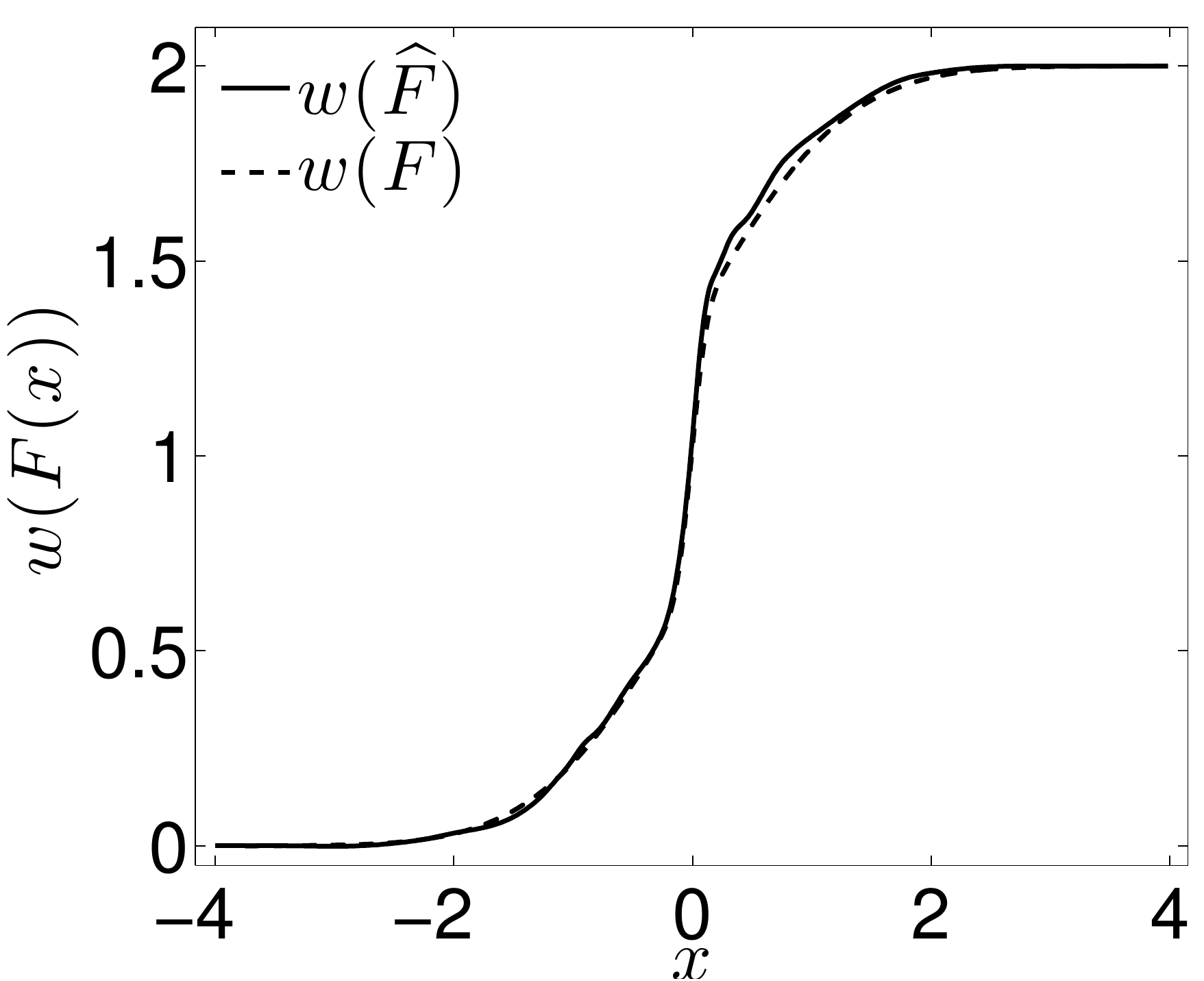}}
\subfigure[]{\includegraphics[width=0.245\textwidth]{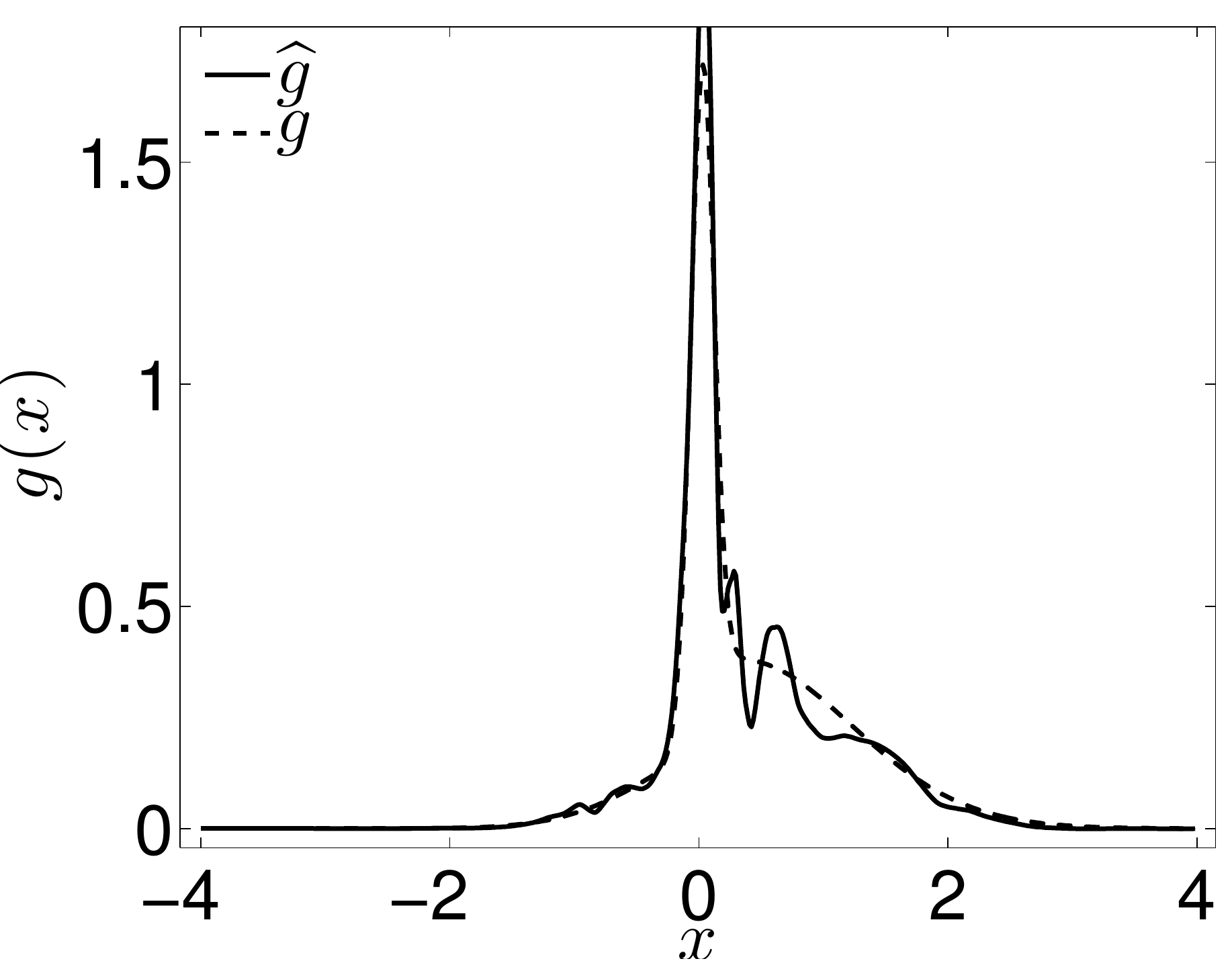}}
\caption{Typical reconstructions from a single simulation with $n=1000$ for the Kurtotic density. The dashed line depicts the original density and the solid one depicts its wavelet block estimate. (a): $\hat f(x)$. (b): $\hat F(x)$. (c): $w(\hat F(x))$. (d): $\hat g(x)=w(\hat F(x))\hat f (x)$}.
\label{fig:single}
\end{figure}

\begin{figure}[t]%[!t]
\centering
\subfigure[Uniform]{\includegraphics[width=0.245\textwidth]{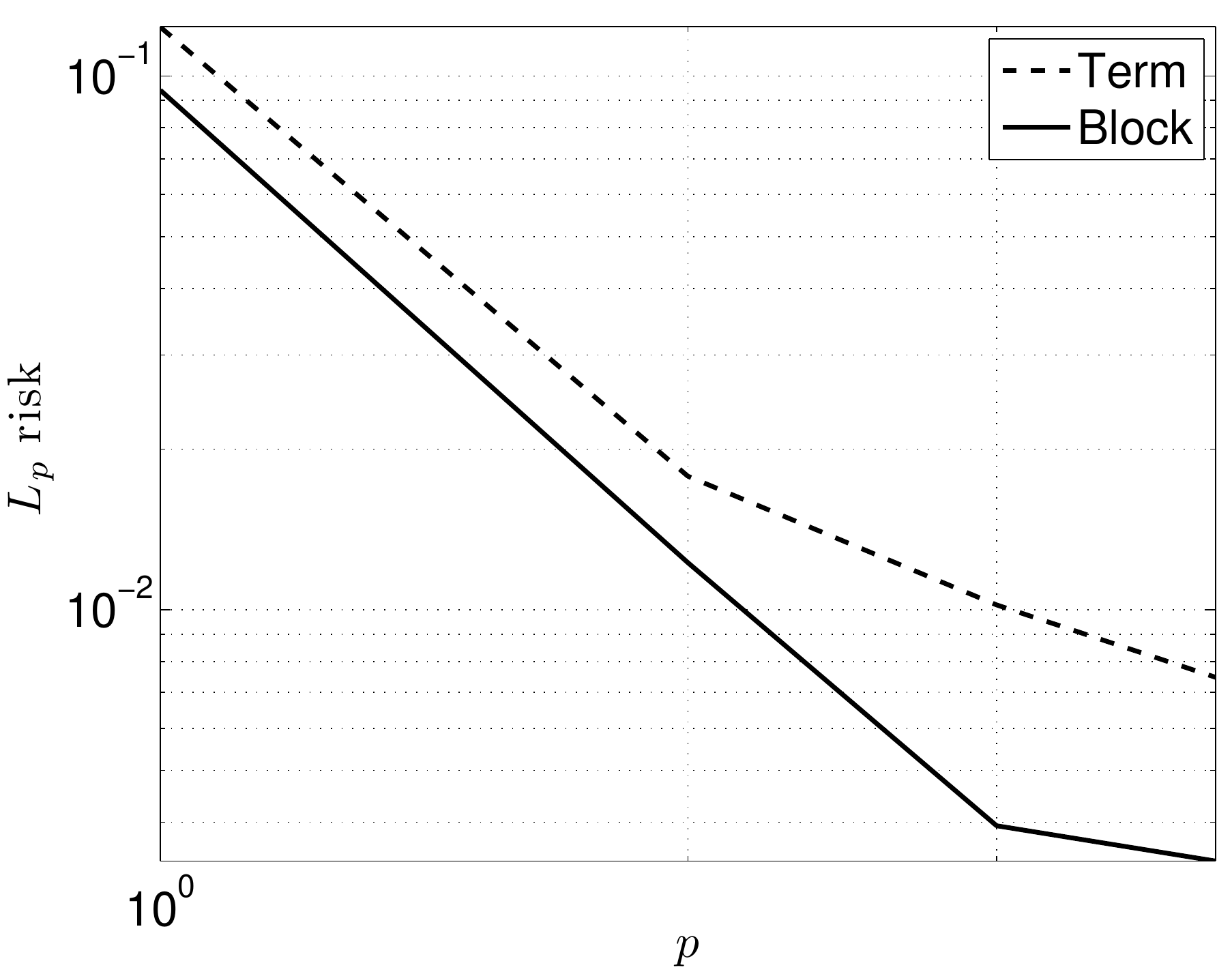}}
\subfigure[SeparatedBimodal]{\includegraphics[width=0.245\textwidth]{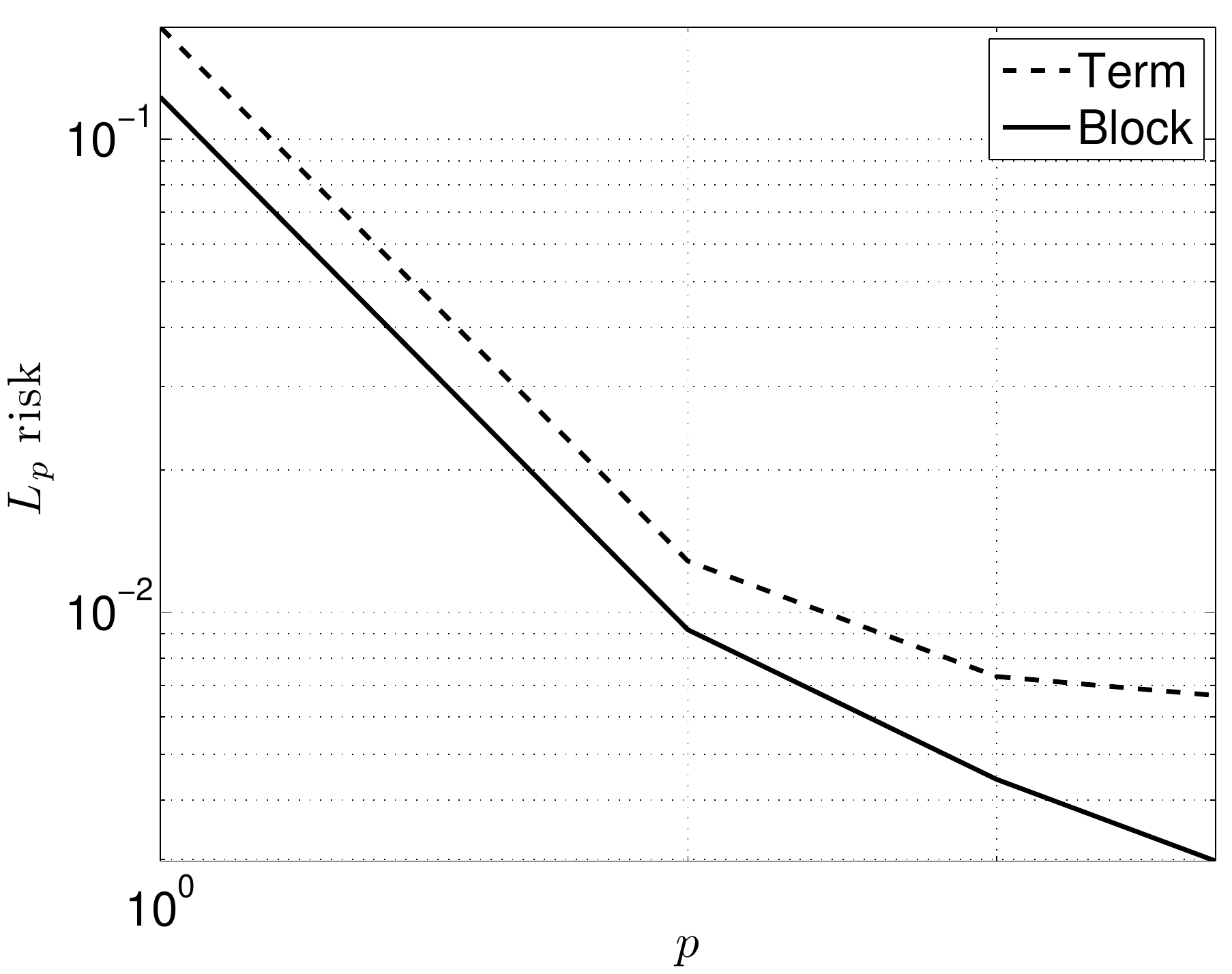}}
\subfigure[Kurtotic]{\includegraphics[width=0.245\textwidth]{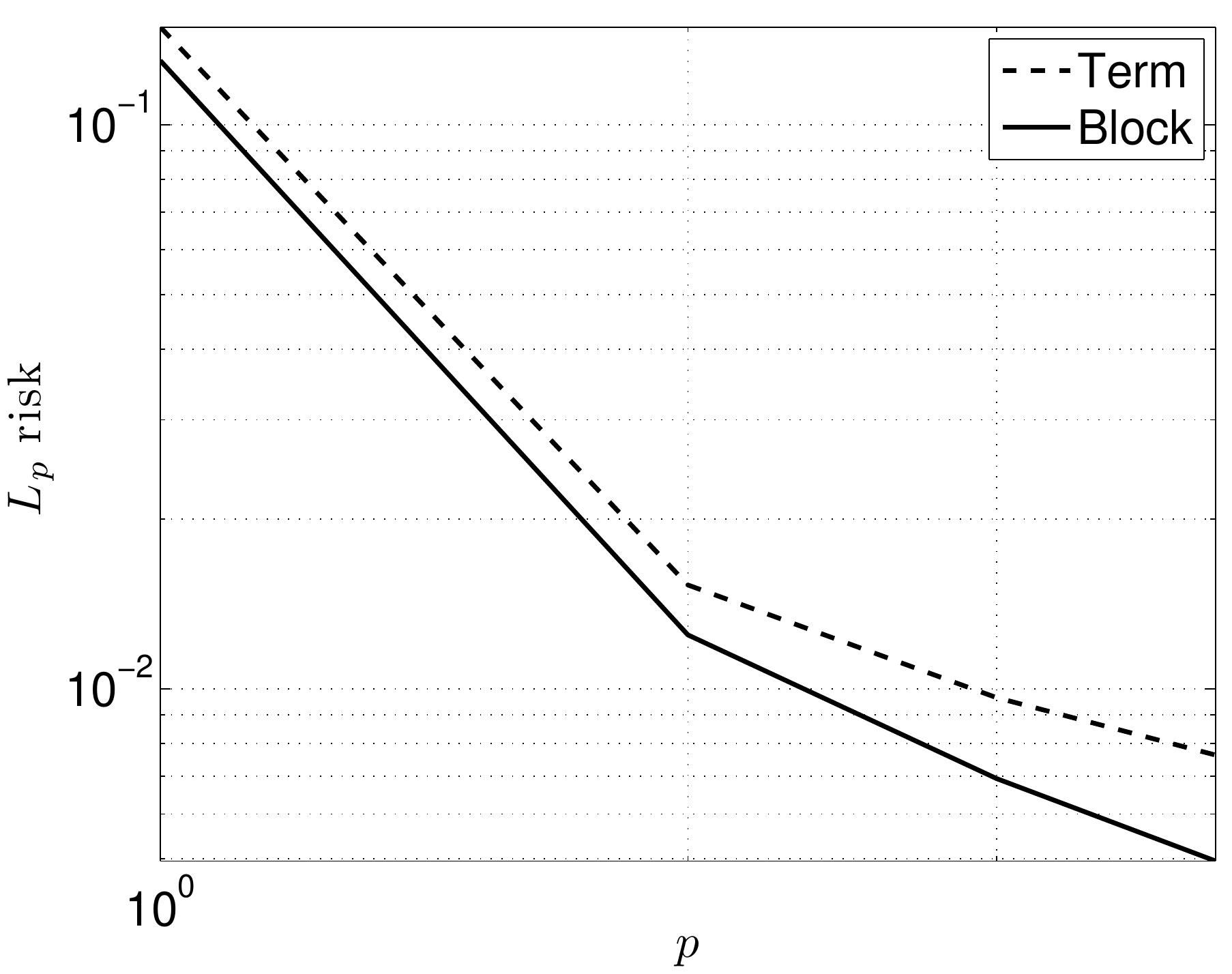}}
\subfigure[StronglySkewed]{\includegraphics[width=0.245\textwidth]{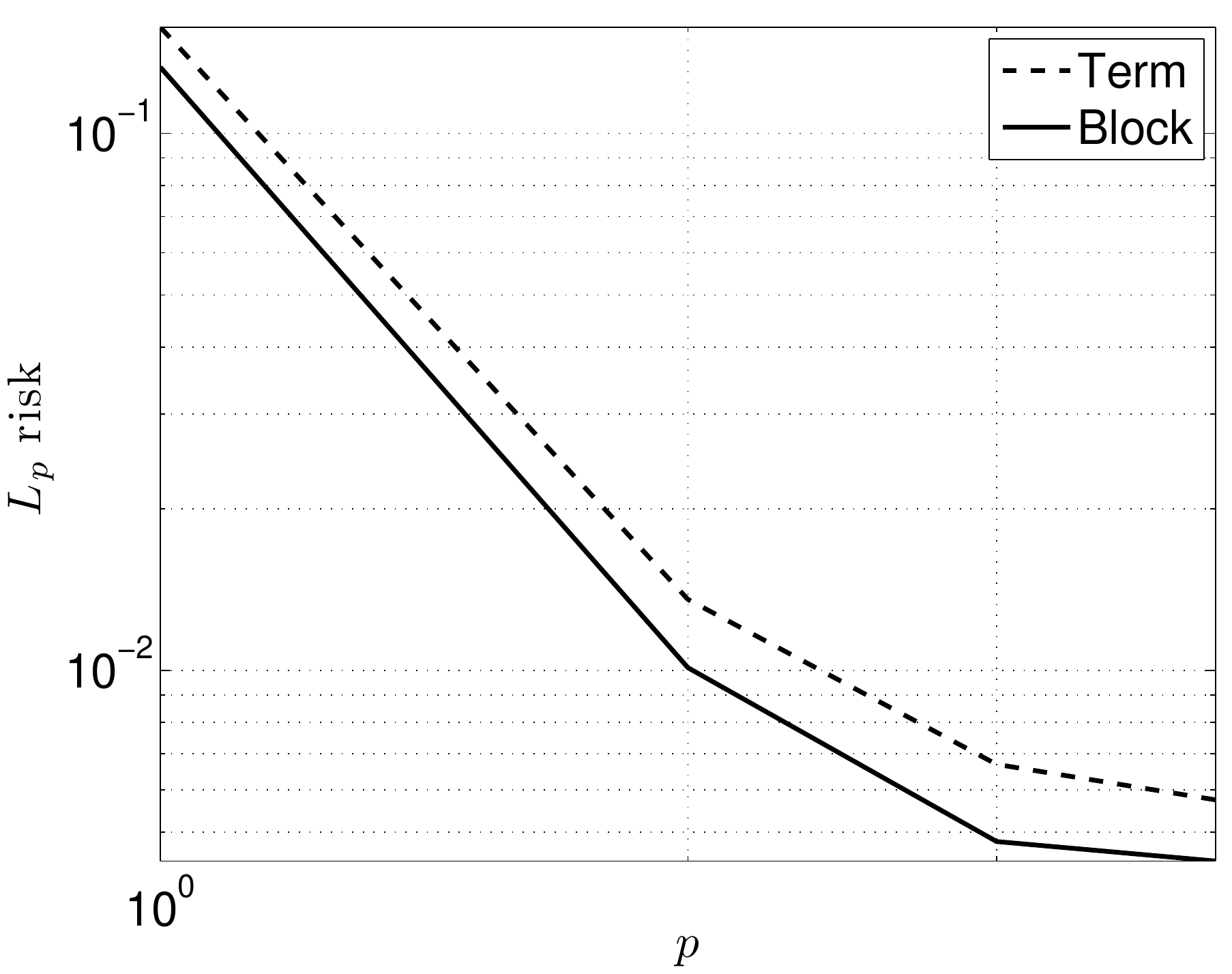}}
\caption{The influence of $p$ in the numerical values of the $\field{L}_p$ risk (in a log-log scale) of Block (solid) and term-by-term (dashed) thresholding ($L = 1$).}
\label{fig:Lprisk}
\end{figure}

\begin{figure}[t]%[!t]
\centering
\subfigure[Uniform]{\includegraphics[width=0.245\textwidth]{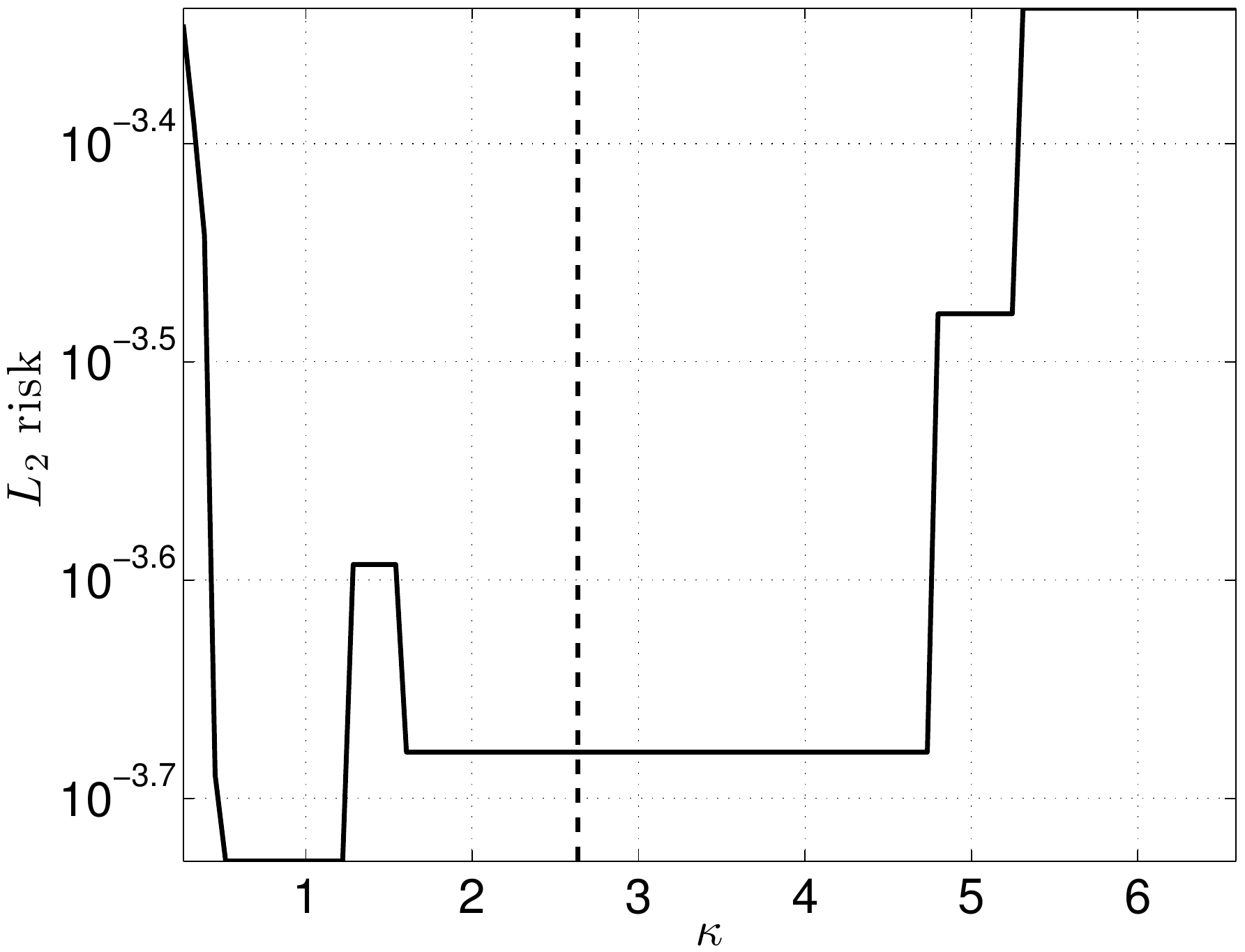}}
\subfigure[SeparatedBimodal]{\includegraphics[width=0.245\textwidth]{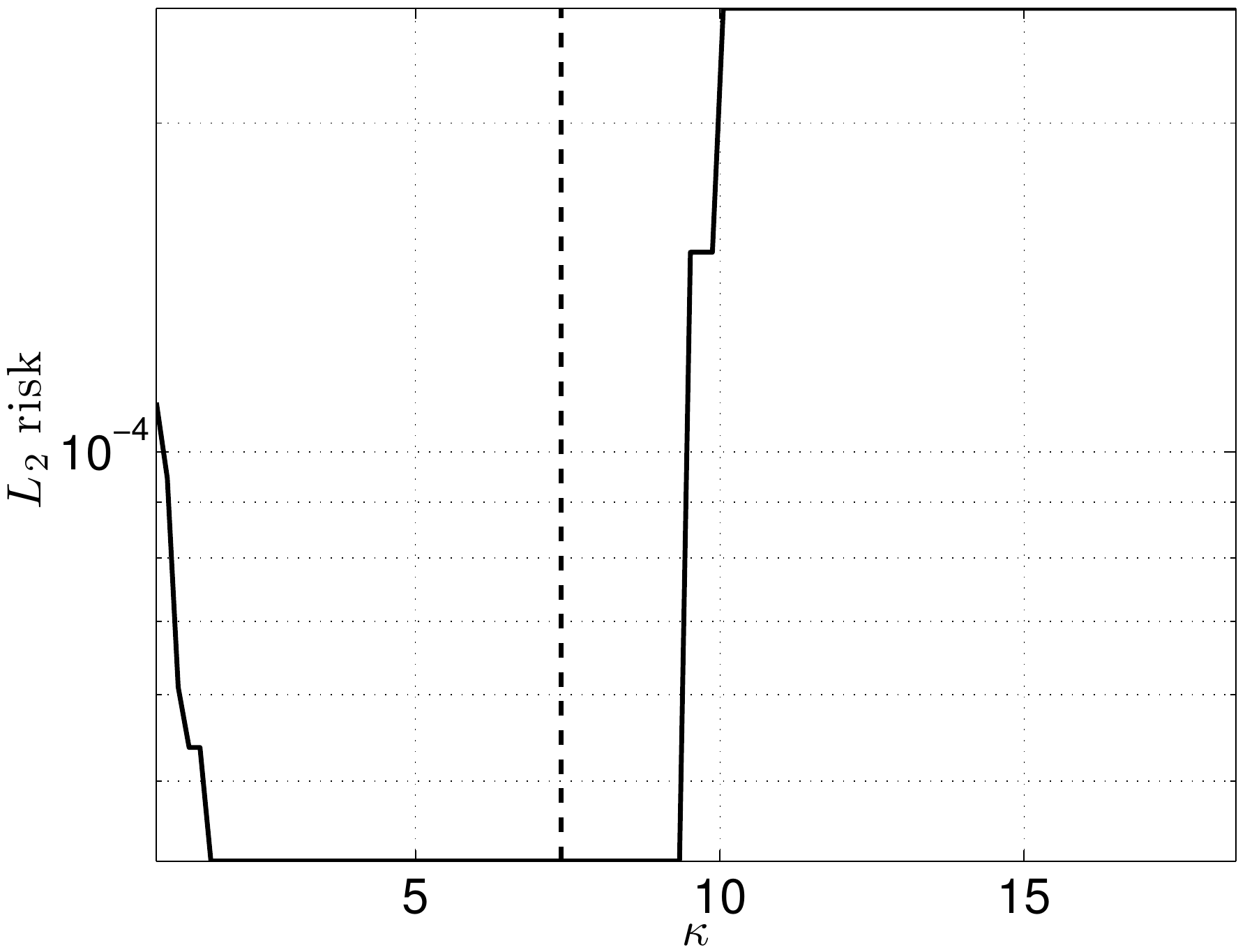}}
\subfigure[Kurtotic]{\includegraphics[width=0.245\textwidth]{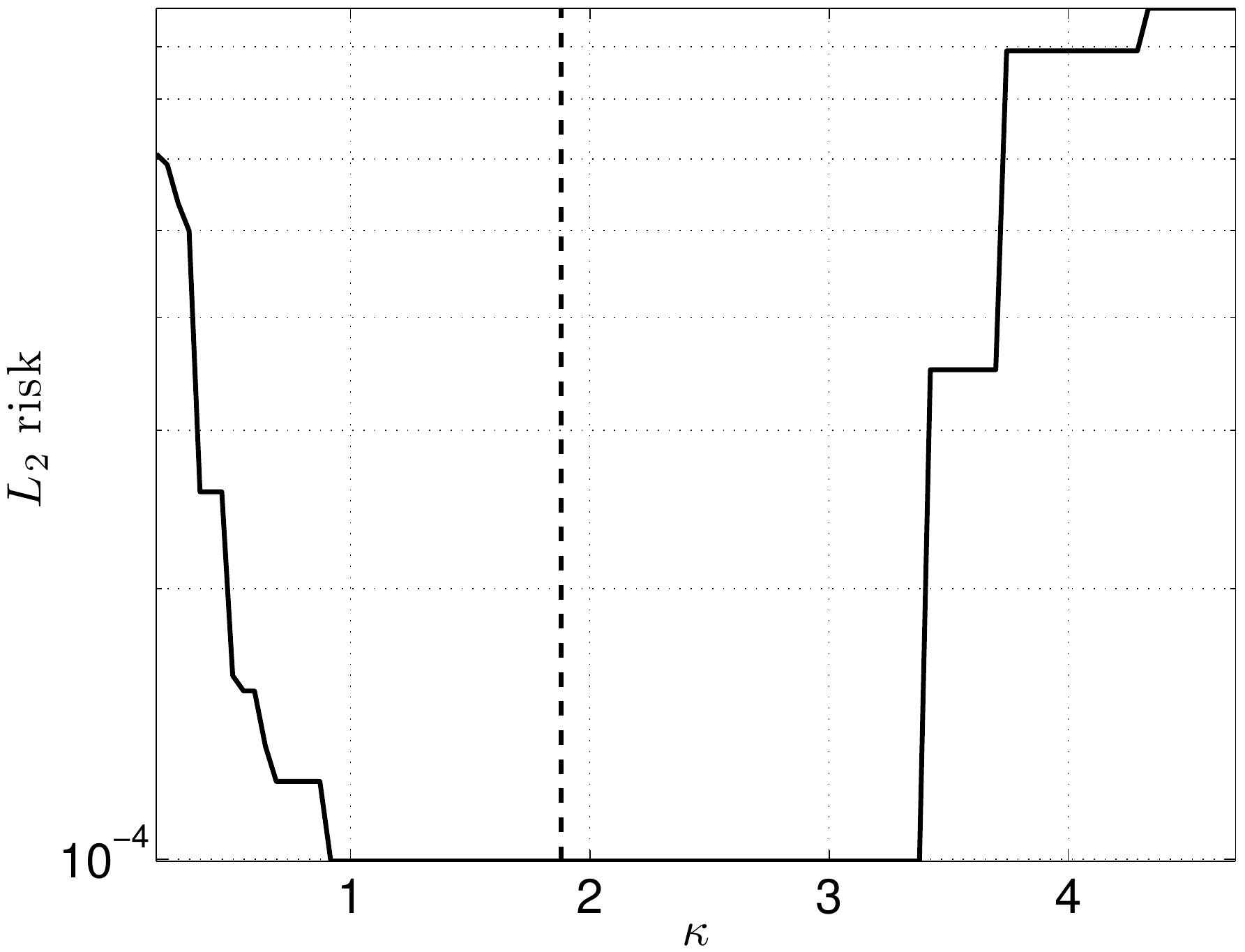}}
\subfigure[StronglySkewed]{\includegraphics[width=0.245\textwidth]{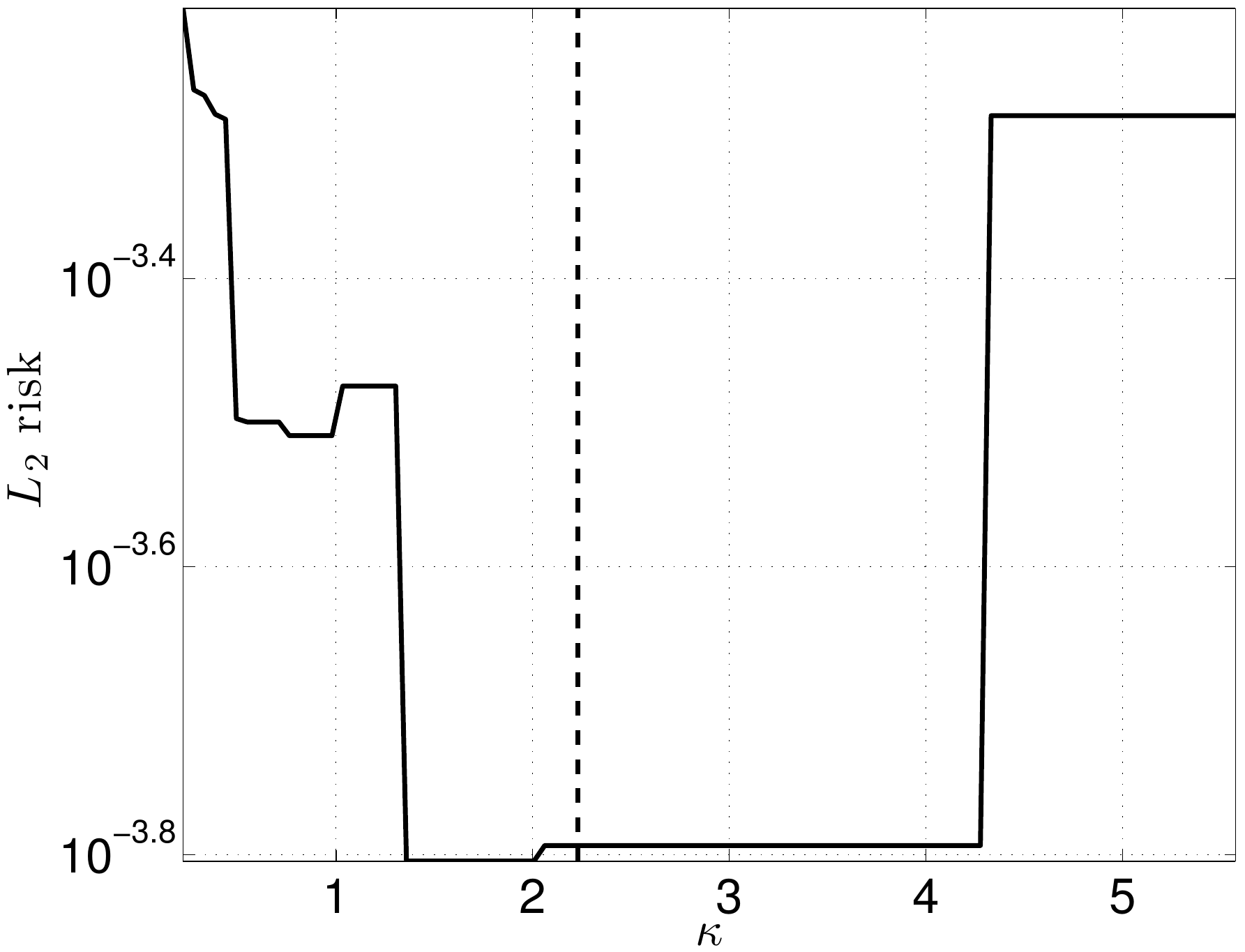}}
\caption{$\field{L}_2$ risk as a function of the threshold level $\kappa$ (in a semi-log scale), the vertical dashed lines represent the universal threshold.}
\label{fig:univ}
\end{figure}

\paragraph{Results and discussion.}
In order to illustrate Theorem~\ref{danzig2}, we study the influence of $p$ on the numerical performances of the Block and the term-by-term ($L=1$) thresholding estimator. Let us first consider a parallel system with $m=2$ identical independent components. Then, the corresponding weighted function is \eqref{fixed} and the goal is to estimate $g$ in \eqref{densityy} from $X_1,\ldots,X_n$ sample simulated from one of the test densities. A typical example of estimation for the Kurtotic density (for $p = 2$), with $n=1000$ is given in Figure~\ref{fig:single}. The mean $\field{L}_p$ risk of $\hat g$ i.e., $R_p(\hat g,g)=(1/T)\sum_{i=-T/2}^{T/2-1}|\hat g(t_i)-g(t_i)|^p$, is obtained with 10 samples for $n=1000$, and it is plotted as a function of $p$ in Figure~\ref{fig:Lprisk}. As predicted by Theorem~\ref{danzig2}, the larger $p$, the smaller $\field{L}_p$ risk of $\hat g$. We can see that our estimation procedure provides better results than the term-by-term thresholding ($L=1$) in all cases. In particular, the risk improvement achieved by the block estimators upon the term-by-term estimator is significant for the non-smooth Uniform density. This is in agreement with the predictions of our theoretical findings. 

To conclude this first example, we illustrate from a single simulation, the fact that the parameters dictated by the theory yield the expected performance. We display in Figure~\ref{fig:univ} the empirical $\field{L}_2$ risk as a function of the threshold level $\kappa$, where the vertical dashed line represents the universal threshold. One can see that the minimum of the $\field{L}_2$ risk is close to the universal threshold for all test densities, thus supporting the choice dictated by our theoretical procedure, although derived in an asymptotic setting.

\begin{figure}[t]
\centering
\subfigure[]{
\includegraphics[width=0.24\textwidth]{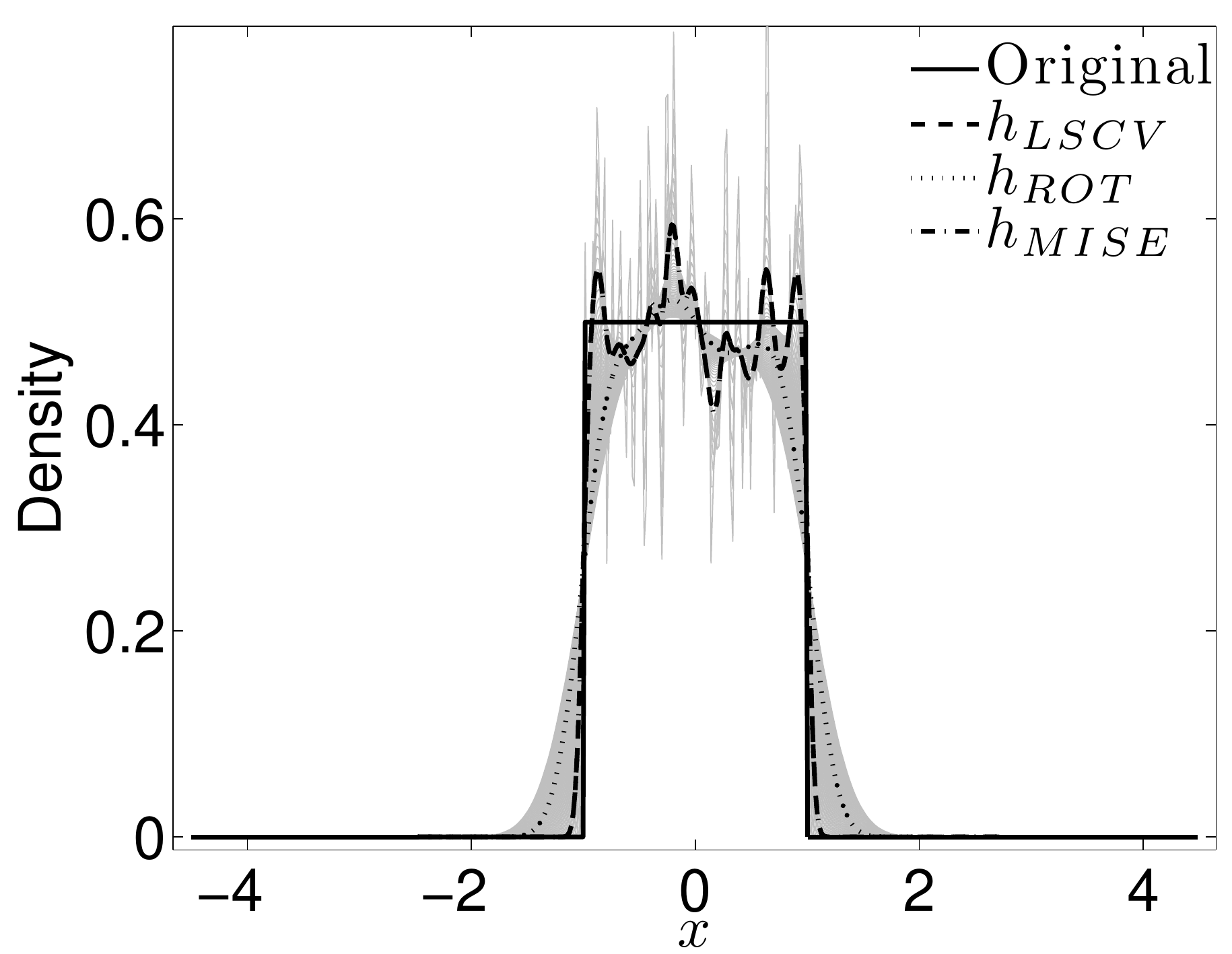}
\includegraphics[width=0.24\textwidth]{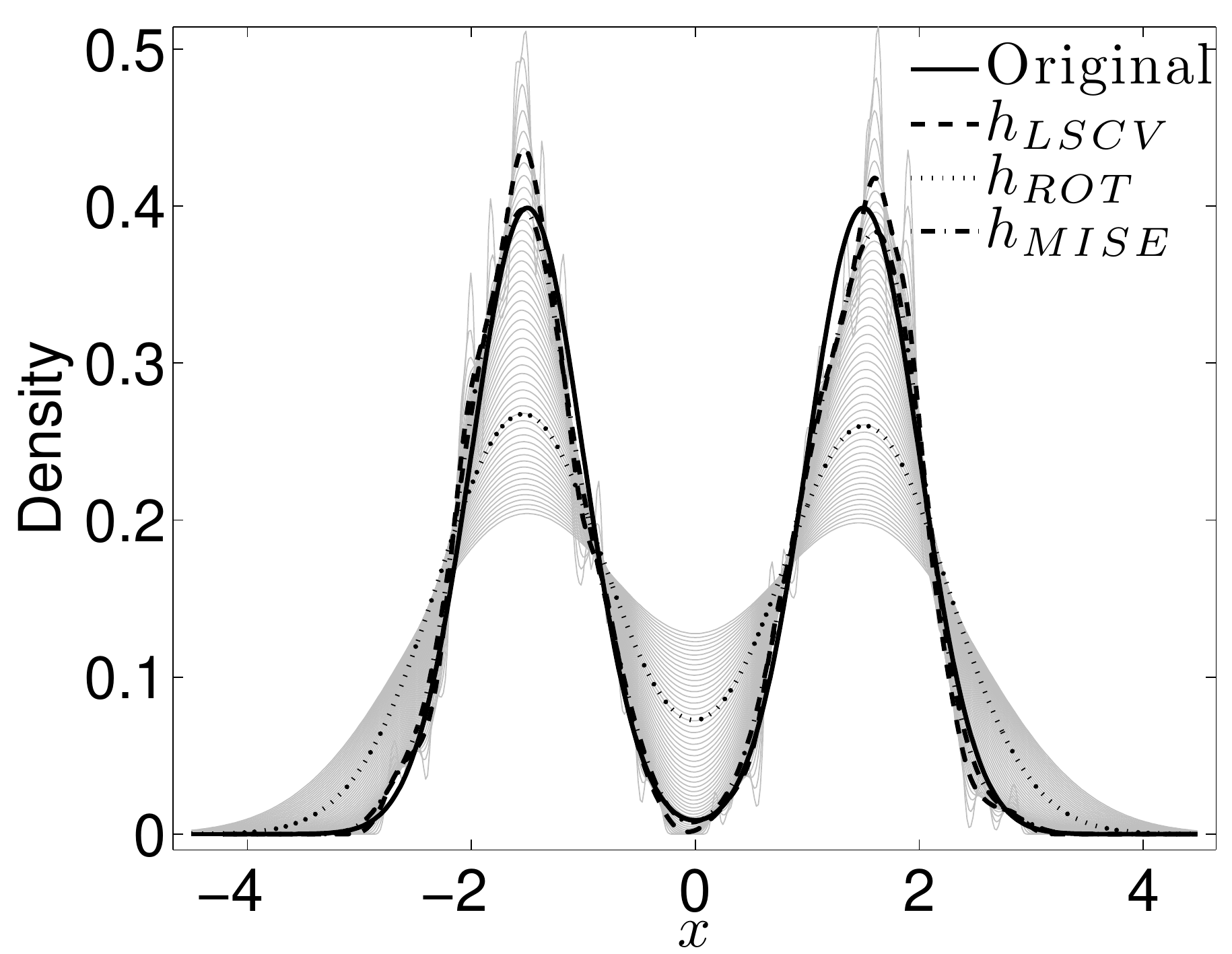}
\includegraphics[width=0.24\textwidth]{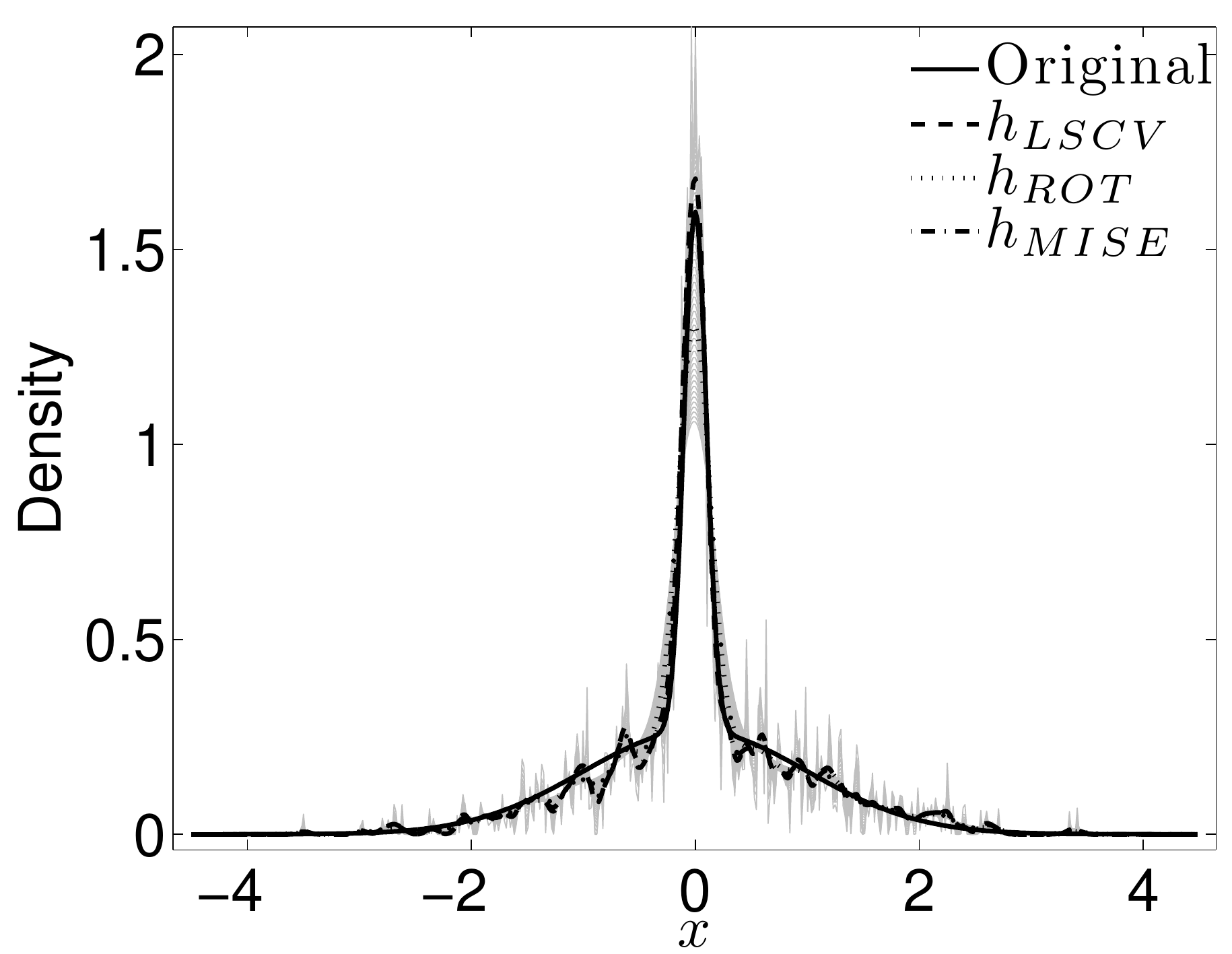}
\includegraphics[width=0.24\textwidth]{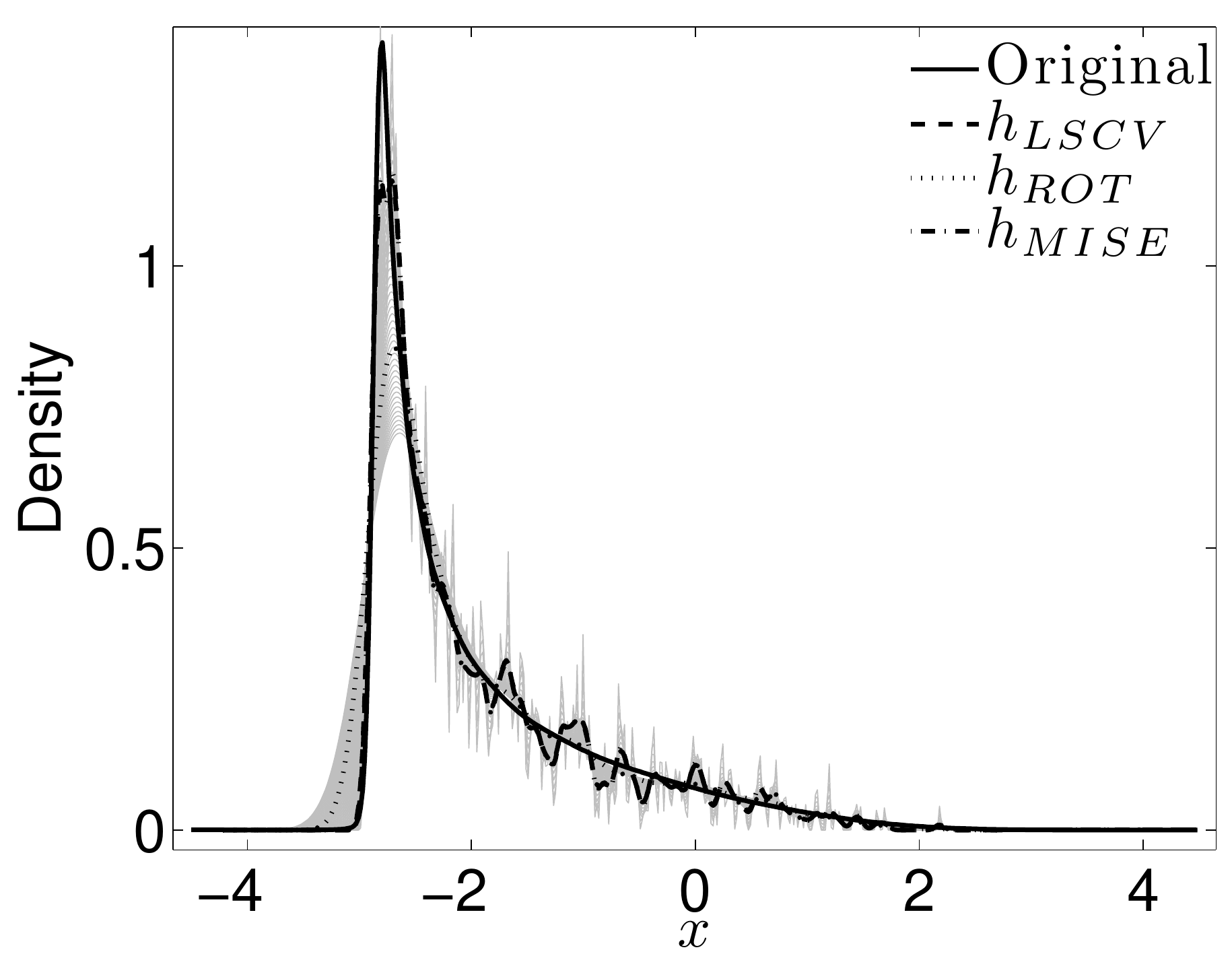}}
\subfigure[]{
\includegraphics[width=0.24\textwidth]{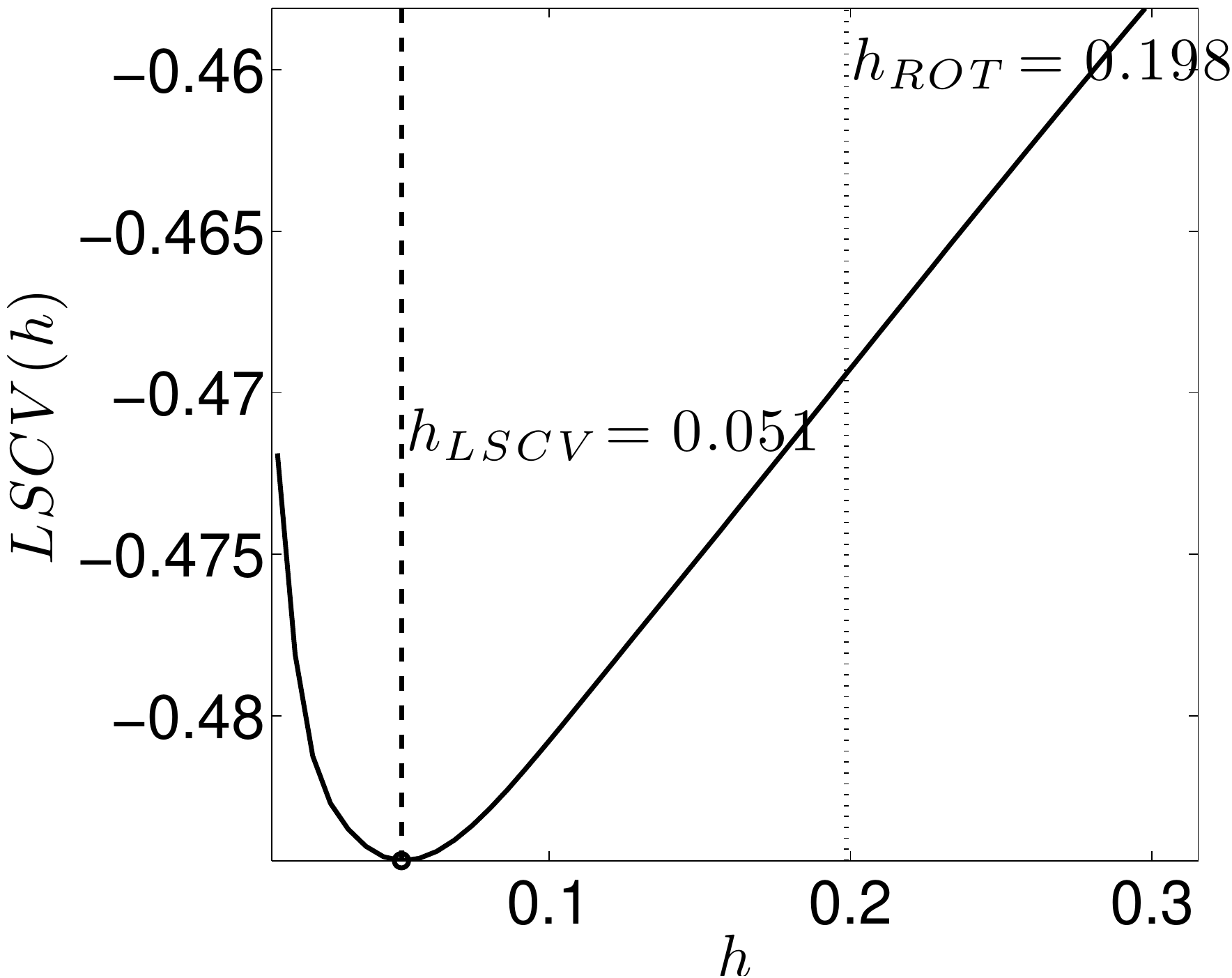}
\includegraphics[width=0.24\textwidth]{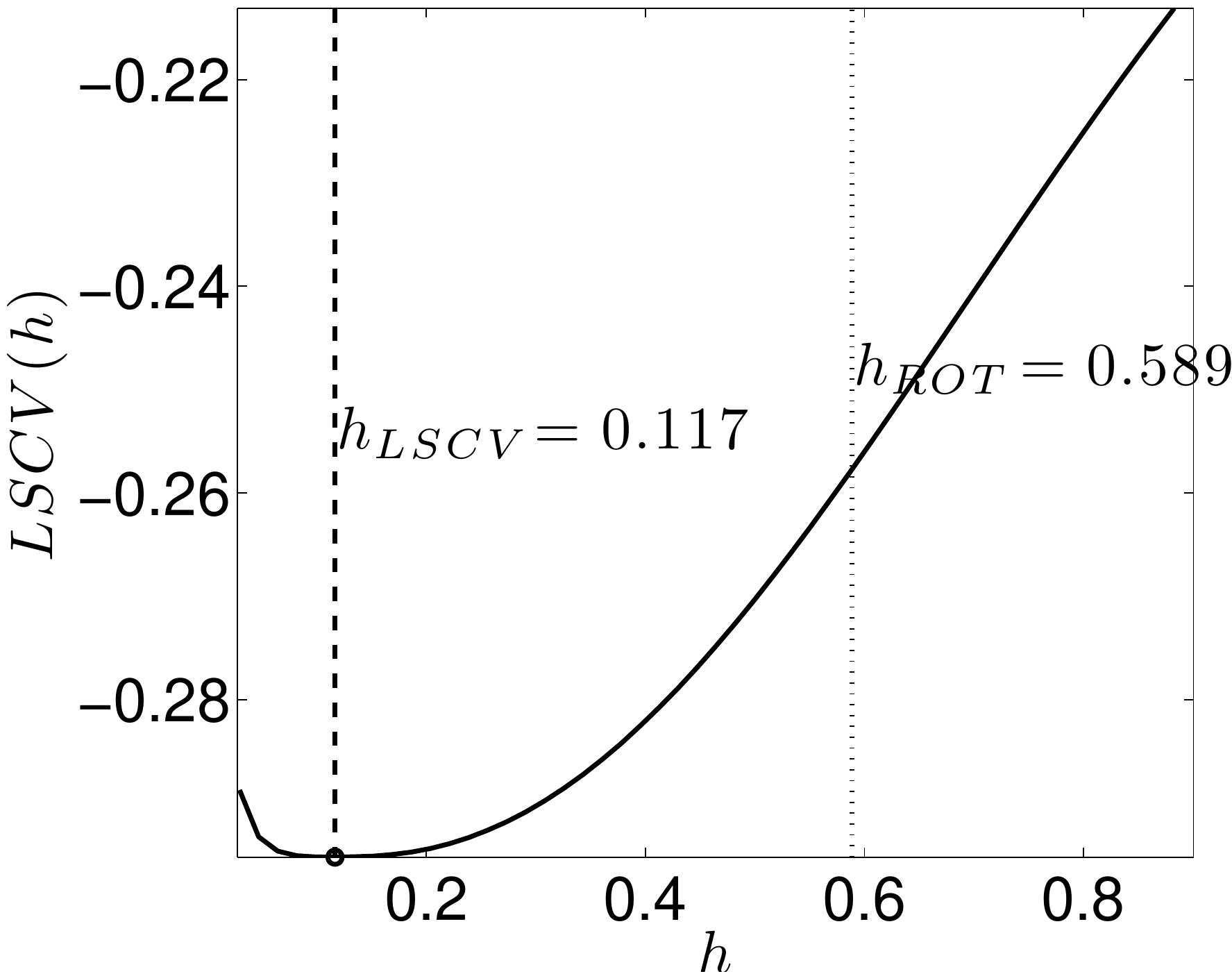}
\includegraphics[width=0.24\textwidth]{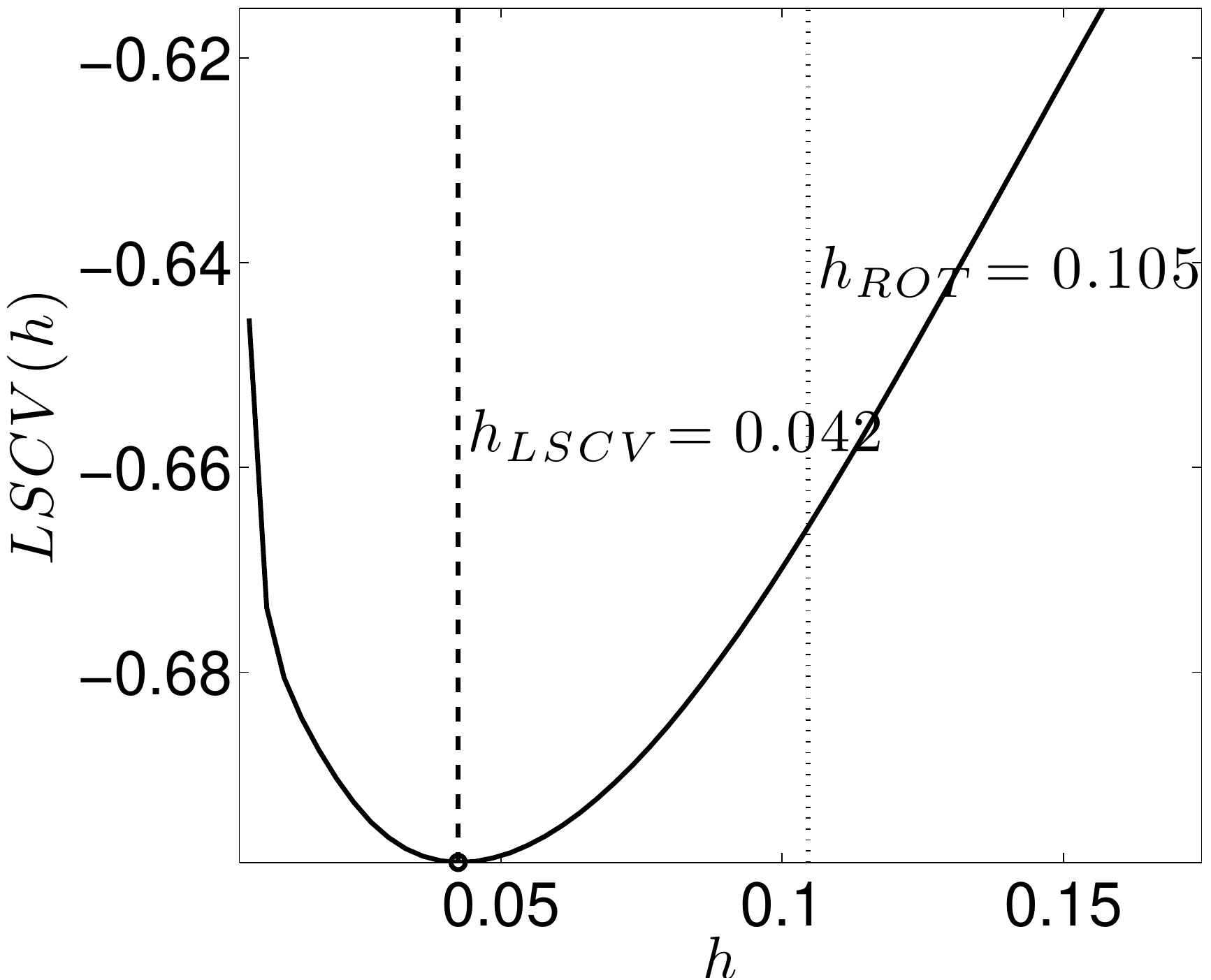}
\includegraphics[width=0.24\textwidth]{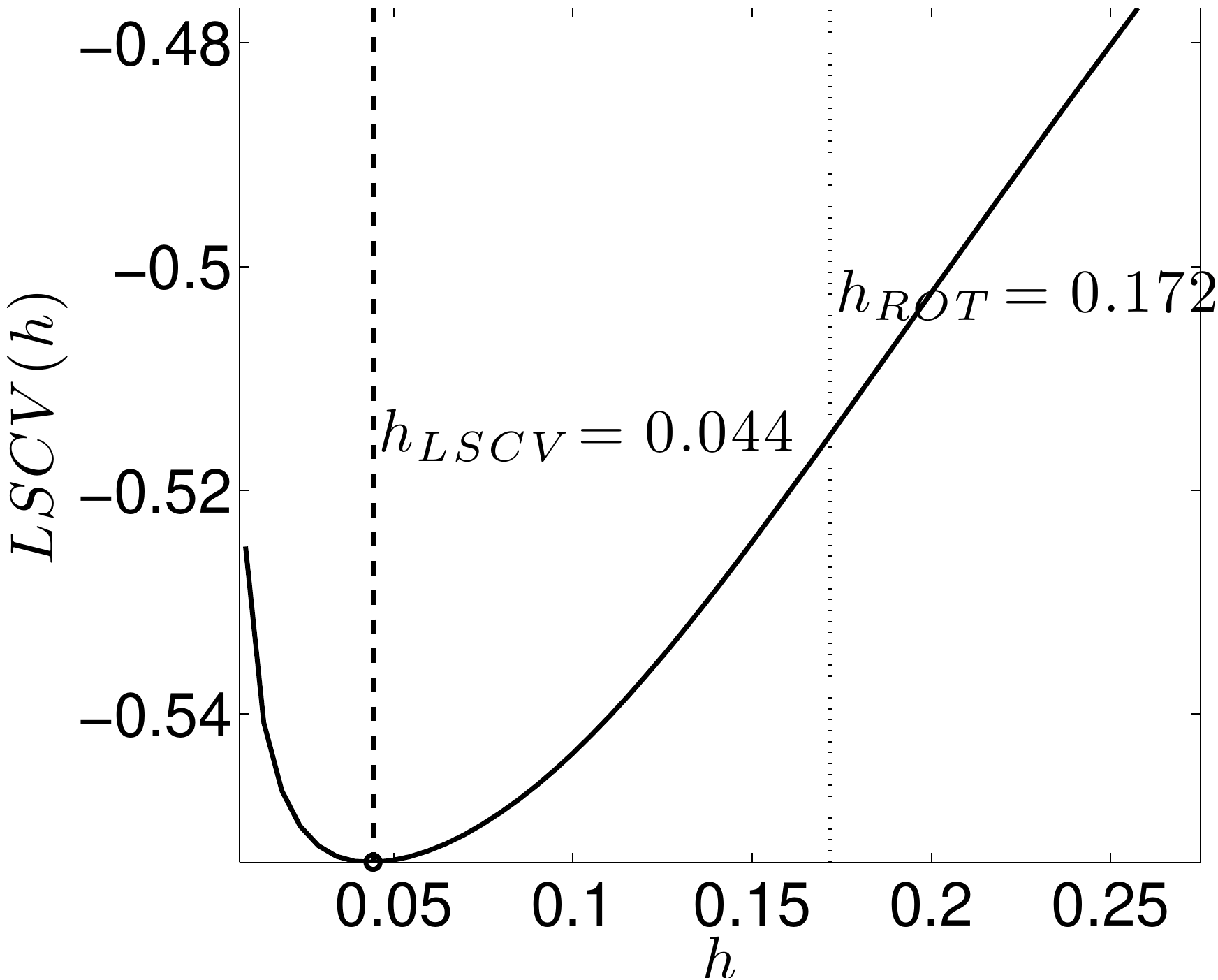}}
\subfigure[]{
\includegraphics[width=0.24\textwidth]{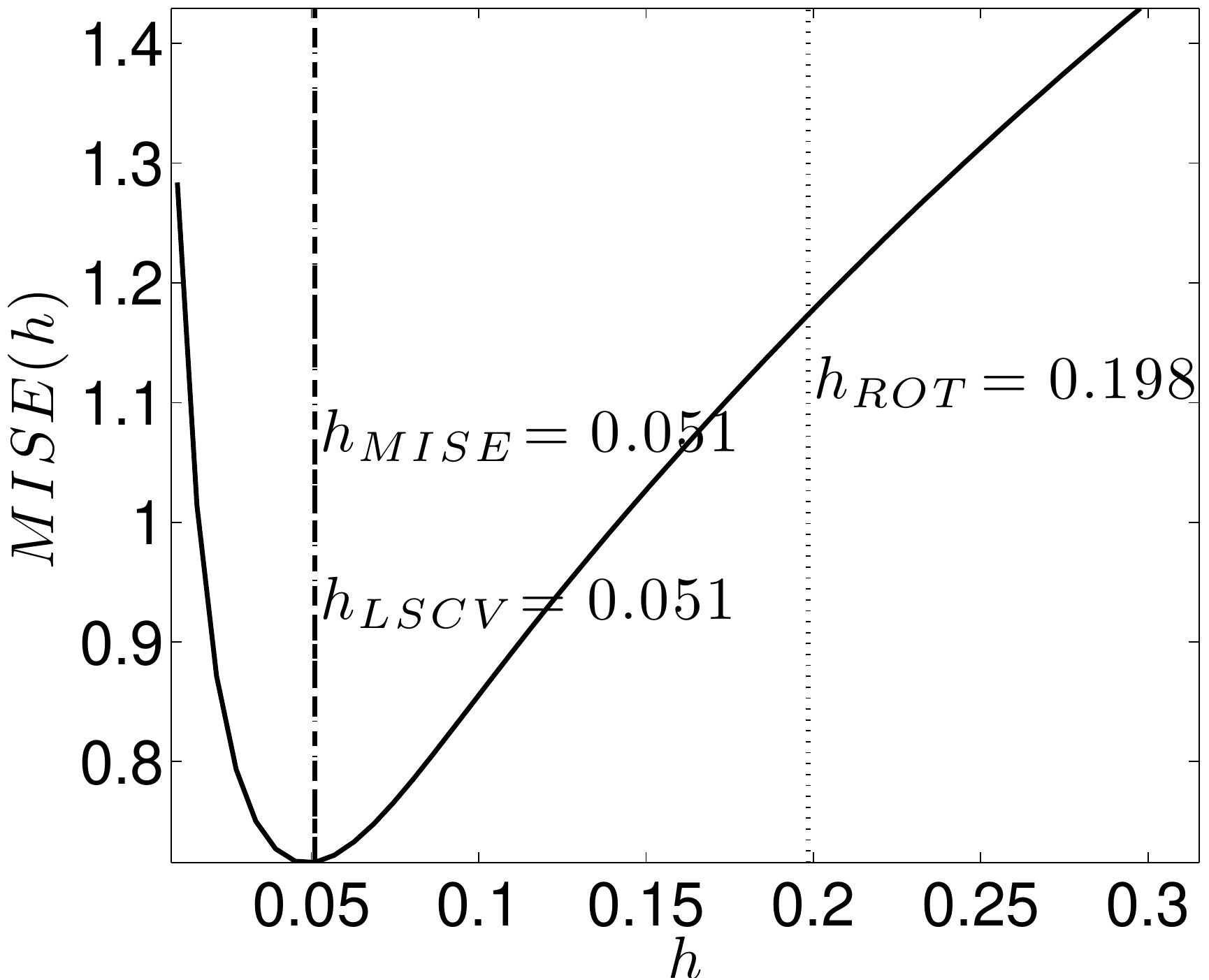}
\includegraphics[width=0.24\textwidth]{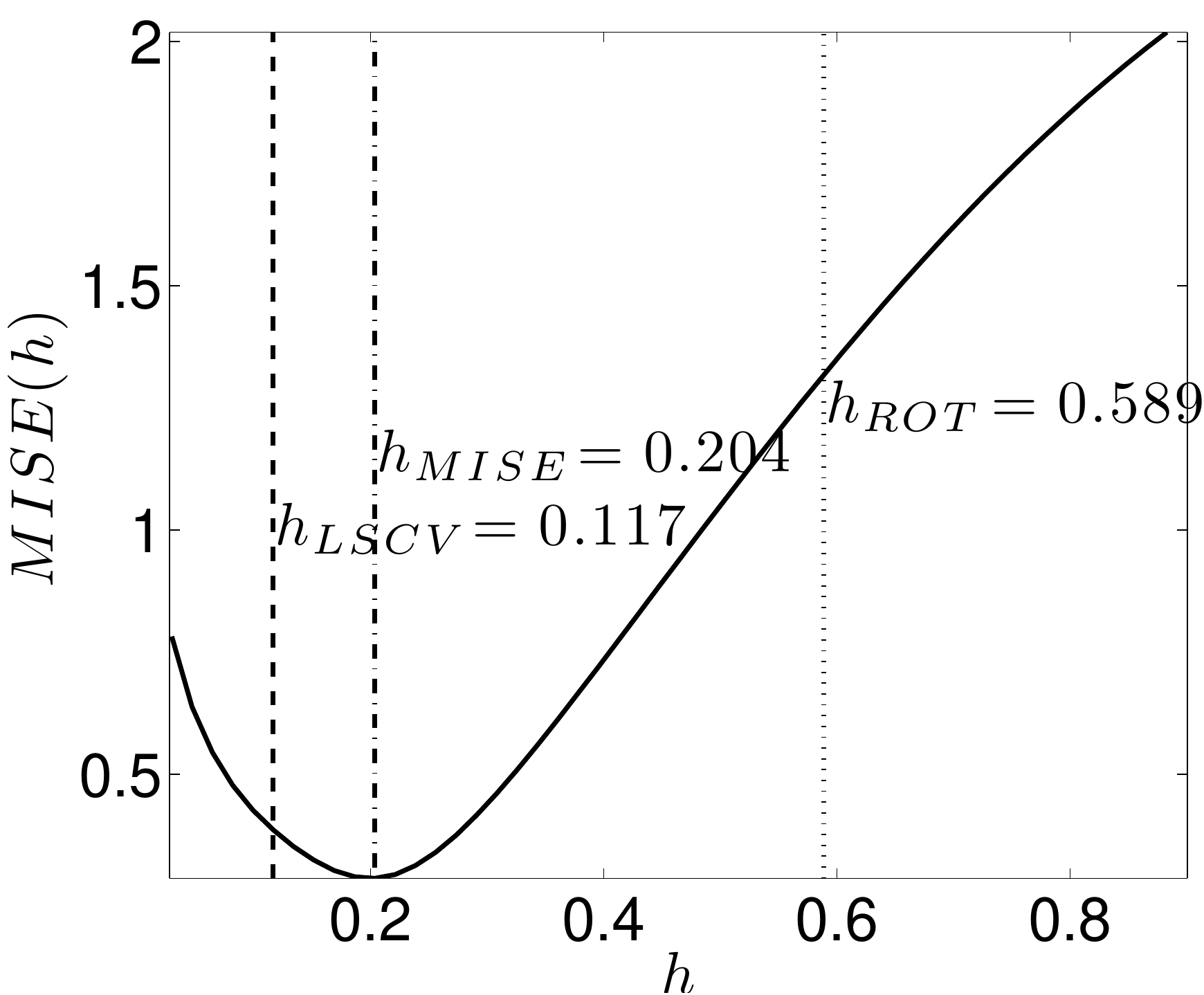}
\includegraphics[width=0.24\textwidth]{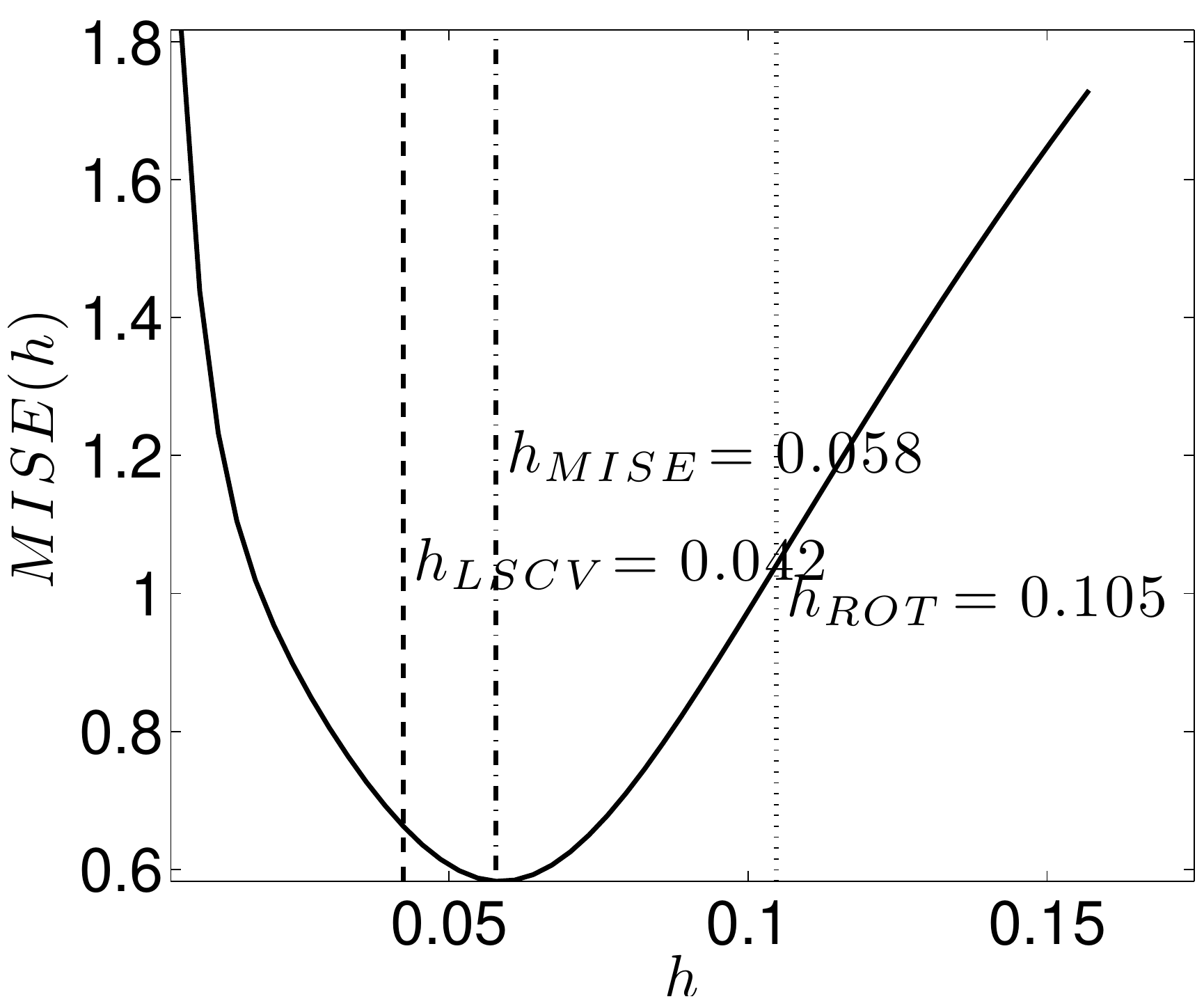}
\includegraphics[width=0.24\textwidth]{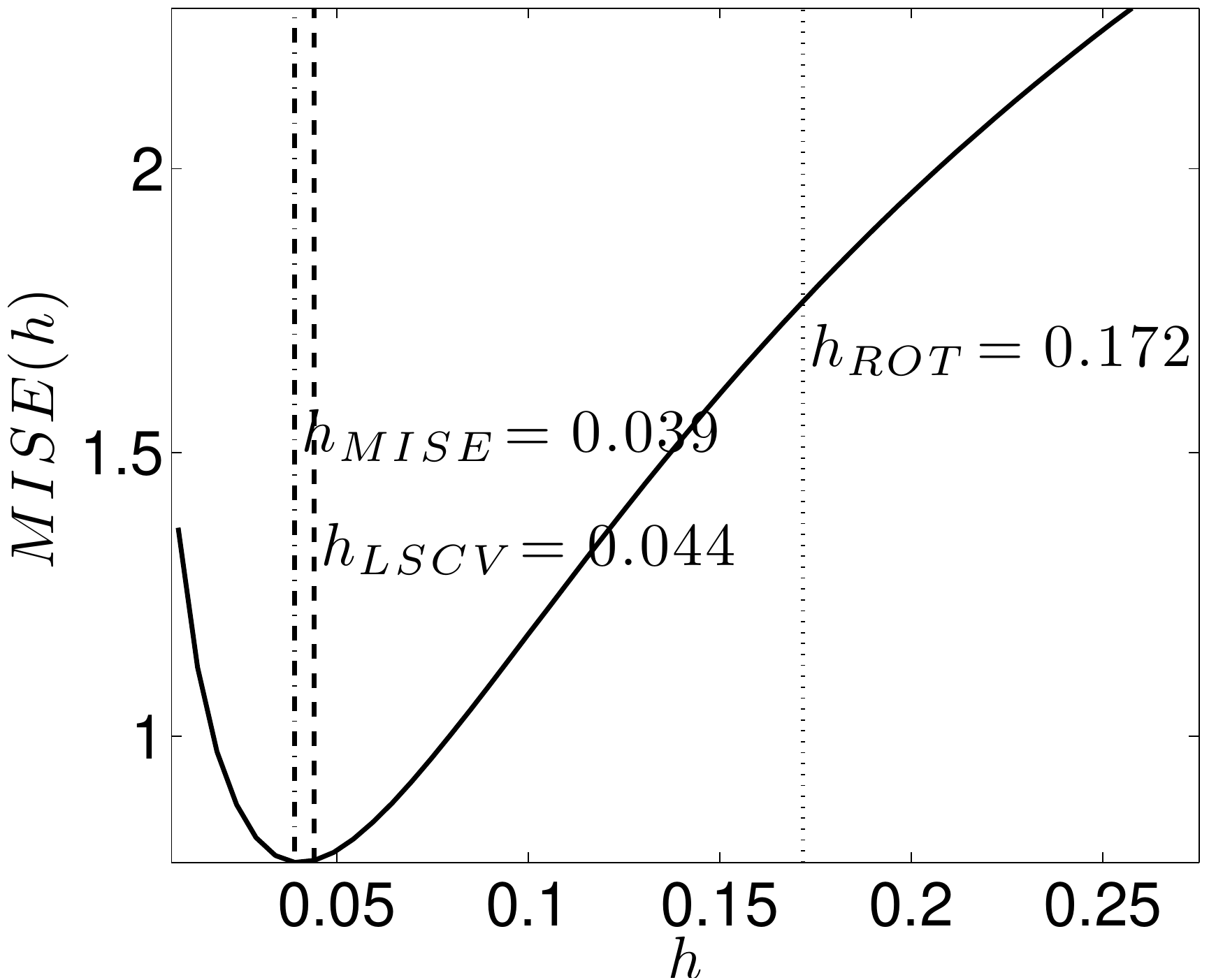}}
\caption{(a): Density estimates. (b): Graph of the the LSCV function versus the kernel bandwidth $h$ for each of the tested densities, the vertical dashed lines represent the value of $h$ that minimizes $\mathrm{LSCV}(h)$. (c): The solid line depicts the estimated MISE as a function of $h$, the vertical dashed-dotted lines represent the true MISE-minimizing bandwidth $h_{\mathrm{MISE}}$ and the vertical dotted lines represent the pilot bandwidth $h_{\mathrm{ROT}}$.}
\label{fig:bandw}
\end{figure}

In our second example, the adaptive estimator described in Section \ref{applications1} is tested when $N$ follows a Geometric distribution, so that the weight function is that given by \eqref{geo}. This example is devoted to a simulation study comparing the performance of the block hard thresholding estimator with that of  the traditional kernel defined as follows
\begin{equation}\label{kernel}
 \hat f_{h}(x) = \frac{1}{nh}\sum_{i=1}^nK\left(\frac{x-X_i}{h}\right),
\end{equation} 
where the positive kernel $K$ satisfies $\int K(x)dx=1$ and the smoothing parameter $h$ is known as the bandwidth.

Many procedures of bandwidth selection for kernel density estimation have been developed in the literature (see, e.g.,~\cite{Sil}).
We use least-squares cross-validation (LSCV) (\cite{rud},~\cite{bow}) where the bandwidth is defined as
\begin{equation*}
h_{\mathrm{LSCV}}= \underset{h}{\operatorname{argmin}}  { \int_{-b}^{b} \hat{f}_h(x)^2dx-2n^{-1}\sum_{i=1}^n\hat f_{-i}(X_i)},
\end{equation*}
and $\hat f_{-i}$ is the \textit{leave-one  out} kernel  estimator constructed from the data without the observation $X_i$. It is motivated by the fact that for independent data
 \begin{equation*}
 \mathrm{LSCV}(h)= \int_{-b}^{b} \hat{f}_h(x)^2dx-2n^{-1}\sum_{i=1}^n\hat f_{-i}(X_i)
 \end{equation*}
 is an unbiased estimator of $\mathrm{MISE}(h)=\int_{-b}^b f^2(x)dx$.
%We use a criterion based on cross-validation (CV): the least-squares cross-validation (LSCV) (\cite{rud},~\cite{bow}). 
One frequently used cross-validation (CV) procedure is the K-fold CV (as described e.g. in \cite[Section 7.10]{hastibfri}) in which the data set $X_1, \ldots, X_n$ is randomly partitioned into K approximately equal-sized and non-overlapping subsets $S1,\ldots, S_K$. To obtain the bandwidths $h_{\mathrm{LSCV}}$, we have performed a $10$-fold CV, using a Gaussian kernel, with a simple ``rule-of-thumb'' pilot bandwidth $h_{\mathrm{ROT}}$. Figure~\ref{fig:bandw}(b) contains a plot of the LSCV function versus the kernel bandwidth $h$ and Figure~\ref{fig:bandw}(c) the estimated MISE as a function of $h$. For each density, it is clear from this figure (Figure~\ref{fig:bandw}(b)), that  the value of $h_{\mathrm{LSCV}}$ is the unambiguous minimizer of $\mathrm{LSCV}(h)$. We see that $h_{\mathrm{LSCV}}$  provides a decent approximation, close to $h_{\mathrm{MISE}}$ for all test densities. For the StronglySkewed density,  the bandwidth which minimizes $\mathrm{MISE}(h)$ in this case is $h_{\mathrm{MISE}}=0.039$ and $h_{\mathrm{LSCV}}=0.044$. In this case, for the Uniform density, $h_{\mathrm{MISE}}=h_{\mathrm{LSCV}}=0.051$.

\begin{figure}[htbp]
\centering
\includegraphics[width=0.24\textwidth]{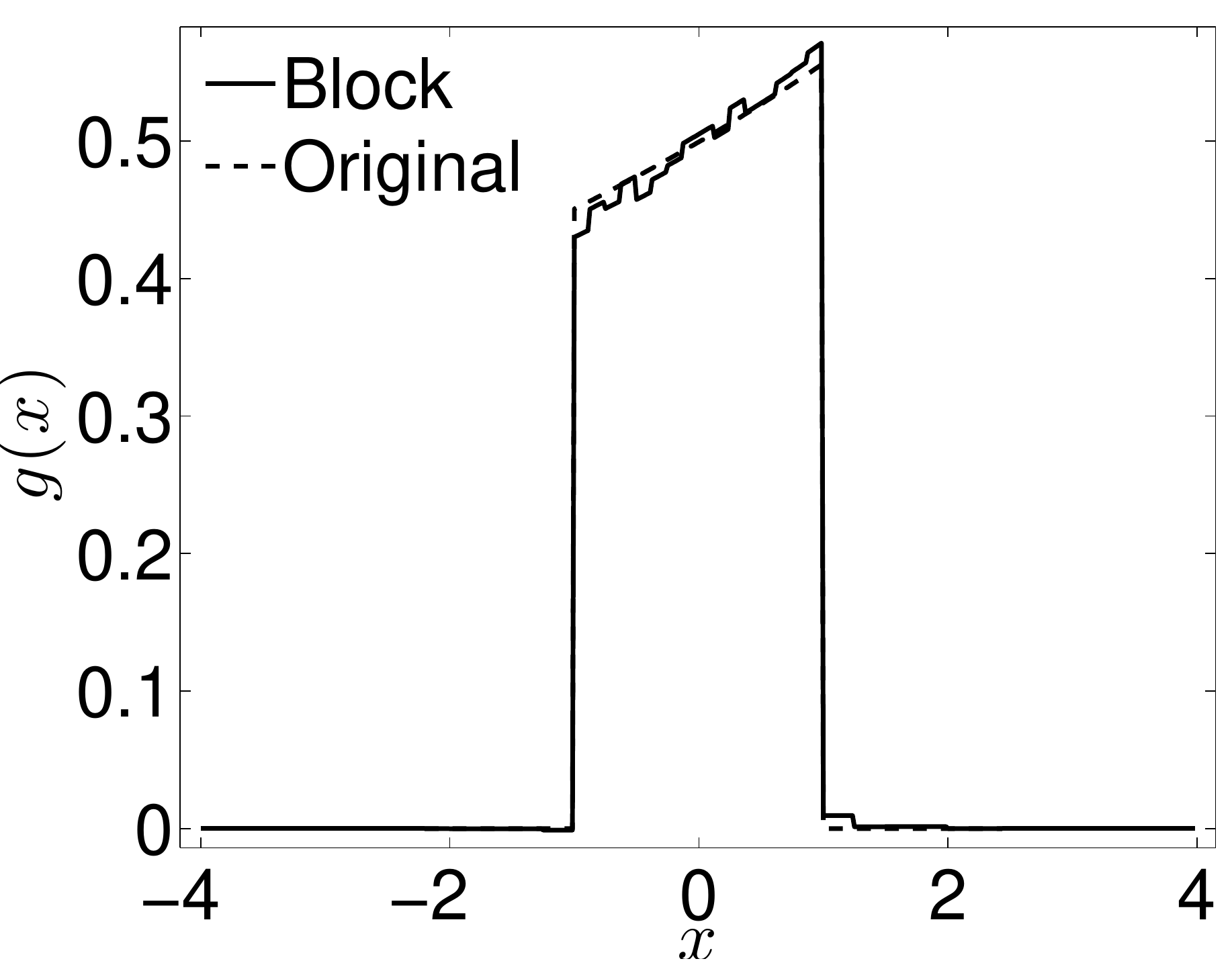}
\includegraphics[width=0.24\textwidth]{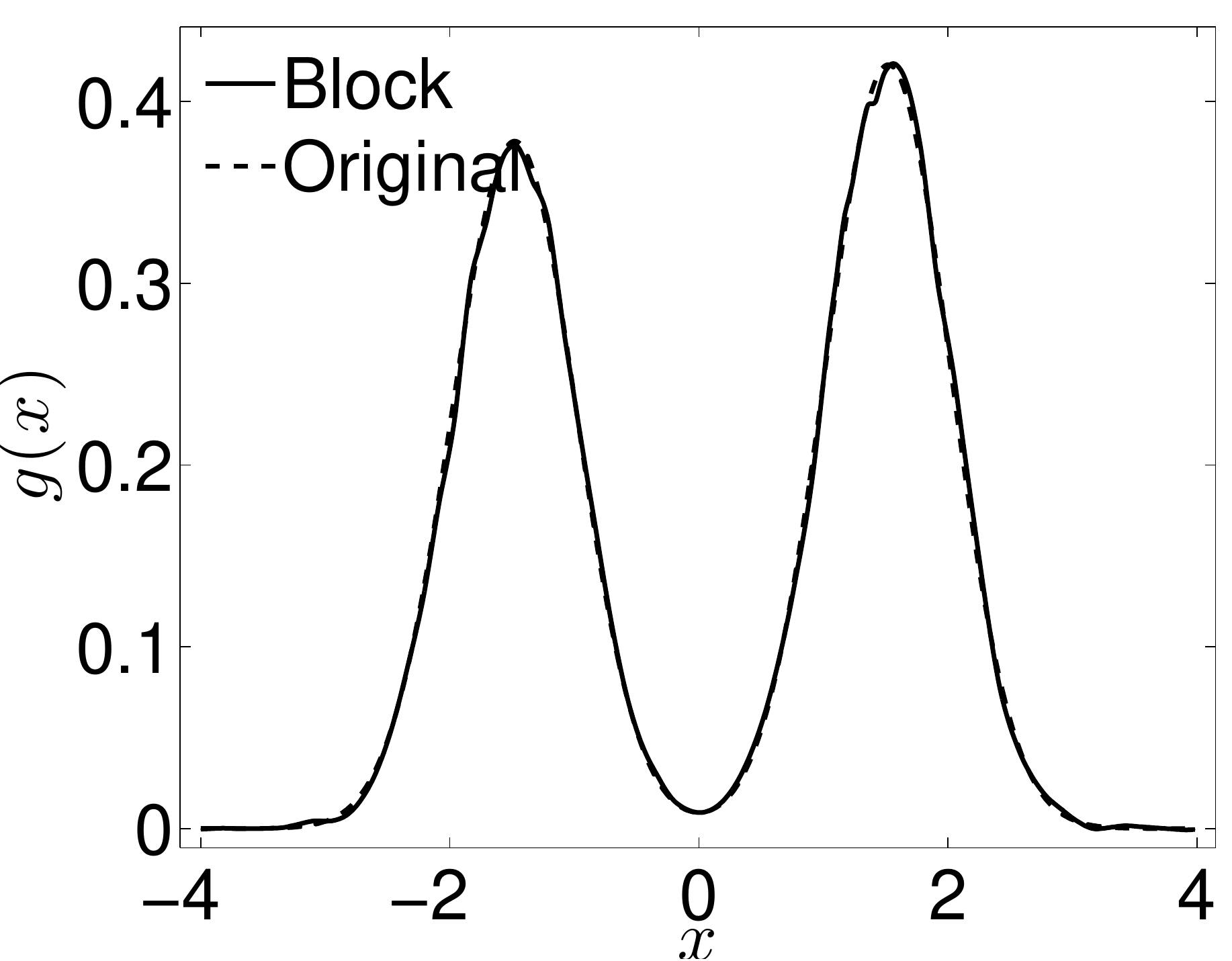}
\includegraphics[width=0.24\textwidth]{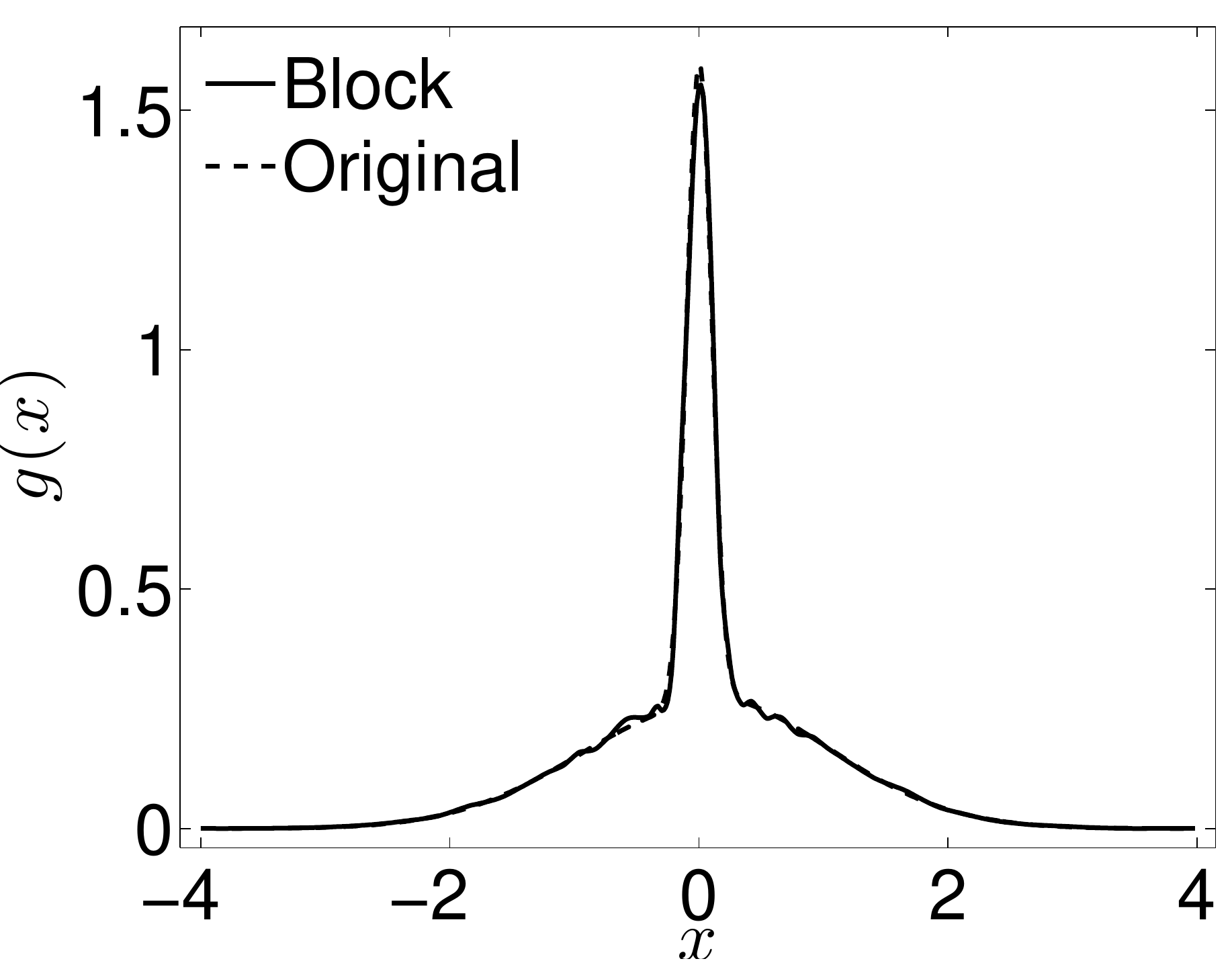}
\includegraphics[width=0.24\textwidth]{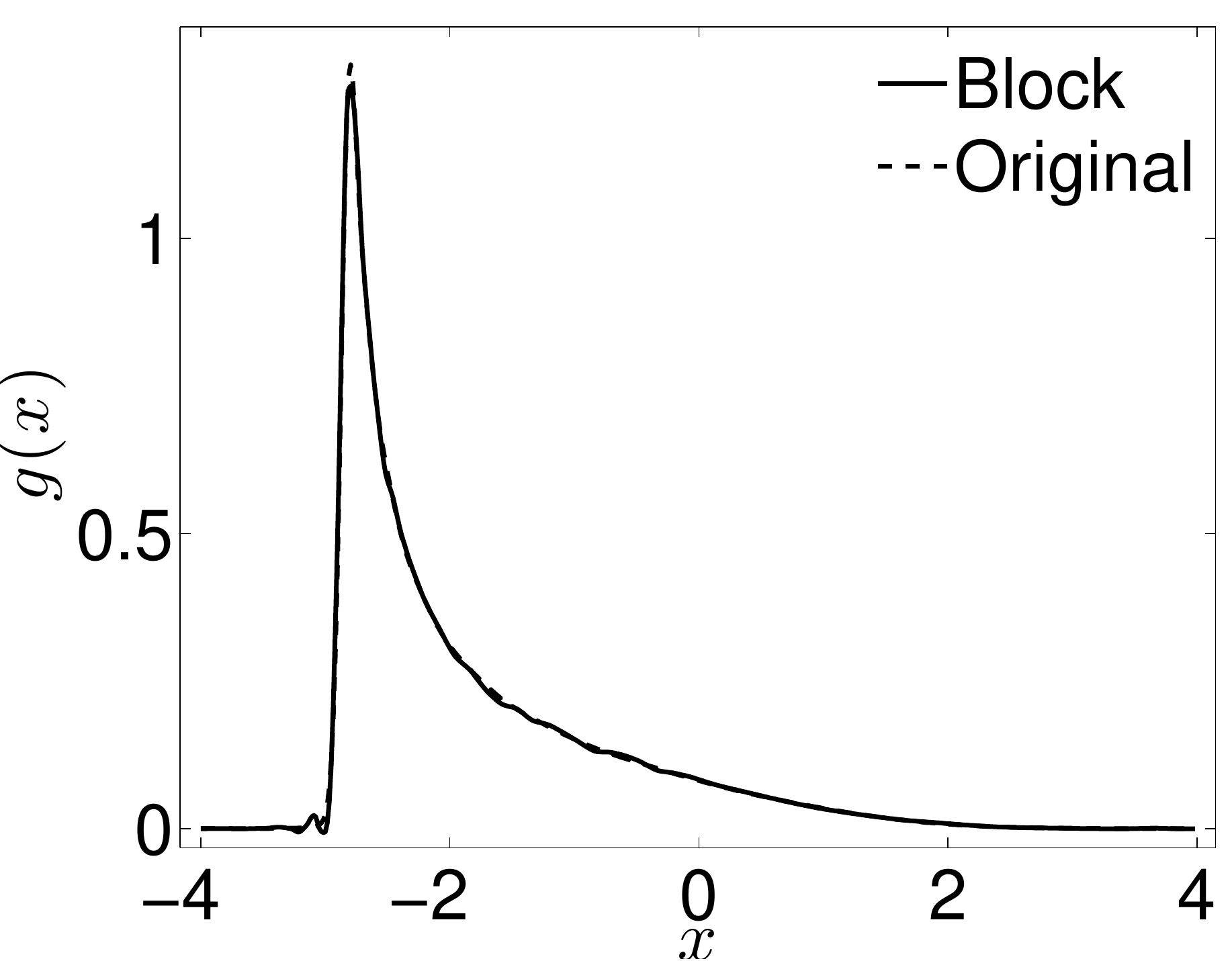}
\subfigure[]{
\includegraphics[width=0.24\textwidth]{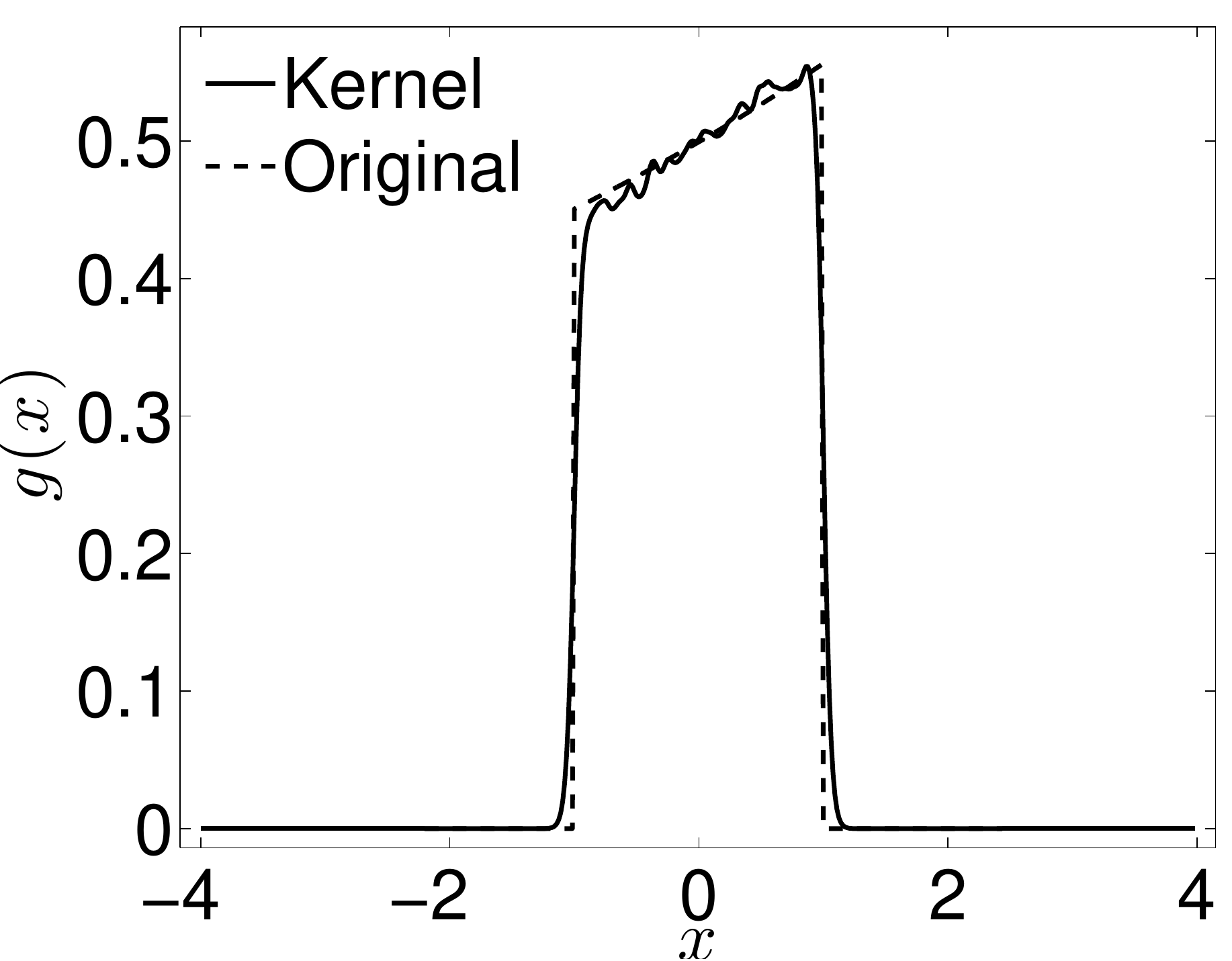}
\includegraphics[width=0.24\textwidth]{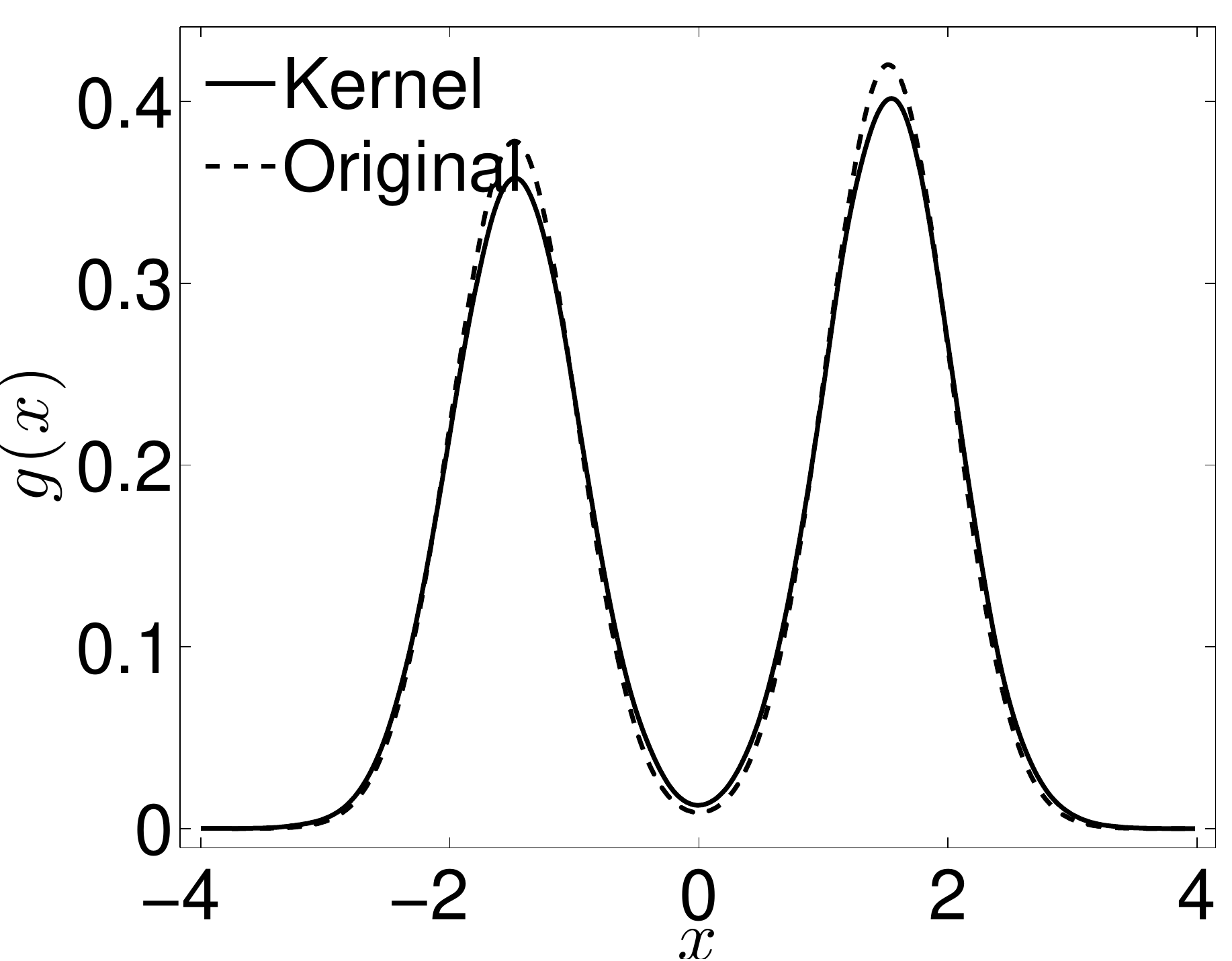}
\includegraphics[width=0.24\textwidth]{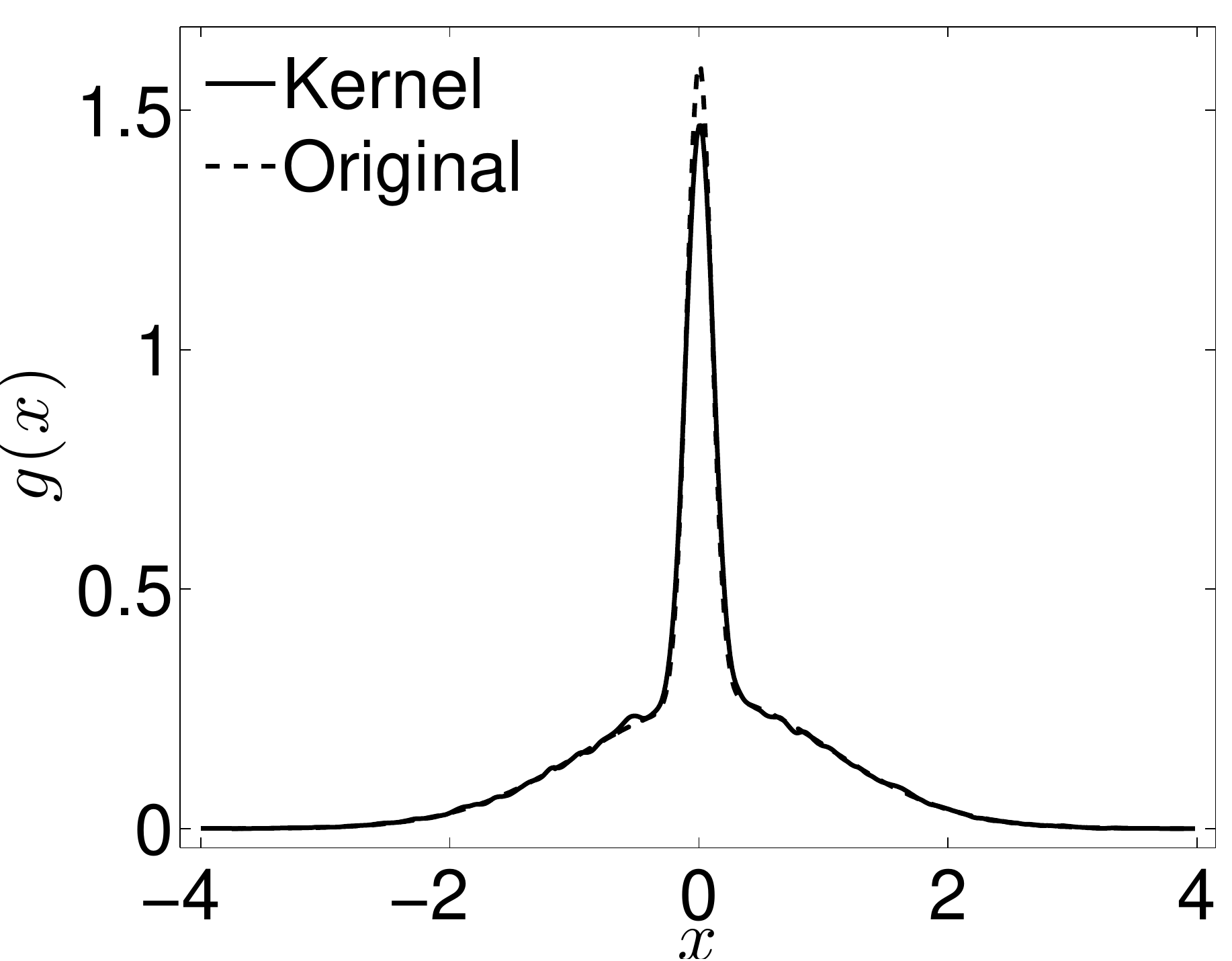}
\includegraphics[width=0.24\textwidth]{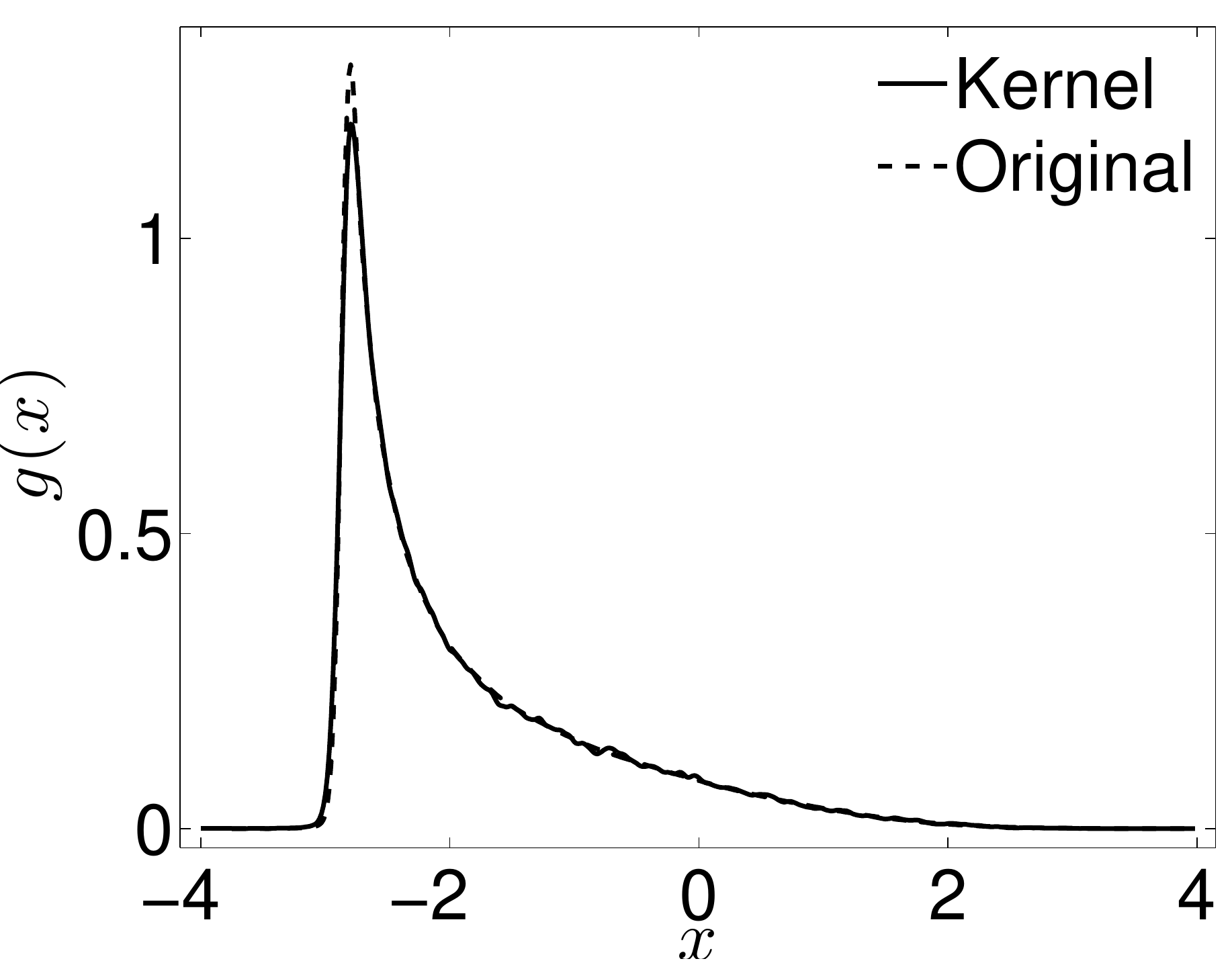}
}
\includegraphics[width=0.24\textwidth]{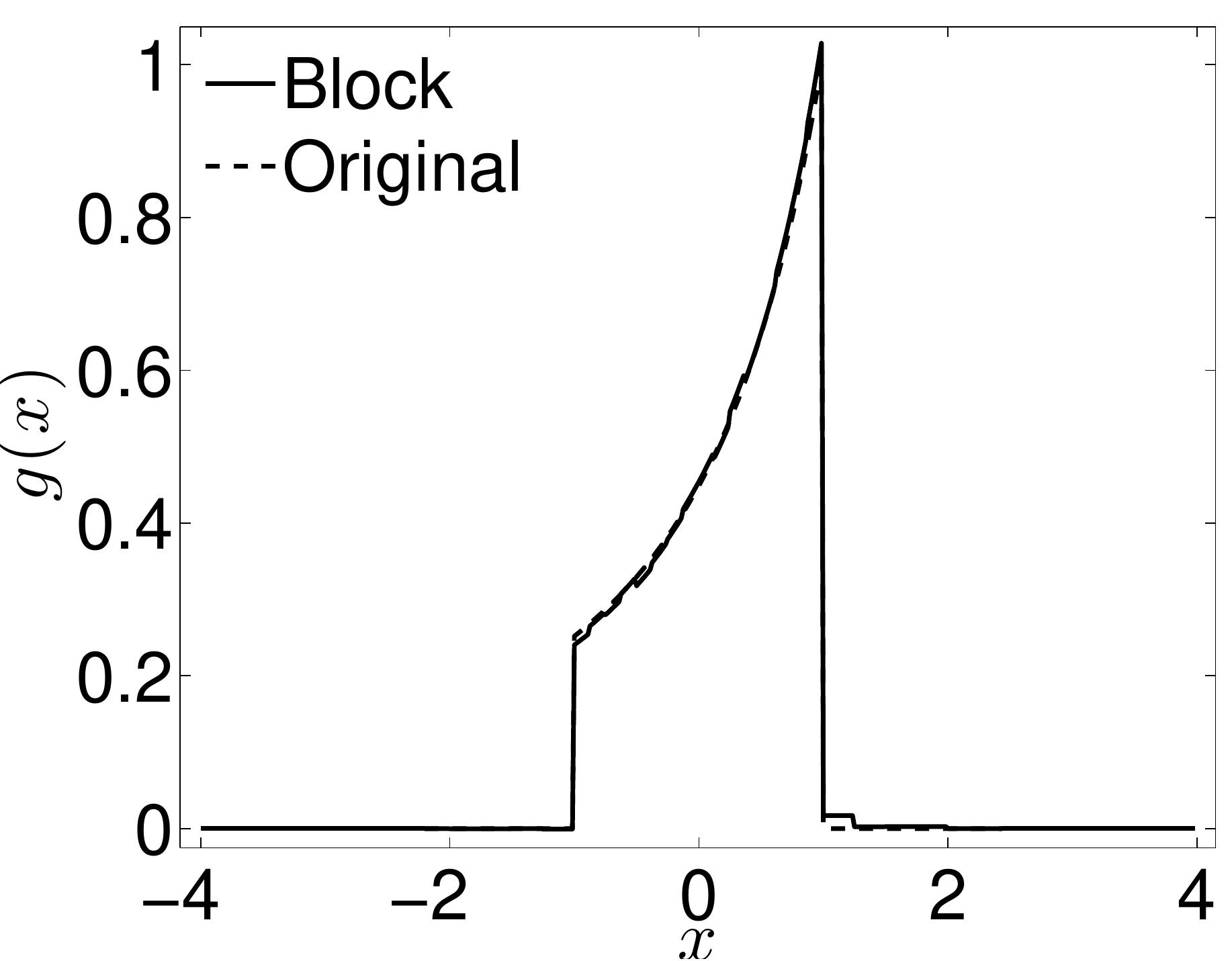}
\includegraphics[width=0.24\textwidth]{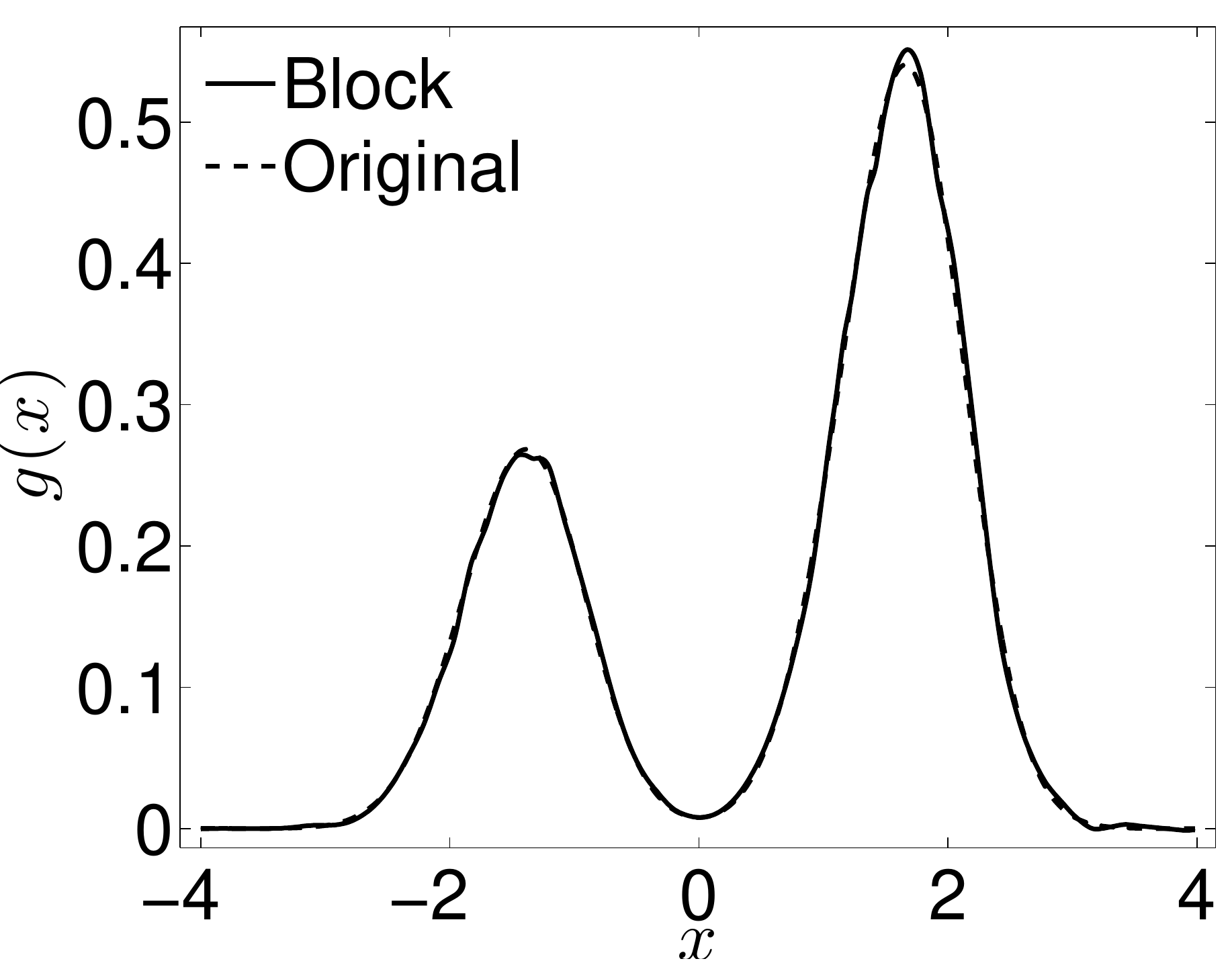}
\includegraphics[width=0.24\textwidth]{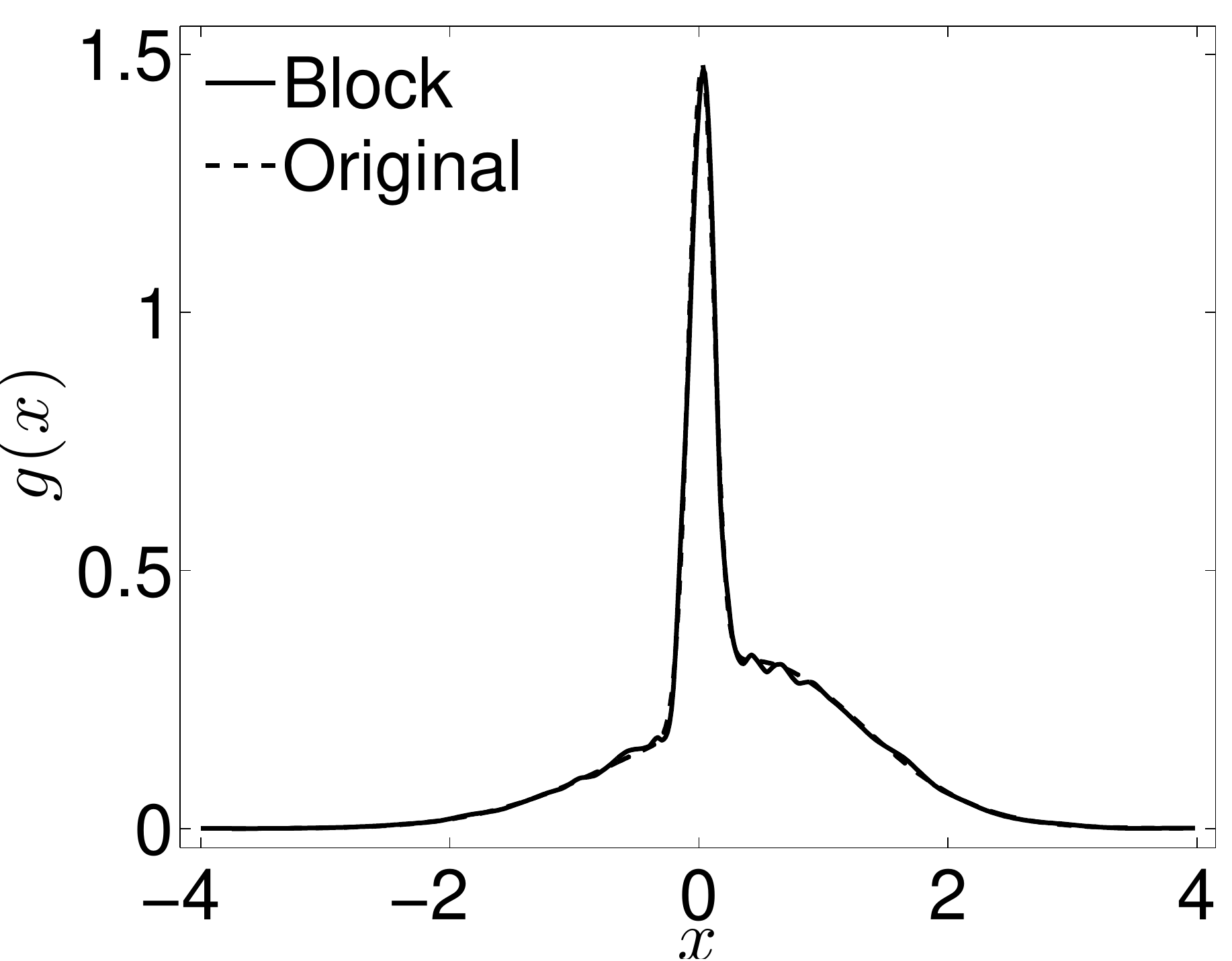}
\includegraphics[width=0.24\textwidth]{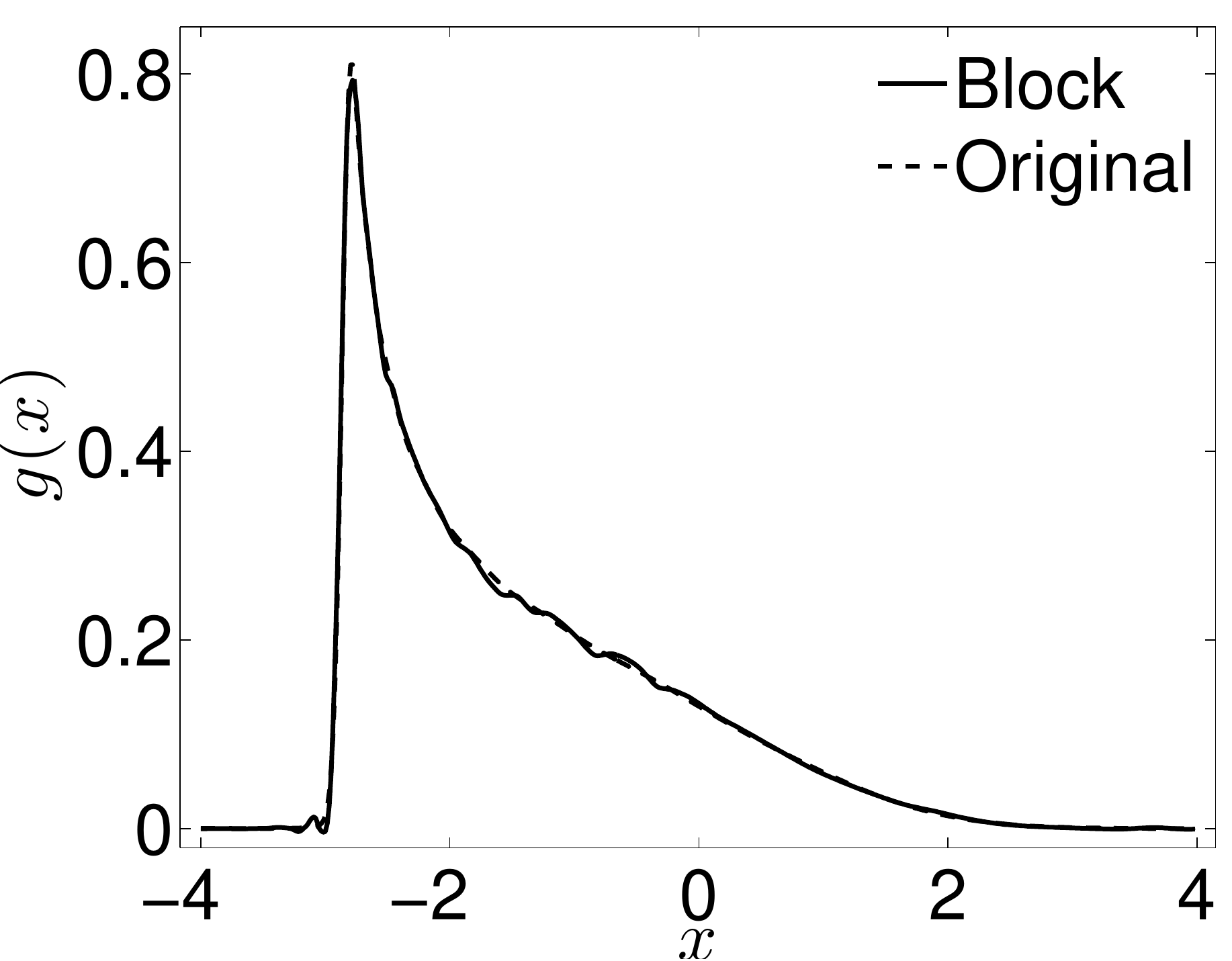}
\subfigure[]{
\includegraphics[width=0.24\textwidth]{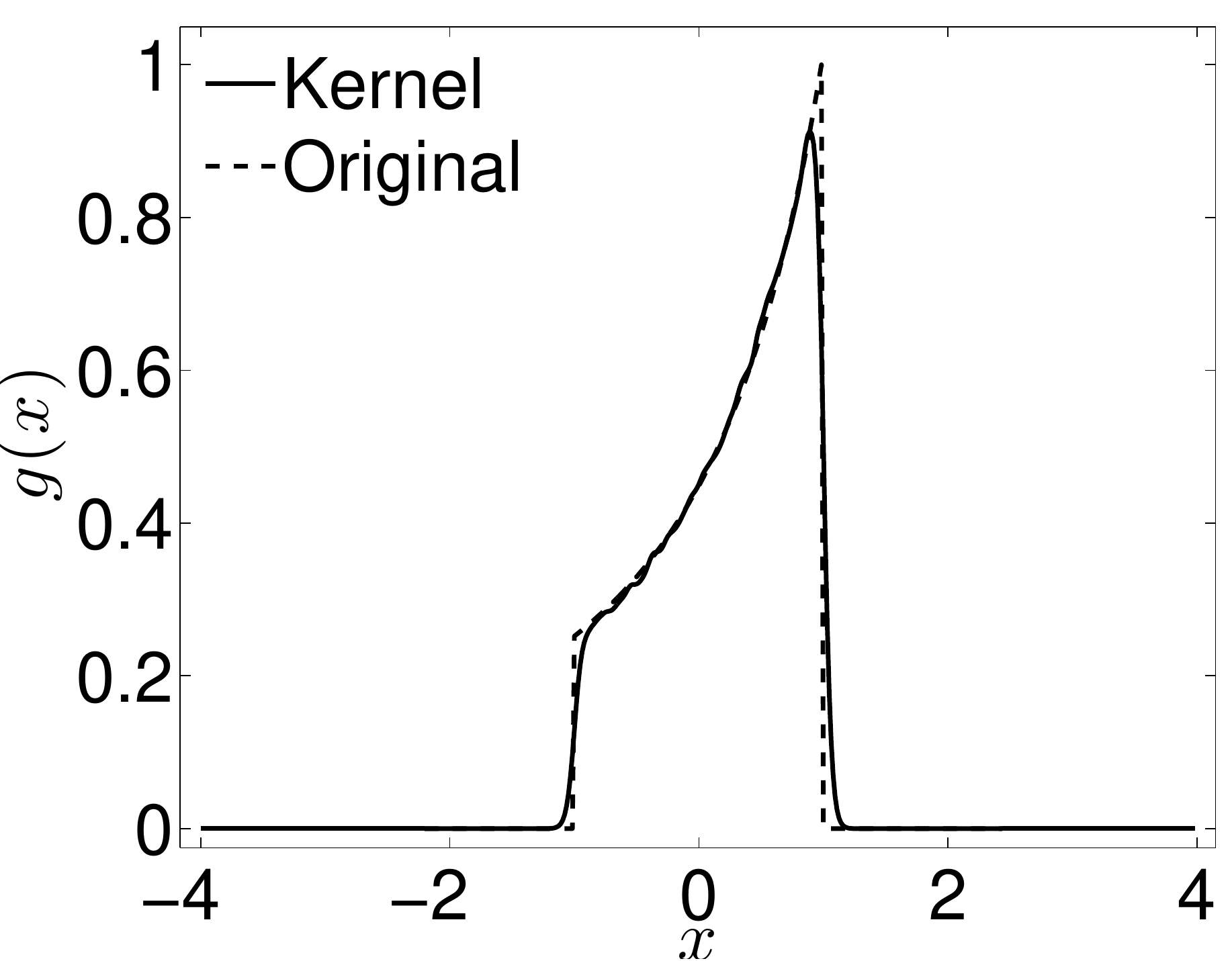}
\includegraphics[width=0.24\textwidth]{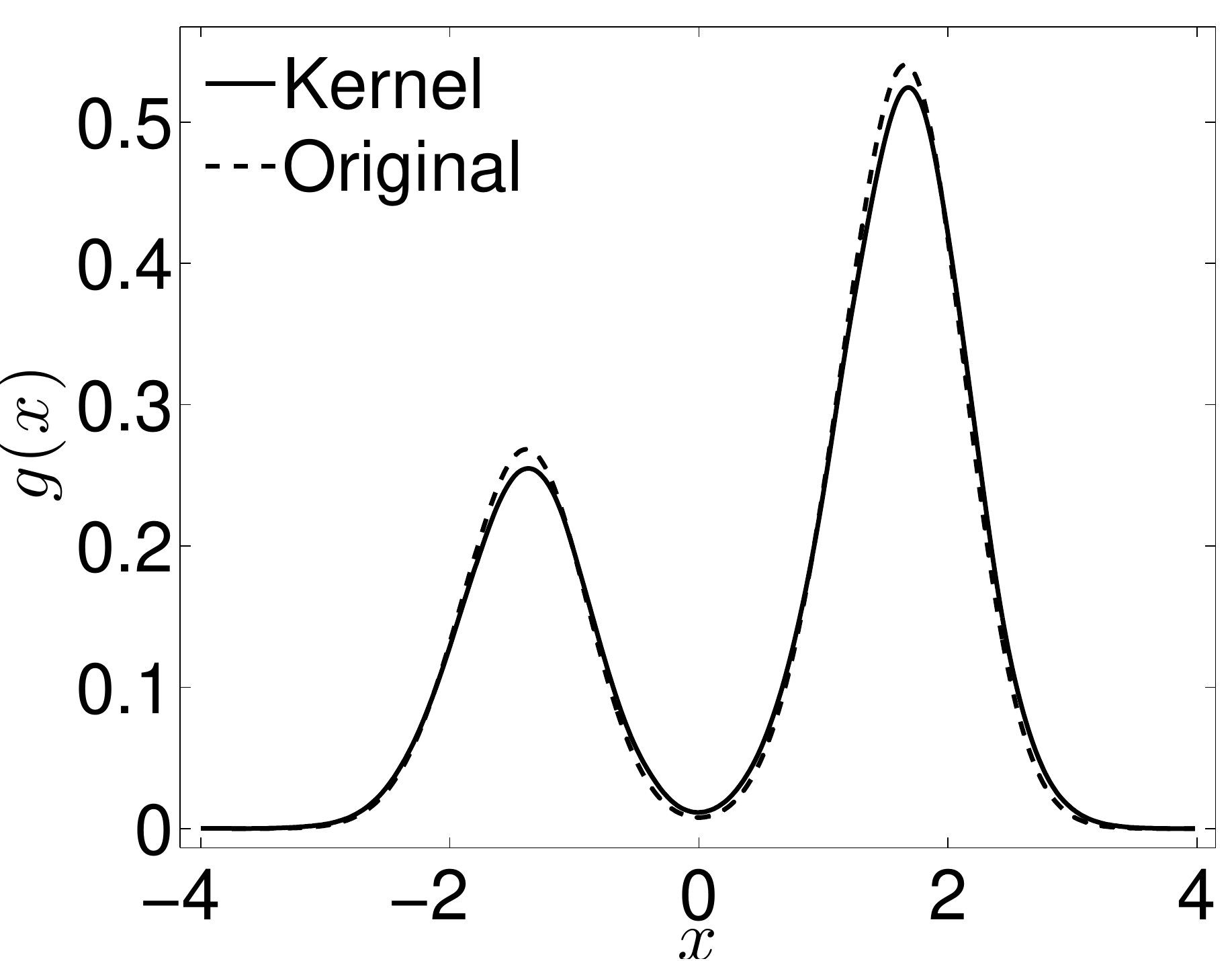}
\includegraphics[width=0.24\textwidth]{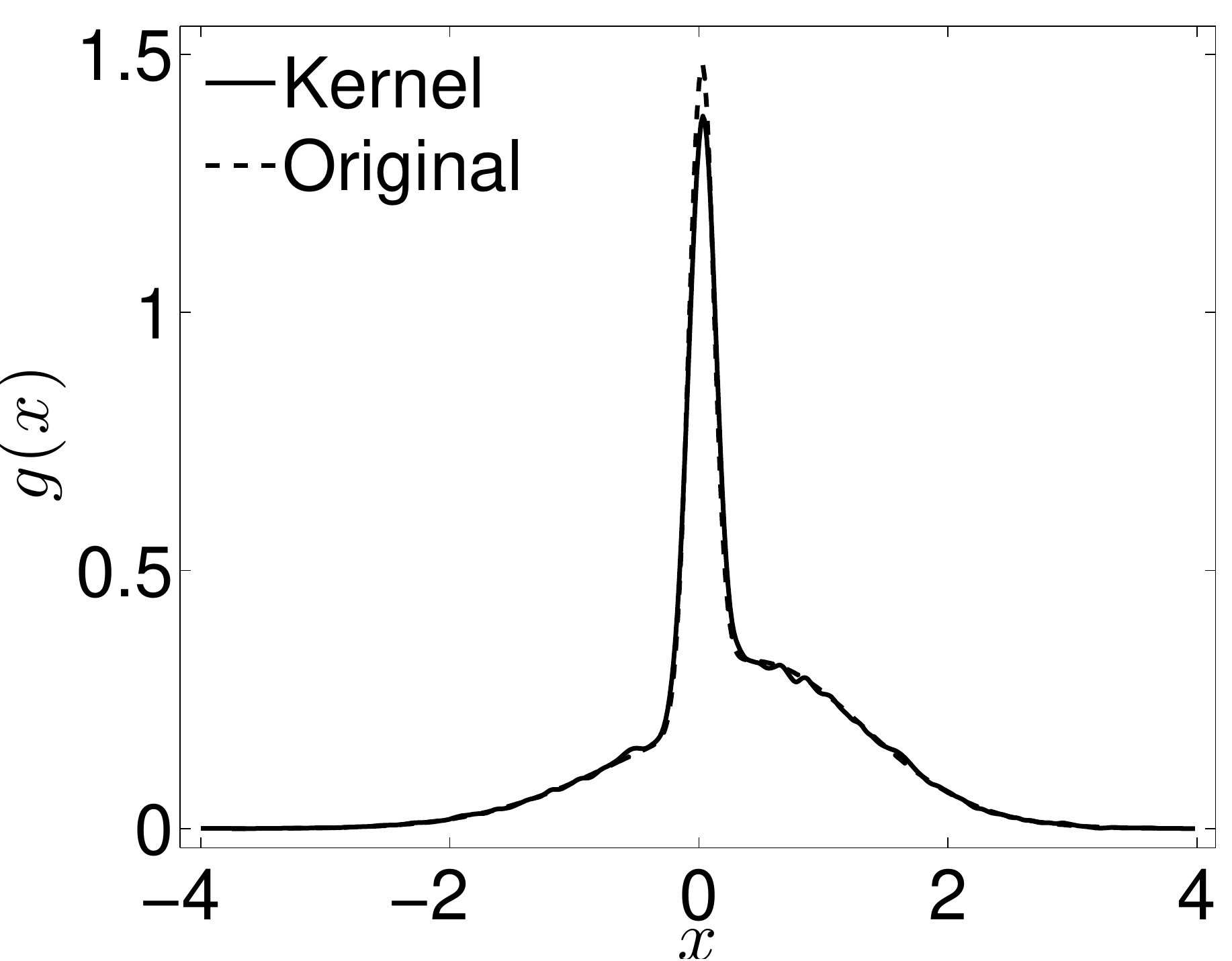}
\includegraphics[width=0.24\textwidth]{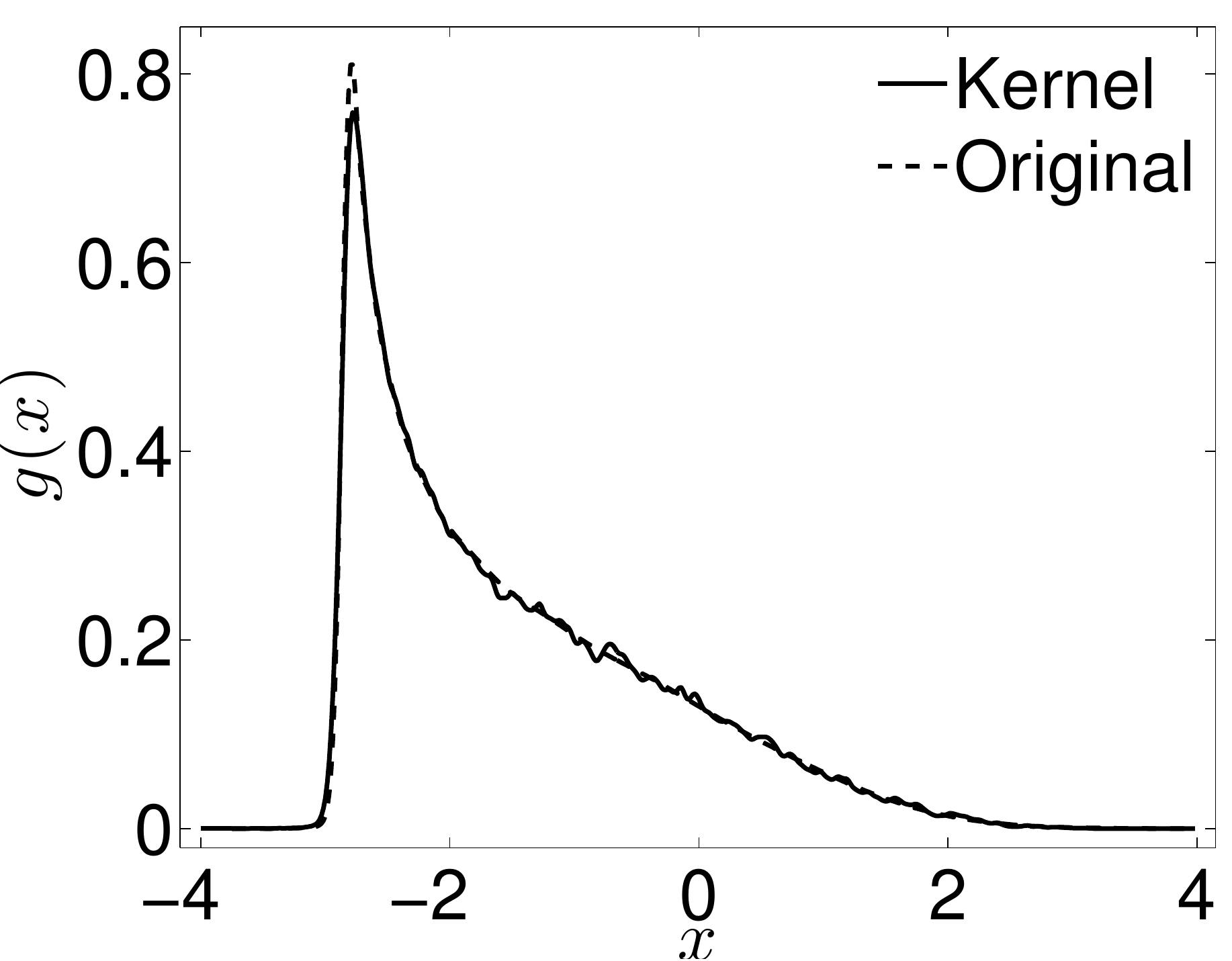}
}
\includegraphics[width=0.24\textwidth]{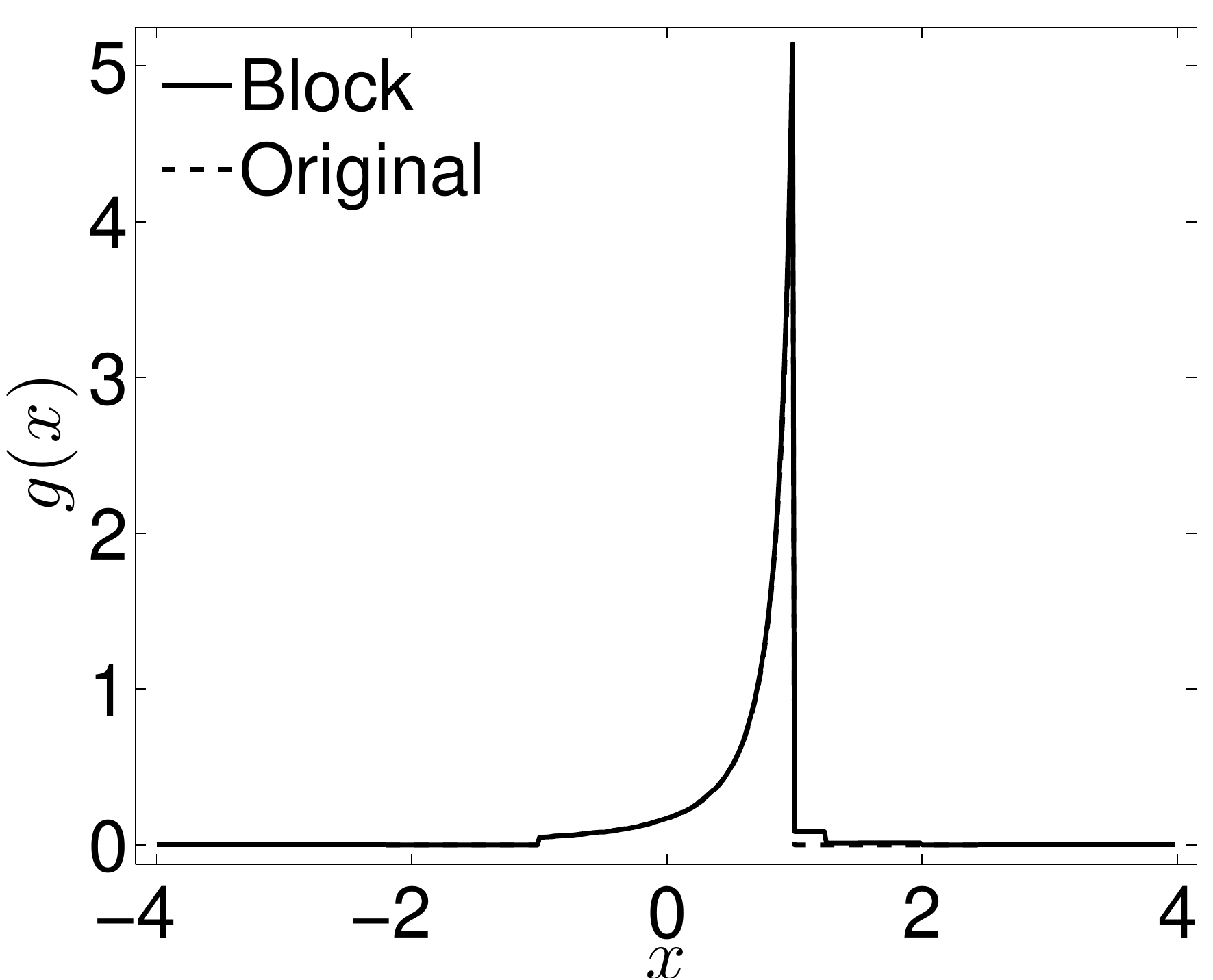}
\includegraphics[width=0.24\textwidth]{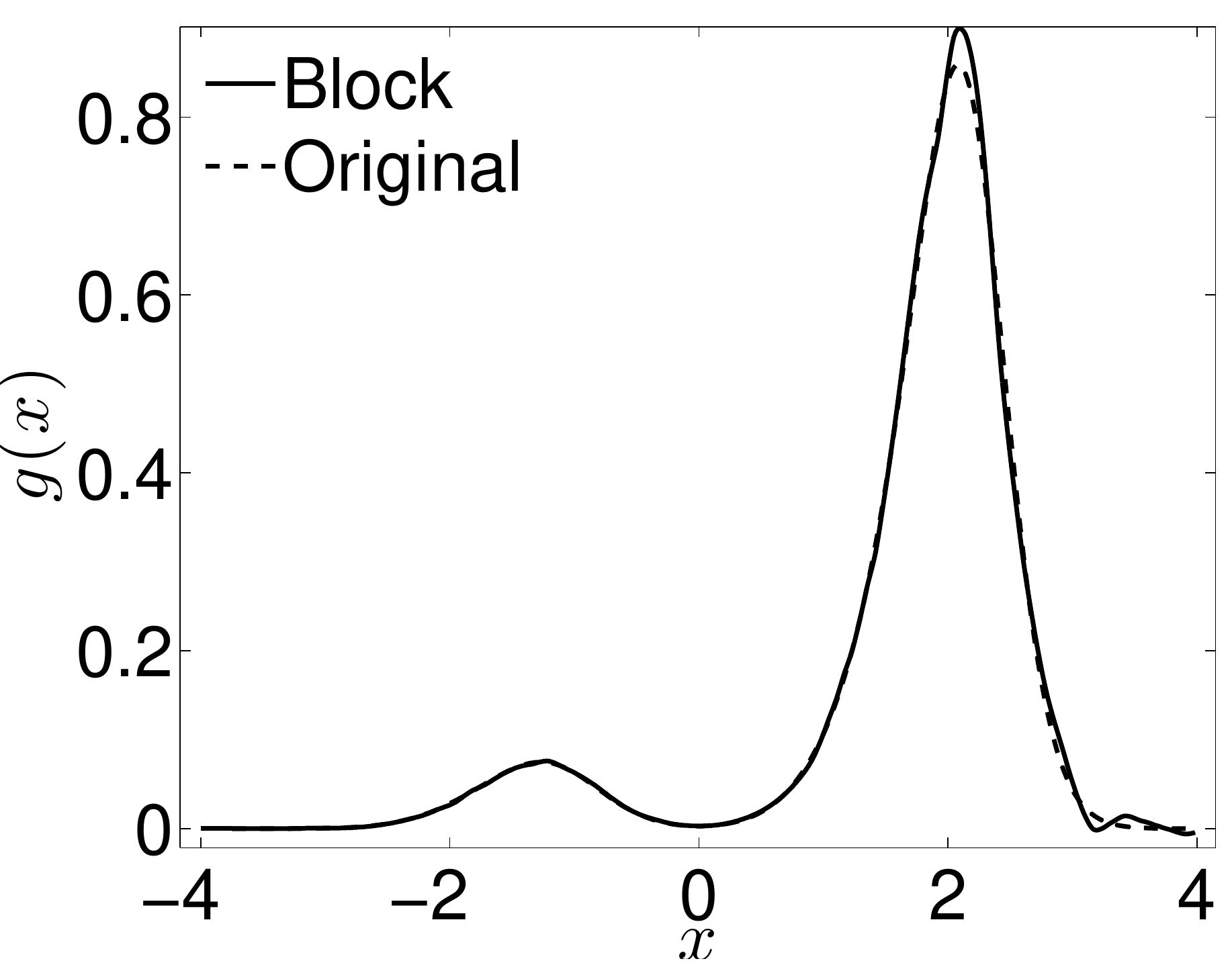}
\includegraphics[width=0.24\textwidth]{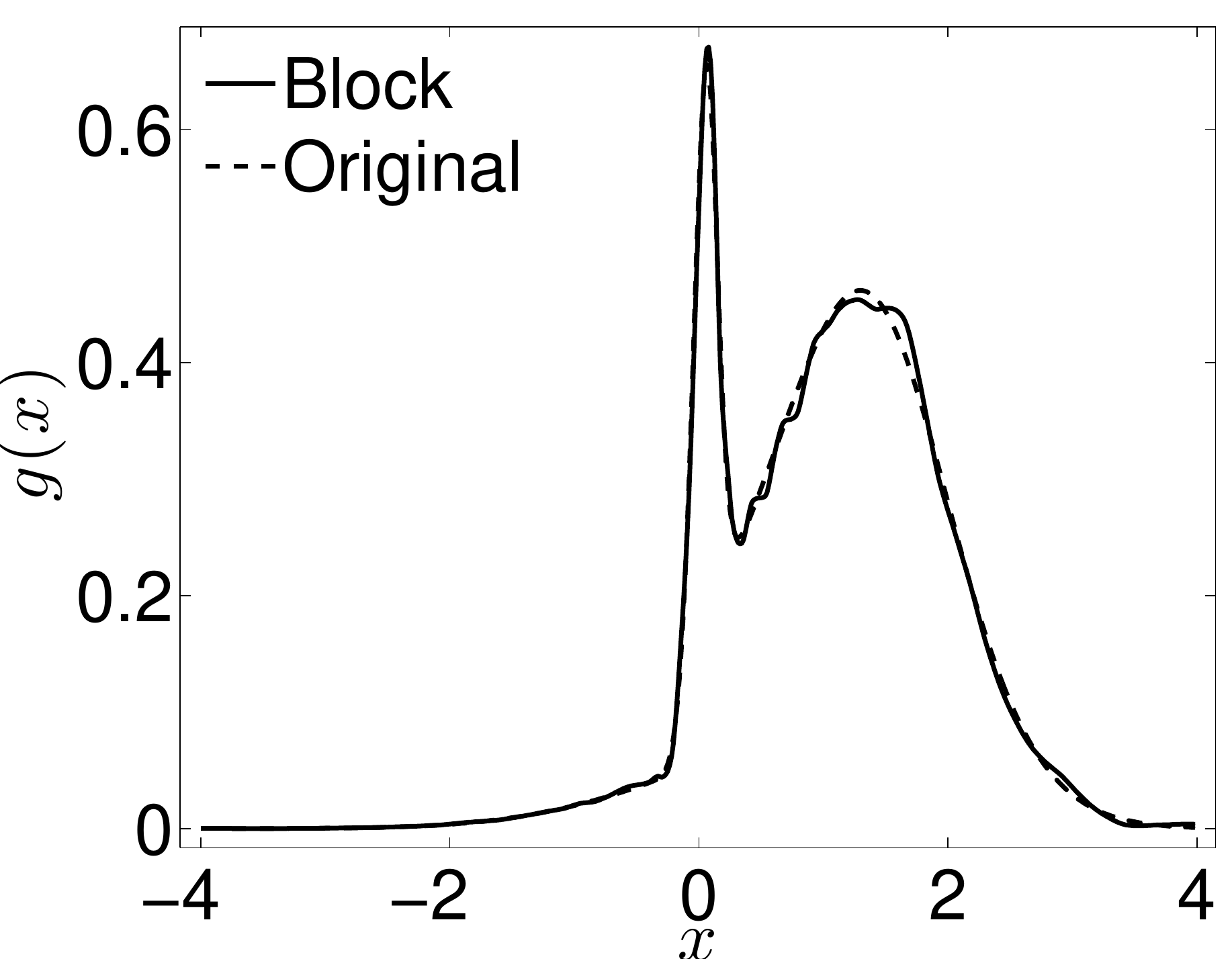}
\includegraphics[width=0.24\textwidth]{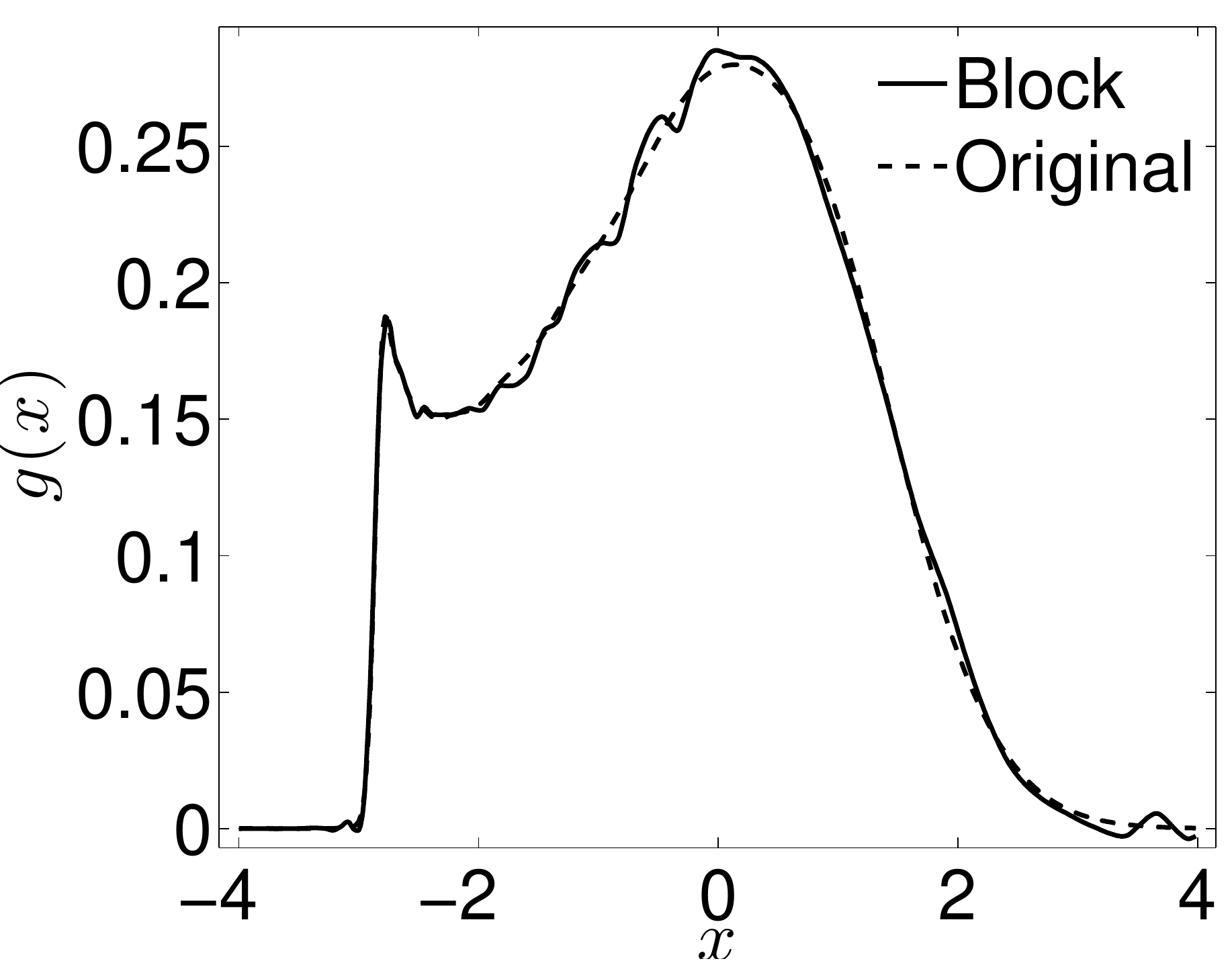}
\subfigure[]{
\includegraphics[width=0.24\textwidth]{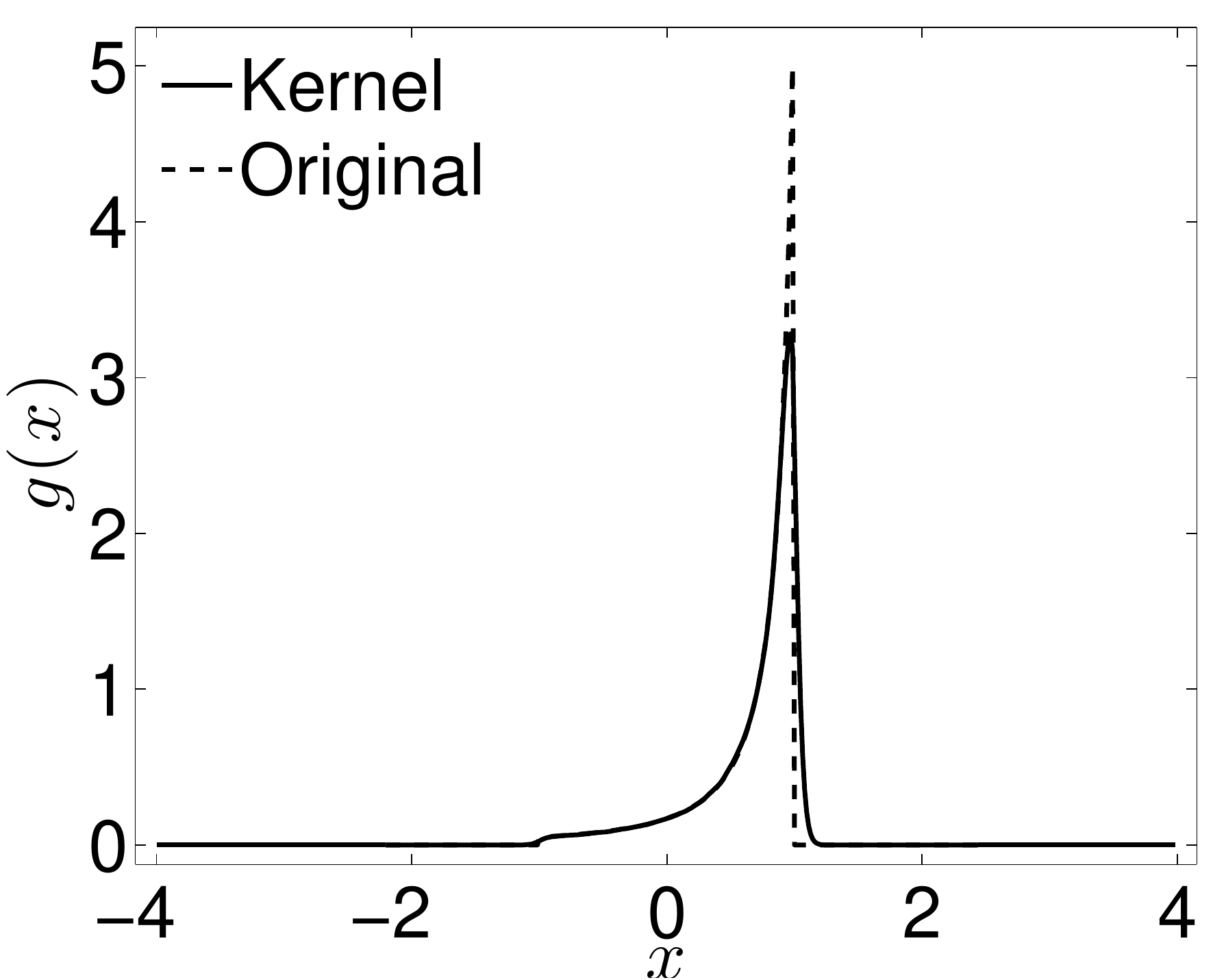}
\includegraphics[width=0.24\textwidth]{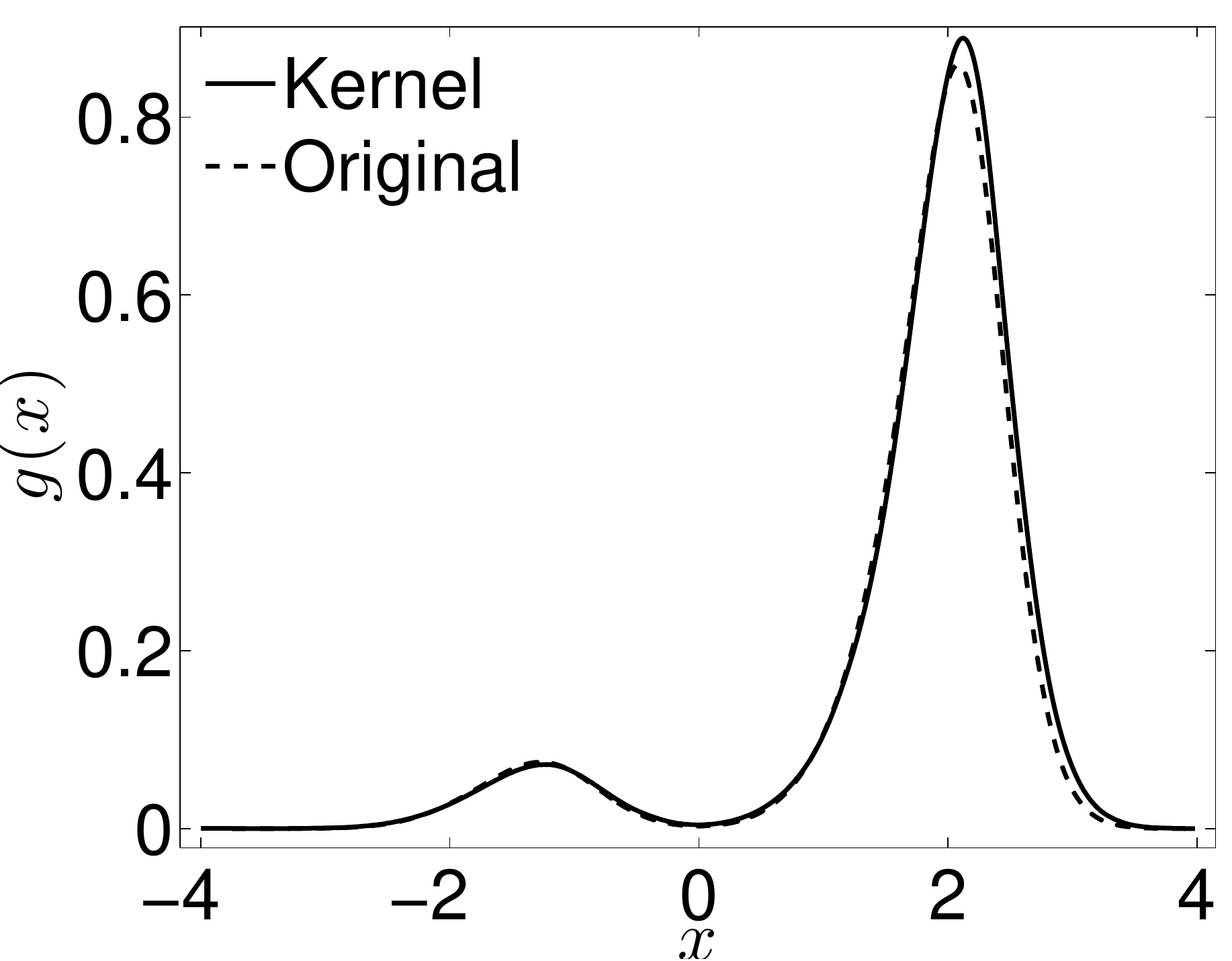}
\includegraphics[width=0.24\textwidth]{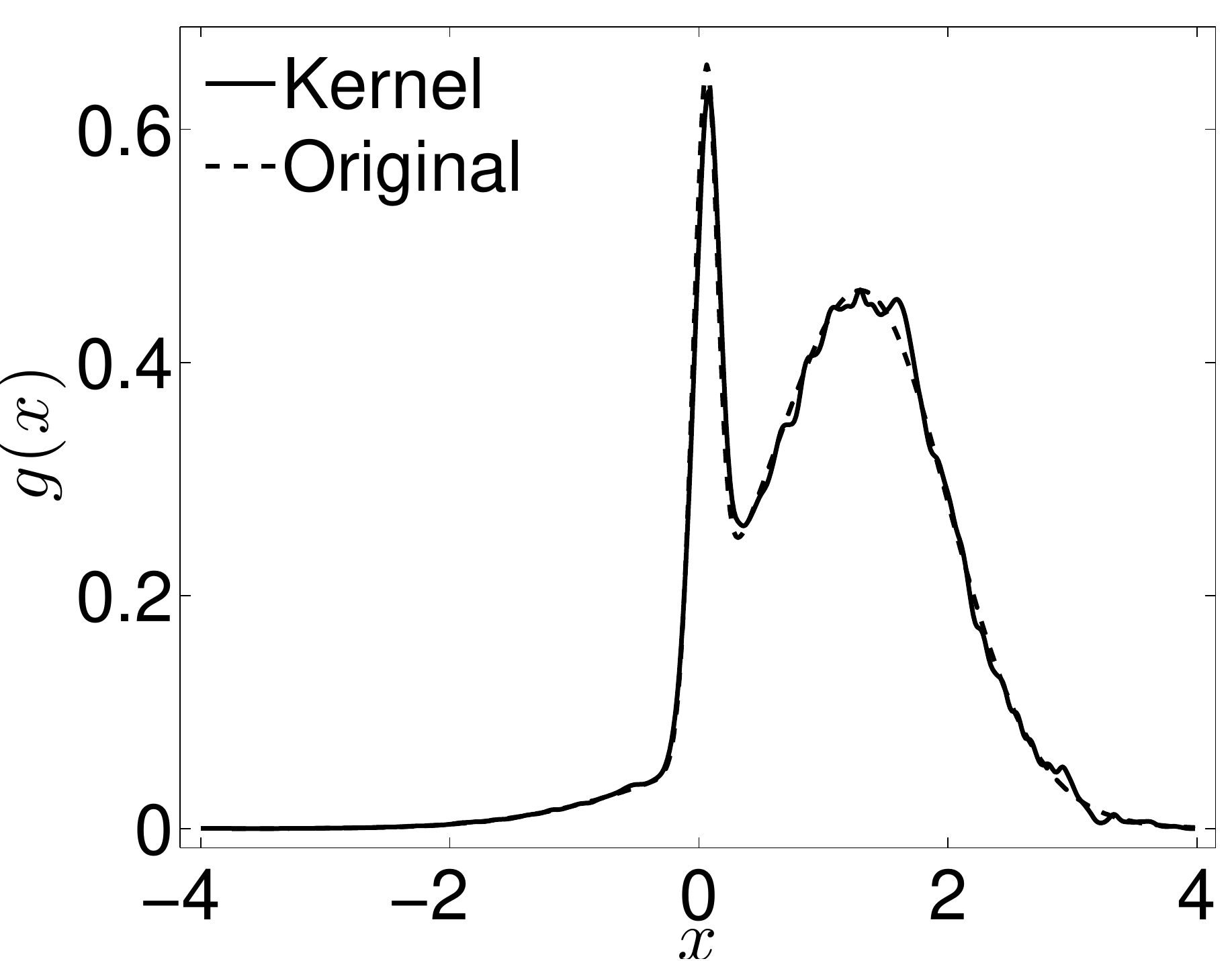}
\includegraphics[width=0.24\textwidth]{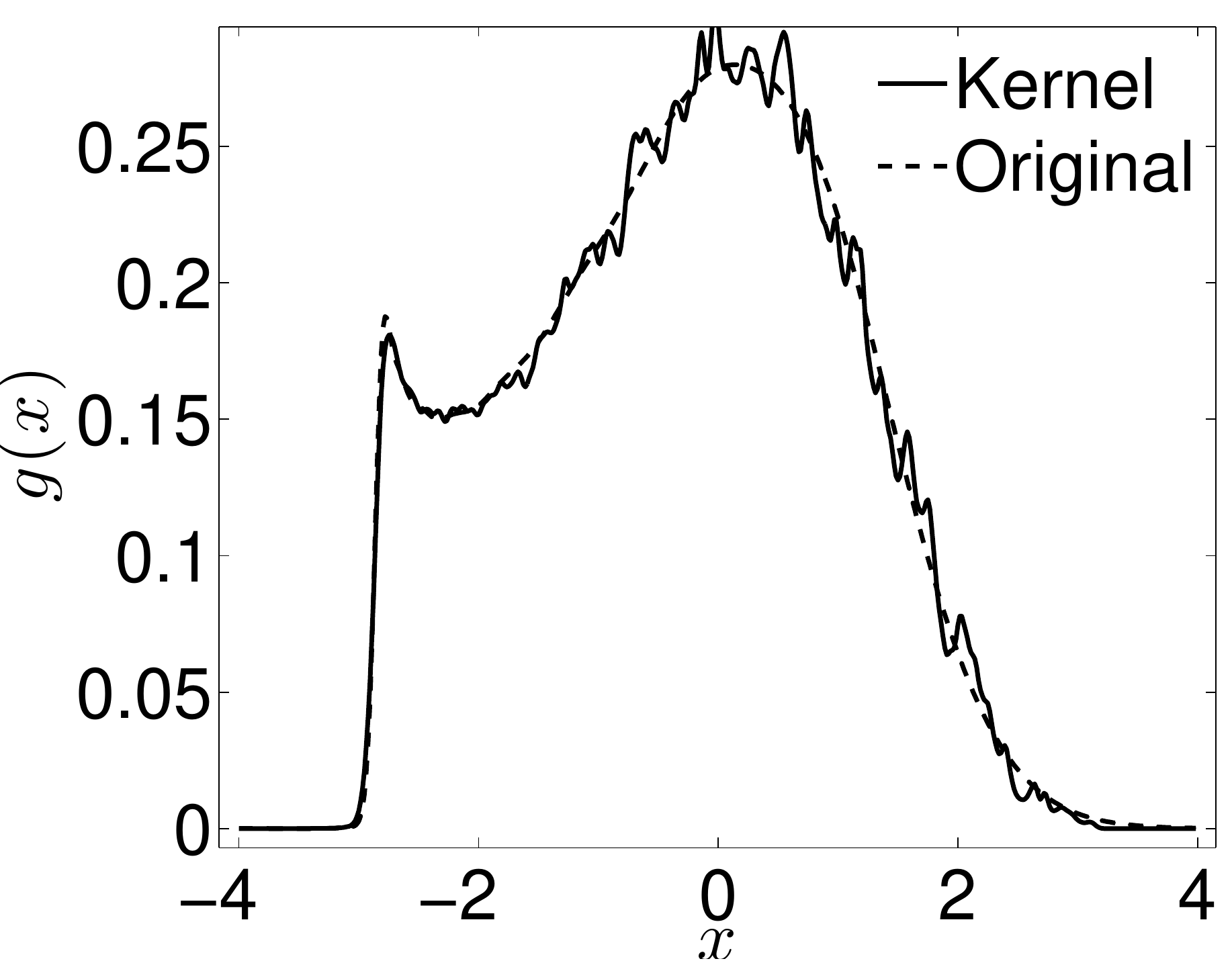}
}
\caption{Original densities (dashed),  Block thresholding estimator $\hat g$ (solid) (1st row), kernel estimator $\hat g_{{\mathrm{LSCV}}}$ (solid) (2nd row) from 50 replications of $n=1000$ samples $X_1,\ldots,X_n$. From left to right Uniform, SeparatedBimodal, Kurtotic and StronglySkewed. $N\sim G(\eta)$, with (a) $\eta=0.9$, (b) $\eta=0.5$ and (c) $\eta=0.1$.} 
\label{fig:MC50}
\end{figure}

We then compared the performance of the Block estimator $\hat g$ with that of the plug-in kernel estimator, say $\hat g_{{\mathrm{LSCV}}}$, given by $\hat g_{{\mathrm{LSCV}}}=w(\hat F(x))\hat f_{\mathrm{LSCV}}(x)$, where  $\hat F$ is defined by \eqref{F} and $\hat f_{\mathrm{LSCV}}$ is given by \eqref{kernel} with $h_{\mathrm{LSCV}}$. Figure~\ref{fig:MC50} shows the results of $\hat g$ and  $\hat g_{\mathrm{LSCV}}$ for $N\sim G(\eta)$, with $\eta=0.9$, $\eta=0.5$ and $\eta=0.1$ respectively. Table \ref{tab:misehi} presents the MISE for samples sizes $n=1000, 2000$ and $5000$. For virtually all cases, the Block estimator consistently showed lower $\field{L}_2$ risk than $\hat g_{{\mathrm{LSCV}}}$, with the exception of the (very smooth) SeparatedBimodal density for which the kernel estimator performs slightly better. This comes at no surprise given that this density is very smooth. Additionally, small discrepancies in the estimate of $f$ may lead to substantial discrepancies for the estimate of $g$ at the locations overweighted by $w(\hat F(\cdot))$. It turns out that this is the case for the Geometric distribution where the weights evolve in $O(1/\eta)$ at high values of $x$, and thus the discrepancies in $\hat{g}$ increase as $\eta$ gets smaller. However, the kernel estimator $\hat g_{{\mathrm{LSCV}}}$ seems to be more concerned (see, Figure~\ref{fig:MC50}(c)), confirming that Block generally provides a better estimate of $f$. Furthermore, as expected, for both methods, and in all cases, the MISE is decreasing as the sample size increases. Without any prior smoothness knowledge on the unknown density, the Block estimator provides very competitive results in comparison to $\hat g_{{\mathrm{LSCV}}}$. 

\begin{table}[htbp]
\caption{$1000\times$MISE values from $50$ replications of sample sizes $n=1000, 2000$ and $5000$, when $N$ follows a Geometric distribution of parameter $\eta$.} %From top to bottom Uniform, SeparatedBimodal, Kurtotic and StronglySkewed.}
\begin{tabular}{lcccccccccccc}\toprule
&\multicolumn{9}{c}{Uniform}\\\hline
 $1.0e$-$03\times$ & \multicolumn{3}{c}{ $\eta=0.9$} &\multicolumn{3}{c}{ $\eta=0.5$} &\multicolumn{3}{c}{ $\eta=0.1$} \\ \hline
  \multicolumn{1}{l}{$n$}          &$~~1000~~$&$~~2000~~$&$~~5000~~$&$~~1000~~$&$~~2000~~$&$~~5000~~$&$~~1000~~$&$~~2000~~$&$~~5000~~$\\
  \hline
  \multicolumn{1}{l}{Block}    &$10.49$&$7.03$&$4.18$&$11.15$&$7.61$&$4.85$&$19.15$&$15.05$&$13.39$\\ 
  \multicolumn{1}{l}{Kernel}  &$13.76$&$11.37$&$9.32$&$17.62$&$14.99$&$12.60$&$70.50$&$63.24$&$55.94$\\\hline\\
 &\multicolumn{9}{c}{SeparatedBimodal}\\\hline
 $1.0e$-$03\times$ & \multicolumn{3}{c}{ $\eta=0.9$} &\multicolumn{3}{c}{ $\eta=0.5$} &\multicolumn{3}{c}{ $\eta=0.1$} \\ \hline
  \multicolumn{1}{l}{$n$}           &$~~1000~~$&$~~2000~~$&$~~5000~~$&$~~1000~~$&$~~2000~~$&$~~5000~~$&$~~1000~~$&$~~2000~~$&$~~5000~~$\\
  \hline
  \multicolumn{1}{l}{Block}    &$8.53$&$6.29$&$3.69$&$9.09$&$6.89$&$3.90$&$15.27$&$11.08$&$7.12$\\ 
  \multicolumn{1}{l}{Kernel}  &$6.30$&$4.98$&$3.54$&$7.00$&$5.32$&$3.82$&$13.26$&$10.86$&$7.79$\\\hline\\
 &\multicolumn{9}{c}{Kurtotic}\\\hline
 $1.0e$-$03\times$ & \multicolumn{3}{c}{ $\eta=0.9$} &\multicolumn{3}{c}{ $\eta=0.5$} &\multicolumn{3}{c}{ $\eta=0.1$} \\ \hline
  \multicolumn{1}{l}{$n$}           &$~~1000~~$&$~~2000~~$&$~~5000~~$&$~~1000~~$&$~~2000~~$&$~~5000~~$&$~~1000~~$&$~~2000~~$&$~~5000~~$\\
  \hline
  \multicolumn{1}{l}{Block}    &$11.31$&$8.03$&$5.52$&$11.93$&$8.04$&$5.57$&$18.08$&$9.22$&$7.01$\\ 
  \multicolumn{1}{l}{Kernel}  &$12.27$&$8.00$&$5.83$&$13.03$&$8.44$&$6.13$&$21.97$&$13.04$&$9.66$\\\hline  \\
 &\multicolumn{9}{c}{StronglySkewed}\\\hline
 $1.0e$-$03\times$ & \multicolumn{3}{c}{ $\eta=0.9$} &\multicolumn{3}{c}{ $\eta=0.5$} &\multicolumn{3}{c}{ $\eta=0.1$} \\ \hline
  \multicolumn{1}{l}{$n$}           &$~~1000~~$&$~~2000~~$&$~~5000~~$&$~~1000~~$&$~~2000~~$&$~~5000~~$&$~~1000~~$&$~~2000~~$&$~~5000~~$\\
  \hline
  \multicolumn{1}{l}{Block}    &$9.91$&$7.78$&$5.12$&$8.69$&$7.02$&$4.70$&$10.12$&$9.03$&$5.93$\\ 
  \multicolumn{1}{l}{Kernel}  &$10.57$&$8.13$&$5.97$&$11.16$&$8.68$&$6.24$&$19.62$&$15.54$&$10.86$\\\hline 
  \botrule
\end{tabular} 
 \label{tab:misehi}
\end{table}

\end{document}